\documentclass[11pt,a4paper]{article}
\usepackage[T1]{fontenc}
\usepackage[utf8]{inputenc}
\usepackage{lmodern}
\usepackage[bidi=default]{babel}
\babelprovide[main,import]{english}

\def\languageshorthands#1{}
\usepackage{amsmath,amssymb,amsthm,mathrsfs}
\usepackage[margin=2.5cm]{geometry}
\usepackage{longtable,array}
\usepackage{enumitem}
\usepackage[hidelinks]{hyperref}

\newtheorem{theorem}{Theorem}[section]
\newtheorem{proposition}[theorem]{Proposition}
\newtheorem{lemma}[theorem]{Lemma}
\newtheorem{corollary}[theorem]{Corollary}
\theoremstyle{definition}
\newtheorem{definition}[theorem]{Definition}
\newtheorem{example}[theorem]{Example}
\theoremstyle{remark}
\newtheorem{remark}[theorem]{Remark}

\numberwithin{equation}{section}
\newcolumntype{L}[1]{>{\raggedright\arraybackslash}p{#1}}

\newcommand{\R}{\mathbb R}
\newcommand{\C}{\mathbb C}
\newcommand{\Z}{\mathbb Z}
\newcommand{\dd}{\,\mathrm d}
\DeclareMathOperator{\Rea}{Re}
\renewcommand{\Cap}{\operatorname{Cap}}
\newcommand{\vol}{\operatorname{vol}}

\title{Green profiles and positive harmonic functions near isolated subanalytic singularities}
\author{Flavio Mosca}
\date{November 27, 2026}

\begin{document}
\maketitle

\begin{abstract}
Let $X\subset\mathbb{R}^N$ be subanalytic, closed near $0$, of pure dimension $d\ge2$, analytic away from $0$, Whitney regular at $0$ and with connected link, endowed with the induced metric. In a previous paper it was proved that, under a uniform Harnack condition on the links along a sequence of radii, positive harmonic functions near $0$ can be written as $qG_X+u_{\mathrm b}$, with $q\ge0$, where $G_X$ is a singular profile, unique up to a positive factor and the addition of a bounded harmonic function, and $u_{\mathrm b}$ is bounded and harmonic. Here we study the size of the profile. Under a uniform annular Harnack condition, $G_X$ is comparable on sufficiently small links to the reciprocal of the capacity of the condenser of the links. We deduce its growth order for generalized horns and for germs intrinsically bi-Lipschitz equivalent to them. For exact horns with exponent $\gamma>1$ we determine the leading coefficient in terms of the flux and, for $d\ge3$, the first correction term, which depends on the embedding of the cross-section. We study germs whose link factors collapse with different exponents and logarithmic terms of the profile. Finally, we construct a semialgebraic germ of dimension $3$, with link diffeomorphic to the sphere, on which a bounded harmonic function has no limit at the singular point, the optimal annular Harnack constants diverge as the radius tends to zero, and the cone of positive harmonic functions vanishing on a sufficiently small outer link has at least two extreme rays.
\end{abstract}



\tableofcontents

\section{Introduction}\label{sec:intro}

Let $X\subset\R^N$ be a subanalytic set, closed near $0$, of pure dimension $d\ge2$, with an analytic regular part away from $0$, satisfying the Whitney condition at $0$ and having connected link. In \cite[Theorem~5.17]{V11} it was proved that, under a uniform Harnack condition on the links along a sequence of radii, every bounded harmonic function near $0$ extends continuously to $0$ and every positive harmonic function near $0$ admits the asymptotic decomposition $u=qG_X+u_{\mathrm b}$. The proof gives no information about the growth order of $G_X$; in \cite{V11} this order is determined for cones, germs with regular spherical blow-up, and rotational cusps.

This paper addresses the following questions; the related known results and elementary consequences are collected in Appendix~\ref{app:noti}.
\begin{enumerate}[label=(Q\arabic*)]
\item Growth order of $G_X$ in terms of the geometry of the germ, in particular for generalized horns, for germs bi-Lipschitz equivalent to them with respect to the intrinsic metrics, and for elliptic cusps (\S\ref{sec:capacita-profilo}, \S\ref{sec:corni}).
\item Leading coefficient of $G_X$ for exact horns and its normalization by flux (\S\ref{sec:corni-esatti}).
\item Subsequent terms in the expansion of $G_X$, their dependence on the embedding of the cross-section, germs with different collapse exponents, and logarithmic terms (\S\ref{sec:correzione}, \S\ref{sec:collasso}, \S\ref{sec:risonanze}).
\item Whether, for connected links, continuous extension of bounded harmonic functions and the decomposition (AD) can fail, and how many minimal positive harmonic functions can exist over the vertex (\S\ref{sec:collo}, \S\ref{sec:martin}).
\item Behavior of harmonic functions on the complex singularities $A_{\mathsf k}$ with the induced metric (\S\ref{sec:Ak} and \S A.3).
\end{enumerate}

The main results are as follows.
\begin{enumerate}[label=(\alph*)]
\item If $X$ satisfies condition (H$_b$), the generator $G_R$ of the cone $\mathcal K_R$, vanishing on the link $L_R$, is comparable on the link $L_r$ to $1/\Cap_X(r;R)$, where $\Cap_X(r;R)$ is the capacity of the condenser between the links $L_r$ and $L_R$ (Theorem~\ref{thm:capacita-profilo}).
\item For a generalized horn $X_{\gamma,L}$ of dimension $d$, with $\beta:=\gamma(d-1)$, we have $G_X\asymp\rho^{1-\beta}$ if $\beta>1$ and $G_X\asymp\log(1/\rho)$ if $\beta=1$ (a cone in dimension $2$); the same holds for every germ bi-Lipschitz equivalent to $X_{\gamma,L}$ with respect to the intrinsic metrics (Theorem~\ref{thm:ordine-corni} and Corollary~\ref{cor:bilip}). In particular, for a germ bi-Lipschitz equivalent to a cone over a smooth connected link of dimension $d-1\ge2$, we have $G_X\asymp\rho^{2-d}$ (Corollary~\ref{cor:coni}).
\item For an exact horn with $\gamma>1$ there exists $c>0$ such that $G_X=c\,\rho^{1-\beta}(1+O(\rho^{\eta}))$ uniformly on the link, for every $0<\eta<\min\{2\gamma-2,\beta-1\}$ (Theorem~\ref{thm:coefficiente}); the flux of $G_X$ is $c(\beta-1)\vol_h(L)$ (Corollary~\ref{cor:flusso} and Remark~\ref{oss:flusso-generale}).
\item For exact horns with $\gamma>1$ and $d\ge3$, the next term of the expansion at fixed height $z$ is $c\,z^{1-\beta+2\gamma-2}f_1$, with $f_1$ explicitly determined by the embedding of the cross-section (Theorem~\ref{thm:correzione}); $f_1$ is not a bi-Lipschitz invariant (Proposition~\ref{prop:non-invarianza}).
\item For germs whose link factors collapse with different exponents, we obtain the decomposition of positive harmonic functions with a remainder, the Picard principle, and the exact profile (Theorem~\ref{thm:multipla} and Proposition~\ref{prop:capacita-multipla}); two germs that are not bi-Lipschitz equivalent can have the same leading-order term of the profile (Proposition~\ref{prop:x24-x33}), and the profile of rotational cusps contains a logarithmic term for infinitely many exponents (Theorem~\ref{thm:risonanza}).
\item For each integer $k\ge1$ we construct a semialgebraic germ $X_k$ of dimension $3$ with link diffeomorphic to $S^2$, on which a bounded harmonic function has no limit at $0$ (Theorem~\ref{thm:discontinuita}). For each $b\ge2$ the optimal annular Harnack constant is bounded below by $c_b r^{-k/2}$, with $c_b>0$, for sufficiently small $r$ (Corollary~\ref{cor:harnack-collo}); the Picard principle, condition (AD), and the Bôcher property fail, and there are at least two minimal positive harmonic functions over the vertex (Corollary~\ref{cor:martin-collo}).
\end{enumerate}

\subsection{Context and perspectives}\label{sec:contesto}
For each section we indicate how its results can be interpreted within current research and what their limitations are. The conjectures and questions formulated in this subsection are not used in the sequel.

\subsubsection*{Section~\ref{sec:prelim}: preliminaries (capacity, flux, bi-Lipschitz transport)}
\emph{Interpretation: the vertex as a potential-theoretic end.} The capacity of the condenser between the links and the flux at the vertex treat the singular point as an end of a manifold in the sense of Grigor'yan, that is, as an object measured by capacity \cite{Gri}. Lemma~\ref{lemma:trasporto} is the technical key of the whole paper: capacity is preserved up to constants by homeomorphisms that are bi-Lipschitz with respect to the intrinsic distance, and all subsequent invariance results rest on this fact.

\emph{Limitations.} These are preliminary tools without independent significance; they make possible the connection with Lipschitz geometry used in Sections~\ref{sec:capacita-profilo}, \ref{sec:corni}, \ref{sec:correzione} and~\ref{sec:collasso}.

\subsubsection*{Section~\ref{sec:capacita-profilo}: capacity and profile}
\emph{Interpretation 1: an analytic problem reduced to a geometric one.} Theorem~\ref{thm:capacita-profilo} states that, under (H$_b$), the profile $G_R$ on the link $L_r$ is comparable to the inverse of the capacity between $L_r$ and $L_R$. In itself this is the adaptation to the vertex of the classical comparison between Green function and capacity on manifolds and metric measure spaces (Remark~\ref{oss:classico}); the novelty lies in the setting. Capacity depends only on the metric and the Riemannian measure, and the growth of the profile becomes computable with the tools of the Lipschitz geometry of singularities.

\emph{Interpretation 2: a prediction for normal complex surfaces.} A normal complex surface decomposes into a thick part, which is essentially metrically conical, and a thin part, which shrinks rapidly towards the origin \cite[Theorems~1.6 and~1.7]{BNP}. The thick part suggests a capacitary comparison with a cone, and one may expect $G_X\asymp\rho^{-2}$ for every normal complex surface with nonempty thick part, provided (HS) holds. This would be a uniform statement on a whole classical class of singularities. A two-sided comparison for the full germ, and its passage from the Harnack sequence to arbitrary radii, are not established here.

\emph{Limitations.} The comparison on all sufficiently small links is proved here under (H$_b$); under (HS) alone it is obtained on the links of the Harnack sequence (Remark~\ref{oss:HS}). For complex surfaces this hypothesis is not verified, and Section~\ref{sec:collo} shows that it may fail in general; the prediction therefore depends on such a verification.

\subsubsection*{Section~\ref{sec:corni}: order of the profile for horns and bi-Lipschitz germs}
\emph{Interpretation 1: a dictionary between Lipschitz invariants and potential-theoretic exponents.} For a generalized horn $G_X\asymp\rho^{1-\gamma(d-1)}$ when $\gamma(d-1)>1$, with $G_X\asymp\log(1/\rho)$ in the exceptional case $d=2$, $\gamma=1$, and the same holds for every germ bi-Lipschitz equivalent to the horn (Corollary~\ref{cor:bilip}). In the inner bi-Lipschitz classification of semialgebraic surfaces with an isolated singularity and connected link the models are the $\beta$-horns, and two horns are equivalent if and only if they have the same exponent \cite{Bir} (see also \cite[Theorem~8.3]{Bir08}). Section~\ref{sec:corni} translates this exponent into a potential-theoretic exponent: the growth exponent of the profile is a numerical invariant of the germ, bi-Lipschitz by construction.

\emph{Interpretation 2: robustness of cones.} Corollary~\ref{cor:coni} gives $G_X\asymp\rho^{2-d}$ for $d\ge3$, and $G_X\asymp\log(1/\rho)$ for $d=2$, for every germ bi-Lipschitz equivalent to a cone over a smooth connected link. The Euclidean growth of the Green function is thus preserved under metric deformations of bounded distortion; this is the potential-theoretic counterpart of the fact that metrically conical germs form the simplest case of the classification.

\emph{Limitations.} Invariance rests on (H$_b$), known for horns from \cite{V11}. It remains open for which classes of subanalytic germs the exponent of the profile is determined by inner invariants alone; without further hypotheses a single singular profile need not describe all positive harmonic functions, as Sections~\ref{sec:collo} and~\ref{sec:martin} show.

\subsubsection*{Section~\ref{sec:corni-esatti}: exact horns, leading coefficient and flux}
\emph{Interpretation 1: a canonical Green function with pole at the vertex.} For exact horns with $\gamma>1$, the profile has the expansion $c\,\rho^{1-\beta}(1+O(\rho^\eta))$, for $0<\eta<\min\{2\gamma-2,\beta-1\}$, and normalization to unit flux fixes the leading coefficient; with zero values on a prescribed outer boundary it determines the Green profile of that domain, a Green function with pole at the singular point, an object that is not defined in general because the pole coincides with the singularity.

\emph{Interpretation 2: an invariant of volume type.} Theorem~\ref{thm:capacita-esatta} gives $\Cap^\Pi(\varepsilon;\zeta)\sim(\beta-1)V_L\,\varepsilon^{\beta-1}$: to first order potential theory sees the cross-section only through its volume, in heuristic analogy with the leading coefficients of Weyl-type laws. The conceptual value appears in comparison with Section~\ref{sec:correzione}: the first order records the volume, the second the embedding.

\emph{Interpretation 3: exactly solvable models.} Proposition~\ref{prop:sfera}, for cross-sections contained in a sphere, gives explicit profiles, which serve as test cases for a general theory.

\emph{Limitations.} The asymptotic formula and the flux normalization concern only exact horns with $\gamma>1$, with the parametrization $\Pi(z,y)=(z^\gamma y,z)$; for an arbitrary profile the coefficient $c$ depends on its normalization, whereas after normalization to unit flux it equals $1/((\beta-1)V_L)$ (Corollary~\ref{cor:flusso}).

\subsubsection*{Section~\ref{sec:correzione}: first angular correction}
\emph{Interpretation 1: potential theory sees more than the inner Lipschitz geometry.} The leading order is a bi-Lipschitz invariant (Corollary~\ref{cor:bilip}); the term $z^\kappa f_1$ is not, and depends on the embedding of the cross-section (Proposition~\ref{prop:non-invarianza}). For normal complex surfaces the outer Lipschitz geometry determines several analytic invariants \cite[Theorem~1.2]{NP}. The expansion of the Green profile lies between the two geometries: the leading growth order is an inner bi-Lipschitz invariant, whereas the first correction, at fixed height $z$ in the given parametrization, is sensitive to the embedding. This leads to an inverse problem: to determine which invariants of the embedding can be recovered from the full expansion of $G_X$.

\emph{Interpretation 2: explicit local formulas.} Theorem~\ref{thm:correzione} expresses $f_1$ through $S-\langle S\rangle$, with $S=\frac12|y|^2$, and through the mean of $T_L=|\nabla_hS|_h^2$, the square of the tangential component of the position vector: these are second-order moments of the position of the cross-section, analogous to an inertia tensor. One expects the subsequent terms to involve moments of higher order; the expansion would then become a sequence of computable local extrinsic invariants.

\emph{Limitations.} Exact horns, $\gamma>1$, $d\ge3$. Taken alone, the result is of moderate scope; its interest depends on the inverse problem it raises.

\subsubsection*{Section~\ref{sec:collasso}: multi-rate collapsing horns}
\emph{Interpretation 1: the profile records several contraction rates.} The germs $X_{2,4}$ and $X_{3,3}$ have the same leading exponent of the profile and of the volume, but their links contain cycles that shrink at different rates (Proposition~\ref{prop:x24-x33}). In the classification of \cite[Theorem~1.9]{BNP} the thin part is refined into pieces distinguished by rates and directions of contraction; Proposition~\ref{prop:x24-x33} shows that, in these two examples, the next terms of the profile distinguish germs whose leading profile and volume exponents agree. This is a second entry, after Section~\ref{sec:corni}, in the dictionary between Lipschitz invariants and the expansion of the profile.

\emph{Interpretation 2: a class in which the theory is complete.} For germs with product structure Theorem~\ref{thm:multipla} gives a mode decomposition, by separation of variables, Lemma~\ref{lemma:edo} and the explicit profile $\varsigma_{\boldsymbol\gamma}$; the decomposition of positive harmonic functions, the Picard principle and the exact profile follow. This is analogous, at a heuristic level, to the analysis of warped products in singular calculi \cite[Chapter~2]{Mel}, and it provides a class in which questions (Q1)--(Q3) are answered for the exact product models considered here.

\emph{Limitations.} The product structure is a strong hypothesis; it remains to be seen how stable the results are under perturbations that do not preserve it.

\subsubsection*{Section~\ref{sec:risonanze}: logarithmic resonances}
\emph{Interpretation 1: polyhomogeneous structure with arithmetic indices.} Theorem~\ref{thm:risonanza} states that the term $\log(1/\rho)$ appears if and only if $\mathsf n(2\gamma-2)=\beta-1$. This is the analogue of the coincidence of indicial exponents in singular calculi, where solutions of the form $r^\gamma(\log r)^N\phi$ arise \cite[Proposition~5.61]{Mel}. The explicit radial integral is already given in \cite[Prop.~5.12]{V11}, and logarithmic terms for regular spherical blow-up appear in \cite[Example~5.11]{V11}; the contribution here is the resonance criterion and its coefficient for collapsing rotational cusps; the coefficient $-81/128$ of Example~8.2 is not visible on the tangent cone.

\emph{Interpretation 2: o-minimal definability of the Green profile.} The Green profile of a semialgebraic germ need not be globally subanalytic near the vertex. In the globally subanalytic setting, finite parametric integrals of subanalytic functions are log-analytic \cite[Théorème d'intégration and proof of Théorème~1]{LR}, and constructible functions, sums of products of globally subanalytic functions and their logarithms, form a class that is stable under integration \cite[Theorem~1.3]{CM}. We conjecture that the Green profiles of subanalytic germs with an isolated singularity have asymptotic expansions in this class; Theorem~\ref{thm:risonanza} provides an explicit example in the collapsing case. A positive answer would connect potential theory on singular spaces with o-minimality.

\emph{Limitations.} The result concerns the collapsing class, in which the profile is an explicit integral: it is a guiding example, not yet a structural theorem.

\subsubsection*{Section~\ref{sec:collo}: a semialgebraic germ with connected link without continuous extension}
\emph{Interpretation 1: the obstruction to extension is metric, not topological.} Theorem~\ref{thm:discontinuita} answers a natural question in the negative: with connected link, even in the semialgebraic category, a bounded harmonic function may have no limit at the vertex. On normal complex analytic germs with an isolated singularity, of complex dimension at least two, with connected regular part and $H^1(Y_{\rm reg};\R)=0$, real-valued pluriharmonic functions extend (Theorem~\ref{th:pluri}); this is a different class of functions. The example exhibits a metric obstruction, produced by the thin neck, which is not detected by the connectedness of the link.

\emph{Interpretation 2: condition (HS) of \cite{V11} is substantial.} Corollary~\ref{cor:harnack-collo} gives $C_{H,b}(r)\ge c_b\,r^{-k/2}$: the optimal Harnack constants grow at least as fast as $r^{-k/2}$, and hypothesis (HS) in the theorems of \cite{V11} is not an artifact of the proof. The spectral gap of the cross-sections (Theorem~\ref{thm:gap}) is the spectral signature of the neck, analogous to the small eigenvalues of dumbbell manifolds.

\emph{Limitations.} The result concerns a single family of models, in dimension $3$. Proposition~\ref{prop:conduttanza-collo} and Theorem~\ref{thm:gap} do not enter the proof of Theorem~\ref{thm:discontinuita}; they acquire weight in the programme of Section~\ref{sec:martin}.

\subsubsection*{Section~\ref{sec:martin}: minimal harmonic functions over the vertex and Martin boundary}
\emph{Interpretation 1: the vertex is not a point for potential theory.} For Euclidean domains the number of Martin points over a boundary point depends on the local geometry and need not be one. At the vertex of a cone generated by a spherical John region the Martin point is unique and minimal \cite[Theorem~1.1]{Hirata}. A John domain has finitely many minimal Martin points at each boundary point, at most two if the John constant exceeds $\sqrt3/2$, in the convention of \cite[Theorem~1.1]{AHL}. Section~\ref{sec:martin} carries the phenomenon over to a setting without boundary, with intrinsic metric and in the semialgebraic category: a connected link does not guarantee a unique Martin point (Corollary~\ref{cor:martin-collo}). In potential-theoretic terms, the same geometric vertex supports at least two distinct minimal Martin kernels.

\emph{Interpretation 2: a thick--thin decomposition for potential theory.} The decomposition of \cite{BNP} is analogous to the Margulis decomposition of negatively curved manifolds. The germ $X_k$ consists of two thick lobes joined by a thin neck. We conjecture that the minimal Martin points over the vertex correspond to the components of the thick part separated by necks whose conductance is negligible with respect to the radial weight; in this programme Proposition~\ref{prop:conduttanza-collo} and Theorem~\ref{thm:gap} provide the quantities that fix the threshold.

\emph{Limitations.} It is not proved that there are exactly two minimal Martin points, and a general criterion is lacking.

\subsubsection*{Section~\ref{sec:Ak}: the singularities \texorpdfstring{$A_{\mathsf k}$}{Ak}}
\emph{Interpretation: a first quantitative step on complex singularities with the induced metric.} The conductance of the balanced band of $A_{\mathsf k}$, for $\mathsf k>1$, is of order $\varrho/((\mathsf k-1)\log(1/\varrho))$ (Proposition~11.1), and in the reduced model of Remark~11.2 a coupling of this order prevents distinct limits in the two channels. The heuristic is that on $A_{\mathsf k}$ bounded harmonic functions extend, with a Harnack constant of logarithmic growth: a possible test case between the two regimes of Section~\ref{sec:martin}.

\emph{Limitations.} The results are partial and programmatic; the passage from the model to the PDE is not proved. An answer to question (Q5) would provide a test case for such a criterion on complex surface singularities.

\subsubsection*{Appendix~\ref{app:noti}: known or elementary results}
The statements of Appendix~\ref{app:noti} are not new in method; their interest lies in delimiting the framework of Sections~\ref{sec:capacita-profilo}--\ref{sec:Ak} by explicit examples.

\paragraph{\S\S A.1--A.2: uniform charts, thickness condition, variable number of charts.}
\emph{Interpretation 1: a thickness criterion complementary to the neck germ.} The chain Definition~\ref{def:T}, Proposition~\ref{prop:T-U}, Lemma~\ref{lemma:U-H}, Theorem~\ref{th:thick} translates a metric condition into (H), hence into continuous extension, the Picard principle, (AD) and the Bôcher property. The metric condition requires every point of the link to be interior to a graph of full dimension, at a scale $s_T(r)$ satisfying the upper bound and the link-volume condition of Definition~\ref{def:T}. This is the thick side of the picture whose thin side is Section~\ref{sec:collo}; the criterion provides the qualitative extension and asymptotic-decomposition properties stated in Theorem~\ref{th:thick}.

\emph{Interpretation 2: a series criterion.} Corollary~A.6 allows a number of charts that grows with the scale, provided $\sum_jC_M^{-n_{\rm cart}(r_j)}=+\infty$. The series criterion is \cite[Theorem~3.3]{V11}, and the contribution here is a geometric sufficient condition; it is a divergence criterion of the same type as Wiener's criteria: extension depends on the accumulation of local effects, not on a uniform condition.

\emph{Limitations.} The graph decomposition of \cite{Stasica} exists for every bounded subanalytic set, but (T) does not follow from it. The Moser constants are not optimal. Between the criterion of Corollary~A.6 and the lower bound $r^{-k/2}$ of Corollary~\ref{cor:harnack-collo} there remains an unquantified gap.

\paragraph{\S A.3: pluriharmonic functions on normal complex singularities.}
\emph{Interpretation: what the complex structure guarantees.} On a normal complex germ with an isolated singularity, of complex dimension at least two, with connected regular part and $H^1(Y_{\rm reg};\R)=0$, every real pluriharmonic function is the real part of a holomorphic function on the germ, hence Lipschitz at the vertex (Theorem~\ref{th:pluri}). This holds for $A_{\mathsf k}$, whose link is a lens space (Corollary~\ref{cor:pluri-AE}). Question (Q5) measures exactly the gap between this class and the much larger class of functions harmonic for the induced metric.

\emph{Limitations.} The vanishing of $H^1(Y_{\rm reg};\R)$ is used to obtain a single-valued holomorphic primitive; the result gives no information on harmonic functions that are not pluriharmonic.

\paragraph{\S A.4: limit values for germs with several branches.}
\emph{Interpretation: the topological count as reference case.} If every branch satisfies (H), the space of bounded harmonic functions modulo those with a common limit is isomorphic to $\R^{n_{\rm br}-1}$ (Theorem~\ref{th:rami}): there are as many independent limits as branches, according to the purely topological count. Section~\ref{sec:martin} shows that without (H) the number of minimal positive harmonic functions over the vertex can exceed the number of branches, even with connected link. In this language, a germ with a single branch satisfying (H) is \emph{harmonically normal}, a terminology proposed here, by analogy with the characterization of normal spaces through the extension of locally bounded holomorphic functions.

\emph{Limitations.} The Harnack condition on each branch is assumed, rather than established for arbitrary germs; a general definition of harmonic normality, not presupposing (H) on the branches, remains to be formulated.

\paragraph{\S A.5: monomial curves.}
\emph{Interpretation 1: an analytic invariant read off from harmonic functions.} For $y^{m_{\rm P}}=x^{n_{\rm P}}$ bounded harmonic functions are Hölder continuous at the vertex with optimal exponent $1/m_{\rm P}$, the inverse of the multiplicity (Theorem~A.10); the profile of unit flux has a correction term with exponent $2(n_{\rm P}-m_{\rm P})/m_{\rm P}$, determined by the Puiseux pair.

\emph{Interpretation 2: intrinsic harmonic functions and ambient restrictions.} A harmonic function on the curve is the restriction of the real part of an ambient holomorphic function if and only if, on the normalization, it is the real part of a function in the ring $\C\{t_{\rm P}^{m_{\rm P}},t_{\rm P}^{n_{\rm P}}\}$: the gap between intrinsic harmonic functions and ambient restrictions is measured by the semigroup of the branch.

\emph{Limitations.} The case is of complex dimension one and classical: through the normalization, harmonic theory reduces to that of the disc.

\paragraph{\S A.6: a real cone over a sextic curve.}
\emph{Interpretation 1: a real example with faster decay than the smooth case.} The unit link of the half-cone has length less than $2\pi$; the exponent $\alpha_{\rm L}=2\pi/\ell_{\rm L}$ is therefore greater than $1$, and every bounded harmonic function satisfies $u-u(0)=o(\rho)$ (Theorem~A.11). Even the simplest cone of real geometry has a harmonic theory different from that of the plane.

\emph{Interpretation 2: branches and extremals in an explicit example.} On the full cone $K_{\rm L}$ every bounded harmonic function has two branch limits, and the cone of nonnegative harmonic functions with zero trace on the outer links has exactly two extremal rays (Corollary~A.12): the count of \S A.4 is realized and coincides with the Martin count.

\emph{Interpretation 3: an additional real algebraic component.} The real locus of the homogenization $Q_{\rm L}$ contains, besides the cone, an additional line. This is an instance of the comparison between real and complex singularities, in which the real algebraic closure adds components that do not belong to the closure of the cone.

\emph{Limitations.} The cone is flat, with metric $d\rho^2+\rho^2d\sigma_{\rm L}^2$, and the harmonic theory is separable; the example illustrates but does not generalize.

\paragraph{\S A.7: metric induced by a homogeneous map.}
\emph{Interpretation: real normalizations as weighted horns.} A homogeneous map of degree $a_{\rm hom}$ with injective differential on the cross-section compares the induced metric with that of the source cone, up to the radial weight $r_{\rm hom}^{2(a_{\rm hom}-1)}$ (Lemma~A.13). Potential theory on the image can therefore be compared, at the level of energies and capacities, with that of a weighted cone, with the same structure as the horns of Section~\ref{sec:corni}. The cubic map from the quadric cone onto $K_{\rm L}$ is an example, and connects real normalizations with the methods of Sections~\ref{sec:capacita-profilo}--\ref{sec:corni-esatti}.

\emph{Limitations.} Injectivity of the differential on the whole cross-section is a strong hypothesis, which excludes normalizations with ramification.

\paragraph{\S A.8: counterexample to the angular estimate for the Gauss map.}
\emph{Interpretation: a limitation of metric approaches to tangent planes.} There is in general no Lipschitz estimate of tangent planes in the angular variable $x/|x|$ alone (Proposition~A.14). The effective directions do not necessarily lie in the link of the tangent cone, and a uniform estimate of this kind is not available.

\emph{Limitations.} The example excludes only this form of the estimate; it gives no information on weaker forms, for instance Hölder estimates or estimates along sequences.

\section{Preliminaries}\label{sec:prelim}

\subsection{Assumptions and notation}

Throughout, $X\subset\R^N$ is a subanalytic set, closed in a neighborhood of $0$, of pure dimension $d\ge2$, such that $X\setminus\{0\}$ is an analytic submanifold of $\R^N$ near $0$ and the pair $(X\setminus\{0\},\{0\})$ satisfies the Whitney condition at $0$. We use the notation of \cite[\S1.5]{V11}, summarized in the following table together with that introduced in this paper.

{\small
\begin{longtable}{L{4.6cm}L{10cm}}
\hline
Symbol & Meaning \\
\hline\endhead
$\rho(x):=|x|$ & Euclidean distance from $0$ \\
$r_0$ & radius in \cite[\S2]{V11} below which $|\nabla_X\rho|\ge\frac12$ and the product structure holds \\
$X^*$ & $X\cap B(0,r_0)\setminus\{0\}$ \\
$L_r$, $\Omega_r$, $A(s,t)$ & $X\cap S(0,r)$, $X\cap\{0<\rho<r\}$, $X\cap\{s<\rho<t\}$ \\
$g$, $\tau$, $\nabla_X$, $\Delta_X$ & induced metric on $X^*$, Riemannian measure, gradient, Laplace--Beltrami operator \\
$d_X$ & intrinsic distance \\
(H$_b$), $C_{H,b}(r)$, $C_H(r)$ & Harnack condition on the annuli $A(r/b,br)$, optimal constant at radius $r$, abbreviation $C_H(r)=C_{H,4}(r)$; cf. \cite[Def.~3.1]{V11} \\
(HS) & Harnack condition along a sequence of radii \cite[Def.~5.13]{V11} \\
$\mathcal P$ & positive harmonic functions on punctured neighborhoods of $0$ \\
$\mathcal K_R,G_R$ & positive harmonic functions on $\Omega_R$ vanishing continuously on $L_R$; positive generator when the Picard principle holds (\cite[Prop.~5.16]{V11}, Theorem 7.3) \\
(AD), $G_X$ & asymptotic decomposition and singular profile \cite[Def.~5.1]{V11} \\
$\Cap(E;O)$ & capacity of the condenser $(E,O)$ (Definition~\ref{def:capacita}) \\
$\Cap_X(r;R)$ & $\Cap(X\cap\{\rho\le r\};X\cap B(0,R))$ \\
$e_{r,R}$ & equilibrium potential of $A(r,R)$ \\
$\mathfrak F(u)$ & flux of a harmonic function toward the vertex (Definition~\ref{def:flusso}) \\
$L$, $h$, $V_L$ & cross-section of a horn, its induced metric, $\vol_h(L)$ \\
$m$ & ambient dimension of the cross-section, $L\subset\R^m\setminus\{0\}$ \\
$\gamma$, $\Pi$ & horn exponent, parameterization $\Pi(z,y)=(z^\gamma y,z)$ \\
$\beta$, $\kappa$ & $\gamma(d-1)$ and $2\gamma-2$; in Section~\ref{sec:collasso}, $\beta=\sum_id_i\gamma_i$ \\
$a_L$, $\Lambda_L$ & $\max_L|y|$ and $(1+a_L^2)^{1/2}$ \\
$\Cap^\Pi(\varepsilon;\zeta)$ & capacity between the levels $\{z\le\varepsilon\}$ and $\{z<\zeta\}$ of a horn \\
$S$, $\mathcal N$, $\vartheta$ & $\frac12|y|^2$, $|y|^2-|\nabla_hS|_h^2$, $1+\gamma^2z^\kappa\mathcal N$ \\
$\sigma_0$, $\eta$, $\sigma_1$ & $1-\beta$, remainder exponent, $\sigma_0+\eta$ \\
$W_\sigma$ & barrier $z^\sigma(1+\gamma\sigma S z^\kappa)$ \\
$K_\Psi$ & bi-Lipschitz constant of a homeomorphism $\Psi$ \\
$T_L$, $\langle\cdot\rangle$, $\mathfrak L_\xi$, $f_1$ & $|\nabla_hS|_h^2$, average on $L$, operator in Lemma~\ref{lemma:sviluppo}, first correction (\S\ref{sec:correzione}) \\
$X_{\boldsymbol\gamma}$, $Y_i$, $\Xi$, $\varsigma_{\boldsymbol\gamma}$ & multi-rate collapsing horn, link factors, radial factor, exact profile (\S\ref{sec:collasso}) \\
$X_k$, $\mathsf H$, $\mathsf B$, $\mathsf r$, $\Phi$, $\Sigma_z$ & thin-neck model and its functions, cross-sections at height $z$ (\S\ref{sec:collo}) \\
$\mathscr G$, $\mathscr K$, $\mathcal B$, $\mathcal T$ & Green function of $\Omega_R$, normalized kernel, normalized harmonic functions in $\mathcal K_R$, limits of the kernels (\S\ref{sec:martin}) \\
$A_{\mathsf k}$ & complex surfaces of type $A$ (\S\ref{sec:Ak}) \\
$\mathsf k$ & integer index of the singularity $A_{\mathsf k}$, distinct from $k$ in Section~\ref{sec:collo} \\
$\mathsf X,\mathsf Y,\mathsf Z$, $\mathsf U$, $\mathsf V$, $\mathsf m$ & ambient complex coordinates, $\mathsf X\pm i\mathsf Y$, $(\mathsf k+1)/2$ \\
$\varrho$, $\psi$, $\phi_1$, $\phi_2$, $\psi_*$ & $|\mathsf Z|$, logarithmic coordinate, angles, half-width of the balanced band \\
$\mathsf S_\varrho$, $\mathsf P$, $\mathsf J$, $\mathsf W$ & cross-section $\{|\mathsf Z|=\varrho\}$, meridional coefficient, determinant of the angular block, energy weight \\
$\mathfrak r$, $\mathfrak v_\pm$, $\mathfrak h$, $\mathfrak c$ & scale, channel variables, half-difference and coupling in the two-channel network model \\
(U), (T), $K_{\rm cart}$, $n_U$, $r_U$, $U^{\rm c}_i$, $\varphi^{\rm c}_i$, $s^{\rm c}_i$ & uniform chart and thickness conditions, constants, charts and scales (Definitions~\ref{def:U} and~\ref{def:T}) \\
$N,d$ & ambient dimension and pure real dimension (Introduction) \\
$r_{\mathrm{ext}},M_{\mathrm{ext}}$ & collar radius and outer neighborhood for capacity (2.3) \\
$E,O,\nu,dS,\mathfrak F_r$ & condenser plates (inner plate and outer domain), oriented normal, surface measure, flux through $L_r$ (2.3--2.5) \\
$\mathcal H^{d-1},B_X(x,s)$ & Hausdorff measure of dimension $d-1$ and intrinsic ball \\
$r_H,C_H$ & threshold and uniform constant of (H$_b$) for fixed $b$ (Theorem 3.1) \\
(H) & abbreviation of (H$_4$) \\
$q,u_{\mathrm b}$ & nonnegative coefficient and bounded remainder of decomposition (AD) (Introduction) \\
$X_{\gamma,L},z,y$ & exact generalized horn, height and point of the cross-section (Definition 2.1) \\
$p_{\rm cap},j,r_*$ & capacity exponents and radius of the estimate (3.4) \\
$\varepsilon,\zeta,\zeta_0$ & heights of the condensers and threshold of the horn chart (Section 4) \\
$\Psi,c_\Psi$ & bi-Lipschitz homeomorphism and $(2K_\Psi)^{-1}$ (4.4) \\
$\mathcal C,\Gamma,h_\Gamma$ & metric cone, unit link and its metric (4.5) \\
$\xi,\chi,\varsigma_L$ & exponent, angular test function and profile for a spherical cross-section (Section 5) \\
$\Theta_\pm,U_\varepsilon,U$ & combined barriers, solutions on truncations and harmonic limit (5.5--5.6) \\
$M_0,M_\delta,e_\delta$ & comparison constants and equilibrium potential of the collar (5.6) \\
$G^{(1)}$ & unit-flux profile, when the construction is specified (5.9) \\
$\Psi_{\xi,\chi},w_0$ & smooth remainder function in the variable $z^\kappa$ and its threshold (6.1) \\
$J_v,J^{(0)},J^{(1)}$ & flux oriented toward increasing heights and coefficients of its expansion (6.2--6.3) \\
$P_1,P_2,c_0,R_1,f_2,\sigma_2,M_1$ & approximants, constant in the first corrector, residual, second corrector and barrier parameters (6.3) \\
$\theta_1,\theta_2,c_2$ & toroidal angles and $\cos\theta_2$ (6.5) \\
$\alpha,b_L,L_{b_L}$ & radius of the spherical cross-section, translation vector and translated cross-section (6.6) \\
$p,n_i,d_i,\alpha_i,h_i,\boldsymbol\gamma$ & number of factors, ambient and intrinsic dimensions, radii and metrics of the factors, vector of exponents (7.1) \\
$\Xi_0,V_Y$ & limit of the radial factor and product of the volumes of the cross-sections (7.1) \\
$w,w_1,c_{\mathsf q},\varkappa$ & resistance coordinate, threshold, lower-bound constant of the ODE potential and decay exponent (7.2--7.3) \\
$\mathsf E_{i,j},\mathsf E_{\mathbf j},\lambda_{i,j},\lambda_{\mathbf j},\lambda_{\min}$ & eigenfunctions of the factors and product, eigenvalues and minimum spectral gap of the factors (7.3) \\
$u_{\mathbf j},\mathsf q_{\mathbf j},A_0,B_0,A,B$ & modal coefficients, ODE potential, coefficients of the constant mode and the average (7.3) \\
$X_{2,4},X_{3,3},\varpi_r$ & multiscale examples and cycle of the second factor (7.5) \\
$[\Xi]_{\beta-1},c_{\log}$ & resonant Puiseux coefficient and logarithmic coefficient (7.7) \\
$\mathsf n$ & integer indexing the resonance (8.1) \\
$k,z_*,t,\theta,\mathsf D,C_\rho$ & exponent and coordinates of the model, determinant of the meridional block, comparison constant between height and radius (9.1) \\
$\tau_0,\mathsf K_{\tau_0}$ & half-width of the neck band and its cross-sectional conductance (9.2) \\
$\varpi,\mathsf w_\varpi,\mathsf c_\varpi,I_0,J_0,V_0$ & parameter $z^k$, mass and energy weights, limiting integrals and mass of a lobe (9.3) \\
$\lambda_1(\Sigma_z),\mathcal E_z,f_\varpi,c_\pm,t_*$ & first positive eigenvalue of $-\Delta_{\Sigma_z}$, cross-sectional energy, test functions, limits on the lobes and reference point (Theorem~\ref{thm:gap}) \\
$\mathcal D_Z,u_Z,u_{Z,\varepsilon},c_Z,\tau_1$ & truncated domain, limiting and approximating harmonic functions, positive limit on the lobe, angular threshold (9.4) \\
$\widehat z,\mathcal R,\mathsf M_s,\mathsf Q_s$ & rescaled height, rectangle and rescaled metric/energy matrices (9.4) \\
$\mathfrak F^z,\mathfrak F^\rho,\ell,\mathcal A(s)$ & fluxes on the half-links toward increasing heights and radii, loss through the meridian and weighted area of the half-link (9.4) \\
$\widehat\Phi,\mathbf V,\mathsf j,\mathsf j_0,\mathfrak m,\bar u,\mathcal V$ & radial direction, radial variation field, density and its limit, mass, average and total variation of the density (9.4) \\
$\mathcal L_+,h_+,U_\infty,\mathsf t_\delta$ & limiting cone, link metric, limiting harmonic function, logarithmic cutoff (9.4) \\
$p_r,p_r^-$ & symmetric points of the lobes (Corollary~\ref{cor:harnack-collo}) \\
$\mathsf a,\mathsf b,\mathsf h$ & channels and half-difference of the reduced model (Remark 9.6) \\
$R_{\rm net}$ & radius of the initial interval in the contraction in Remark 9.6 \\
$r_c,o,\mathsf p,\mathscr G_\varepsilon$ & normalization radius and point, collar potential, Green function of the annulus (10.2) \\
$\Upsilon,\mathsf s_1,\mathsf s_2,u_\star,\mathscr K^\pm$ & cutoff, compactly supported sources, discontinuous harmonic function and limits of the kernels (Theorem 10.3) \\
$\omega,\upsilon_\varepsilon,\widetilde\omega_\varepsilon,\mu_\varepsilon,\pi_\varepsilon$ & normalized harmonic function, harmonic replacement, patched function, jump measure and representing probability measure (10.4) \\
$\varrho_0,K_{A_{\mathsf k}}$ & threshold and conductance of the balanced band (11.1) \\
$\mathfrak r_0,\mathfrak h_0$ & radial threshold and hypothetical limit of the half-difference (11.2) \\
$G_{\rm gr,\iota},U_{\rm gr,\iota},D_{\rm gr,\iota},\varphi_{\rm gr,\iota},\pi_{\rm gr,\iota}$ & graphs, planes, domains, functions and projections of Stasica's decomposition (\S A.1) \\
$K_S,c_T,r_T,s_T,V^*,\omega_d$ & constants, thickness scale, upper bound for link measure and volume of the unit ball (\S A.1) \\
$n_{\rm cart},K_{\rm cart},C_M,g_{\rm cart}$ & number of charts, metric constant, Moser constant, local metric matrix (\S A.2) \\
$V_{\rm an}$, $n_{\rm amb}$, $d_{\C}$, $r_{\rm an}$, $Y_{\rm reg}$, $\alpha_{\rm hol}$, $f_{\rm hol}$, $x_{\rm hol}$ & complex space, ambient and complex dimensions, radius and regular part, closed form, primitive and base point (\S A.3) \\
$\Omega^{(j)}$, $n_{\rm br}$, $\mathcal H_{\rm bd}$, $\mathcal H_{\rm zero}$, $\mathcal H_{\rm cont}$, $\ell_{{\rm br},j}$, $\operatorname{ev}_{\rm br}$, $V_{\rm quad}$ & branches, their number, spaces of harmonic germs, limits and evaluation, normal quadric (\S A.4) \\
$m_{\rm P},n_{\rm P},V_{\rm P},\nu_{\rm P},t_{\rm P},s_{\rm P},\theta_{\rm P}$ & exponents, curve, normalization, complex parameter, modulus and argument (\S A.5) \\
$G_{\rm P}^{(1)},f_{\rm P},\mathcal S_{\rm P},j_{\rm P}$ & profile, holomorphic function, numerical semigroup and integer exponent (\S A.5) \\
$\Gamma_{\rm L},q_{\rm L},\theta_{\rm L},s_{\rm L},K_{\rm L},K_{\rm L,\pm},Q_{\rm L}$ & real curve, parameterization, parameters, cone and its half-cones, homogenized polynomial (\S A.6) \\
$\ell_{\rm L},\sigma_{\rm L},\alpha_{\rm L},\omega_{\rm L},\zeta_{\rm L},c_{\rm L},G_{\rm L,+}^{(1)}$ & length, arc length, exponent, normalized link, conformal coordinate, cosine and unit-flux profile (\S A.6) \\
$Y_{\rm hom}$, $g_{Y_{\rm hom}}$, $\nu_{\rm hom}$, $a_{\rm hom}$, $r_{\rm hom}$, $\widehat E_{\R}$ & source cone, metric, homogeneous map, degree, radius and real cone in the example (\S A.7) \\
$\delta_{\rm Gr},\omega_{\rm G},n_{\rm G}$ & Grassmannian distance, direction of the spherical cross-section and unit normal (\S A.8) \\
\hline

$r,R,s,t$ & working and truncation radii, when used as radial arguments (\S2) \\
$z_0,z_1$ & working levels of the horns and outer level of the construction, fixed in the respective statements (Sections 5--7) \\
$F$ & radial metric factor $1+\gamma^2a_L^2 z^\kappa$ in the spherical-cross-section calculation (Proposition 5.3); not the flux $\mathfrak F$ \\
$\mathsf q,\mathsf b_{\rm ODE},\ell_\infty$ & scalar potential, barrier and limit of the solution in Lemma 7.2 \\
$\mathbf j,\Lambda$ & modal multi-index and spectral threshold (Theorem 7.3) \\
$\Xi_{2,4},\Xi_{3,3},\varsigma_{2,4},\varsigma_{3,3}$ & radial factors and profiles in the two examples of Proposition 7.5 \\
$\rho_{2,4},\rho_{3,3}$ & respective Euclidean distances (Proposition 7.5) \\
$\mathsf w_0,f_{\rm lim}$ & limiting density on the lobe and piecewise limiting function (Theorem~\ref{thm:gap}) \\
$Z,\mathsf f$ & truncation height and odd boundary data (Theorem~\ref{thm:discontinuita}) \\
$N_0,N_1,\mathsf p_i$ & number of polynomials, uniform bound on the number of monotonicity intervals and sign polynomials (Theorem~\ref{thm:discontinuita}, Step 5) \\
$I,\eta_1,s_n$ & lobe interval, rescaled threshold and blow-up scales (Theorem~\ref{thm:discontinuita}, Steps 6--7) \\
$\mathsf v$ & function $1-u_Z$ in Step 6 of 9.4, where it appears as an auxiliary function \\
$r_-,r_+,A_K,\mathscr G_{A_K},C_K$ & radii, fixed annulus, its Green function and constant in the local estimates (Proposition~\ref{prop:green}) \\
$\mathsf y_n^\pm,\ell_\pm,q_\pm,\mathsf b_\pm,\mathsf b_\varepsilon$ & sequences of poles, cluster values, coefficients, bounded remainders and harmonic residual of the truncation (Theorem 10.3) \\
$\mathscr P_{\delta,\varepsilon},\mathscr P_\varepsilon$ & potentials of the jump measure with respect to the truncated and global Green functions (Theorem 10.4) \\
$K_n,\pi,a_\upsilon$ & increasing exhaustion, limiting probability measure and value of the dominated function at $o$ (Theorem 10.4) \\
$\mathsf a_{\rm M},\lambda_{\rm M},\Lambda_{\rm M}$ & matrix and ellipticity bounds in the formulation of Moser's Harnack inequality (\S A.1) \\
$\mathsf a_{{\rm cart},i},\mathcal V_i^{\rm c},\mathcal V_x^{\rm c}$ & divergence-form matrix in the chart and subcharts of radius one (Lemma~\ref{lemma:U-H}, Proposition~\ref{prop:T-U}) \\
$G_x,u_x,x_i$ & local graph, projection of the center and separated centers (Proposition~\ref{prop:T-U}) \\
\hline
\end{longtable}
}

For two positive functions $f,g$ we write $f\asymp g$ if there exists $C\ge1$ with $C^{-1}g\le f\le Cg$ on the set under consideration. The letter $C$ without subscripts denotes a positive constant whose value may change from line to line. Dummy integration variables, summation indices, local coordinates and auxiliary constants in a proof have the meaning specified in their respective contexts. The unsubscripted parameter $b$ in (H$_b$) is real and at least two. For $0<br<r_0$ we denote by $C_{H,b}(r)$ the optimal constant in the inequality on the link $L_r$ for positive harmonic functions on the annulus $A(r/b,br)$; we set $C_H(r)=C_{H,4}(r)$. The symbol $C_H$ without an argument denotes, when (H$_b$) holds with fixed $b$, a uniform constant for this inequality. When referring to the positive cone and its extreme rays we consider $\mathcal K_R\cup\{0\}$; its normalized section is $\mathcal B$. The brackets $\langle f\rangle$ denote the average on $L$, and brackets with two arguments denote the inner product.

\subsection{Classical tools on compact manifolds}\label{sec:classici}

We use three classical facts on smooth compact Riemannian manifolds. They follow from the Euclidean statements of \cite{GT}, applied in finitely many coordinate charts with a subordinate partition of unity: in a chart the Laplace--Beltrami operator is in divergence form, with smooth coefficients and uniformly elliptic principal part. The sign convention is $\Delta\le0$.
\begin{enumerate}[label=(T\arabic*)]
\item \label{T:dirichlet} \emph{Dirichlet problem.} Let $A$ be a compact connected Riemannian manifold with smooth nonempty boundary, $f\in C^\infty(A)$ and $\varphi\in C^\infty(\partial A)$. There is a unique $u\in C^\infty(A)$ with $\Delta u=f$ in the interior and $u=\varphi$ on $\partial A$; if $f=0$, $u$ minimizes the Dirichlet energy among the functions of $W^{1,2}(A)$ with trace $\varphi$. Writing $u=\Phi+w$, with $\Phi\in C^\infty(A)$ extending $\varphi$, the weak solution $w\in W^{1,2}_0(A)$ is given by the Lax--Milgram theorem \cite[Theorem~5.8]{GT}; coercivity is the Poincaré inequality on $W^{1,2}_0(A)$, obtained by contradiction from the compactness of the embedding $W^{1,2}\subset L^2$ in charts \cite[Theorem~7.26]{GT}. Smoothness up to the boundary follows from \cite[Theorems~8.10 and~8.13]{GT}, and uniqueness from the energy identity.
\item \label{T:poisson} \emph{Poisson equation on closed manifolds.} Let $Y$ be a compact connected Riemannian manifold without boundary and $g\in C^\infty(Y)$ with $\int_Yg\,d\vol=0$. There is $f\in C^\infty(Y)$, unique up to an additive constant, with $\Delta f=g$. The proof is that of~\ref{T:dirichlet} on the closed subspace of $W^{1,2}(Y)$ of functions with zero mean, where the Poincaré--Wirtinger inequality follows from the same compactness argument and from the connectedness of $Y$; smoothness follows from \cite[Theorem~8.10]{GT}.
\item \label{T:sobolev} \emph{Sobolev embedding.} Let $Y$ be as in~\ref{T:poisson} and let $M$ be an integer with $4M>\dim Y$. There is $C>0$ such that $\sup_Y|f|\le C\|(1-\Delta)^Mf\|_{L^2(Y)}$ for every $f\in C^\infty(Y)$. Indeed, $\sup_Y|f|\le C\|f\|_{W^{2M,2}(Y)}$ by \cite[Theorem~7.26]{GT} in charts; the interior estimate of \cite[Theorem~8.10]{GT} gives $\|v\|_{W^{j+2,2}(Y)}\le C\bigl(\|(1-\Delta)v\|_{W^{j,2}(Y)}+\|v\|_{L^2(Y)}\bigr)$, and $\|v\|_{L^2(Y)}\le\|(1-\Delta)v\|_{L^2(Y)}$ because $-\Delta\ge0$; iterating $M$ times gives the claim.
\end{enumerate}

\subsection{Results from \texorpdfstring{\cite{V11}}{[V11]} used below}

We recall the results of \cite{V11} on which this paper relies.
\begin{enumerate}[label=(R\arabic*)]
\item \label{R:distanza} \cite[Lemma~2.4]{V11}: $|x|\le d_X(x,0)\le2|x|$ for $x\in X^*$.
\item \label{R:picard} \cite[Prop.~5.16]{V11}: if $X$ has connected link and satisfies (HS), for every $R\in(0,r_0]$ there exists $G_R\in\mathcal K_R$ such that $\mathcal K_R=\{tG_R:t>0\}$, and $G_R\to+\infty$ uniformly as $\rho\to0$.
\item \label{R:AD} \cite[Theorem~5.17]{V11}: if $X$ has connected link and satisfies (HS), then $X$ satisfies (AD) with $G_X=G_{r_0}$, bounded harmonic functions extend continuously to $0$, and the Bôcher property holds. Condition (H$_b$) implies (HS).
\item \label{R:corni} \cite[Prop.~4.10, Cor.~4.12 and proof of Theorem~4.15]{V11}: generalized horns, germs bi-Lipschitz equivalent to a generalized horn with respect to the intrinsic metrics, and germs bi-Lipschitz equivalent with respect to the intrinsic metrics to a cone over a smooth compact connected link satisfy (H$_b$) for a suitable $b\ge2$.
\end{enumerate}

\begin{definition}\label{def:corno}
Let $m\ge d\ge2$, let $L\subset\R^m\setminus\{0\}$ be a compact connected analytic submanifold of dimension $d-1$, and let $\gamma\ge1$ be rational. The \emph{generalized horn} of exponent $\gamma$ and cross-section $L$ is
\[
X_{\gamma,L}:=\{(z^\gamma y,z)\in\R^m\times\R:\ y\in L,\ z\ge0\},
\]
with parameterization $\Pi(z,y):=(z^\gamma y,z)$.
\end{definition}

For $m=d$ this is the definition in \cite[Def.~4.8]{V11}.

\begin{remark}\label{oss:m-maggiore}
The proof of \cite[Prop.~4.10]{V11} uses, for the cross-section $L$, only the bound $|y|\le a_L$ and the Riemannian manifold $(L,h)$. It therefore holds without modification for $m>d$, as do \cite[Theorem~4.15]{V11} and condition (H$_b$). The set $X_{\gamma,L}$ is subanalytic, since, writing $\gamma=p/q$ with positive integers $p,q$, the graph of $(s,y)\mapsto(s^py,s^q)$ on $[0,1]\times L$ is compact and semianalytic, so that its projection is subanalytic near $0$; it is analytic away from $0$; it satisfies the Whitney condition at $0$, since the vectors $\Pi(z,y)/|\Pi(z,y)|$ and $\partial_z\Pi(z,y)/|\partial_z\Pi(z,y)|=(\gamma z^{\gamma-1}y,1)/|(\gamma z^{\gamma-1}y,1)|$ have distance $O(z^{\gamma-1})$ if $\gamma>1$ and coincide if $\gamma=1$.
\end{remark}

\subsection{Capacity}

\begin{definition}\label{def:capacita}
Fix $r_{\mathrm{ext}}>r_0$ such that $X$ is closed in a neighborhood of $\overline B(0,r_{\mathrm{ext}})$, smooth away from $0$, with $|\nabla_X\rho|\ge\frac12$ on $\{0<\rho<r_{\mathrm{ext}}\}$; this is possible by reducing $r_0$. Set $M_{\mathrm{ext}}:=X\cap B(0,r_{\mathrm{ext}})$. Let $E\subset O\subset M_{\mathrm{ext}}$, with $E\setminus\{0\}$ relatively closed in $M_{\mathrm{ext}}\setminus\{0\}$, $O$ relatively open in $M_{\mathrm{ext}}$ and $\overline O\subset\overline B(0,r_0)$. The admissibility conditions concern only the regular part; capacity is unchanged if the vertex is added to or removed from the plates, provided the resulting pairs satisfy these conditions. A function $f\colon M_{\mathrm{ext}}\setminus\{0\}\to\R$ is \emph{admissible} for $(E,O)$ if it is locally Lipschitz, equals $1$ on $E\setminus\{0\}$ and $0$ on $(M_{\mathrm{ext}}\setminus\{0\})\setminus O$. The \emph{capacity} of the condenser $(E,O)$ is
\[
\Cap(E;O):=\inf\Bigl\{\int_{M_{\mathrm{ext}}\setminus\{0\}}|\nabla_Xf|^2\,d\tau:\ f\ \text{admissible for }(E,O)\Bigr\}.
\]
Capacity does not depend on the choice of $r_{\mathrm{ext}}$, since admissible functions vanish outside $O$.
\end{definition}

It follows from the definition that $\Cap(E;O)$ increases when $E$ is enlarged and decreases when $O$ is enlarged.

\begin{lemma}\label{lemma:dirichlet}
Let $0<r<R\le r_0$. There exists a unique function $e_{r,R}$ continuous on $\overline{A(r,R)}$, harmonic on $A(r,R)$, with $e_{r,R}=1$ on $L_r$ and $e_{r,R}=0$ on $L_R$; it is smooth up to the boundary and, if $\nu$ denotes the outward unit normal to $A(r,R)$,
\[
\Cap_X(r;R)=\int_{A(r,R)}|\nabla_Xe_{r,R}|^2\,d\tau=\int_{L_r}\partial_\nu e_{r,R}\,dS ,
\qquad \partial_\nu e_{r,R}>0\ \text{on }L_r .
\]
\end{lemma}

\begin{proof}
Since $|\nabla_X\rho|\ge\frac12$ on $X^*$, the links $L_r$ and $L_R$ are smooth submanifolds and $\overline{A(r,R)}$ is a compact manifold with smooth boundary; by the product structure every connected component meets both $L_r$ and $L_R$, and existence, uniqueness and boundary regularity follow from~\ref{T:dirichlet} applied on each component. The function equal to $e_{r,R}$ on $\overline{A(r,R)}$, to $1$ on $\{0<\rho\le r\}$ and to $0$ on $M_{\mathrm{ext}}\cap\{\rho\ge R\}$ is admissible, with energy $\int_{A(r,R)}|\nabla_Xe_{r,R}|^2d\tau$. If $f$ is admissible, its restriction to $\overline{A(r,R)}$ is Lipschitz, equals $1$ on $L_r$ and $0$ on $L_R$; by the Dirichlet principle in~\ref{T:dirichlet}, $\int_{M_{\mathrm{ext}}\setminus\{0\}}|\nabla_Xf|^2d\tau\ge\int_{A(r,R)}|\nabla_Xf|^2d\tau\ge\int_{A(r,R)}|\nabla_Xe_{r,R}|^2d\tau$. Thus $\Cap_X(r;R)$ is the energy of $e_{r,R}$. Green's formula gives $\int|\nabla_Xe_{r,R}|^2\,d\tau=\int_{\partial A(r,R)}e_{r,R}\partial_\nu e_{r,R}\,dS=\int_{L_r}\partial_\nu e_{r,R}\,dS$. Since $0<e_{r,R}<1$ in $A(r,R)$ by the strong maximum principle and $e_{r,R}$ attains its maximum on $L_r$, the Hopf lemma \cite[Lemma~3.4]{GT} gives $\partial_\nu e_{r,R}>0$ on $L_r$.
\end{proof}

\begin{definition}\label{def:flusso}
Let $u$ be harmonic on $\Omega_R$. For $0<r<R$ set $\mathfrak F_r(u):=\int_{L_r}\partial_\nu u\,dS$, where $\nu$ is the unit normal to $L_r$ directed toward $0$. By the divergence theorem on $A(r_1,r_2)$, $\mathfrak F_r(u)$ is independent of $r$ and is denoted by $\mathfrak F(u)$.
\end{definition}

The same definition, with the same notation, is used for level surfaces of a smooth function with nonvanishing gradient, for example the levels $\{z=\varepsilon\}$ of a horn: the same flux is obtained through any smooth separating level hypersurface that bounds, together with an inner link, a compact annulus, with normal directed toward the vertex.

\begin{lemma}[bi-Lipschitz transport]\label{lemma:trasporto}
Let $X$ and $Y$ be as in this section, of dimension $d$, and let $\Psi$ be a homeomorphism from a neighborhood $U$ of $0$ in $X$ onto a neighborhood of $0$ in $Y$, with $\Psi(0)=0$, that is $K_\Psi$-bi-Lipschitz with respect to the intrinsic metrics. Suppose $U$ is chosen sufficiently small that $U\setminus\{0\}\subset X^*$ and $\Psi(U)\setminus\{0\}\subset Y^*$. If $f$ is locally Lipschitz on $Y^*\cap\Psi(U)$, then $f\circ\Psi$ is locally Lipschitz on $U\setminus\{0\}$ and
\[
\int_{U\setminus\{0\}}|\nabla_X(f\circ\Psi)|^2\,d\tau_X\le K_\Psi^{d+2}\int_{\Psi(U)\setminus\{0\}}|\nabla_Yf|^2\,d\tau_Y .
\]
\end{lemma}

\begin{proof}
On the regular part, the infinitesimal metric of the intrinsic distance is the induced Riemannian metric. In regular charts $\Psi$ and $\Psi^{-1}$ are locally Lipschitz and their infinitesimal metric constants are at most $K_\Psi$. By Rademacher's theorem \cite[3.1.6]{Fed} in local charts, $\Psi$ is differentiable almost everywhere, with $|d\Psi|\le K_\Psi$ and $|d\Psi^{-1}|\le K_\Psi$; since both maps are locally bi-Lipschitz in regular charts, the exceptional sets can be pulled back and pushed forward without changing their null-measure property, and at corresponding points of differentiability the differentials are inverse linear maps with singular values in $[K_\Psi^{-1},K_\Psi]$ for the induced Riemannian metrics; it follows that $|\nabla_X(f\circ\Psi)|\le K_\Psi|\nabla_Yf|\circ\Psi$ almost everywhere and, by the area formula \cite[3.2.3]{Fed}, $\int\varphi\circ\Psi\,d\tau_X\le K_\Psi^d\int\varphi\,d\tau_Y$ for every measurable $\varphi\ge0$. Combining these yields the estimate.
\end{proof}

\section{Capacity and profile}\label{sec:capacita-profilo}

\begin{theorem}\label{thm:capacita-profilo}
Suppose $X$ has connected link and satisfies (H$_b$) for some $b\ge2$, with constant $C_H$ and radius $r_H$. Let $R\in(0,r_0]$ and let $G_R$ be as in~\ref{R:picard}. Then, for $0<r<\min\{r_H,R/b\}$ and $x\in L_r$,
\[
\frac{\mathfrak F(G_R)}{C_H\,\Cap_X(r;R)}\le G_R(x)\le\frac{C_H\,\mathfrak F(G_R)}{\Cap_X(r;R)} ,
\qquad \mathfrak F(G_R)>0 .
\]
\end{theorem}

\begin{proof}
Apply Green's second identity to $G_R$ and $e:=e_{r,R}$ on the compact manifold with boundary $\overline{A(r,R)}$; both functions are harmonic and smooth up to the boundary, since $G_R$ is smooth on $\Omega_R$ and, by boundary regularity, up to $L_R$. This gives
\[
\int_{L_r\cup L_R}\bigl(G_R\,\partial_\nu e-e\,\partial_\nu G_R\bigr)\,dS=0 .
\]
On $L_R$ both functions vanish; on $L_r$ we have $e=1$, and the outward normal to $A(r,R)$ is directed toward $0$. Hence
\begin{equation}\label{eq:media-pesata}
\int_{L_r}G_R\,\partial_\nu e\,dS=\mathfrak F(G_R).
\end{equation}
By Lemma~\ref{lemma:dirichlet}, $\partial_\nu e>0$ on $L_r$ and $\int_{L_r}\partial_\nu e\,dS=\Cap_X(r;R)$. Thus $\mathfrak F(G_R)/\Cap_X(r;R)$ is a positively weighted average of $G_R$ on $L_r$, and $\mathfrak F(G_R)>0$ since $G_R>0$. We have
\[
\min_{L_r}G_R\le\frac{\mathfrak F(G_R)}{\Cap_X(r;R)}\le\max_{L_r}G_R .
\]
Since $br<R$, the function $G_R$ is positive and harmonic on $A(r/b,br)$, and (H$_b$) gives $\max_{L_r}G_R\le C_H\min_{L_r}G_R$. Combining the two inequalities proves the statement.
\end{proof}

\begin{remark}\label{oss:classico}
Theorem~\ref{thm:capacita-profilo} adapts a classical comparison between Green function and capacity to the vertex: on a manifold, for a regular pole, the reciprocal of the capacity of a condenser lies between the minimum and maximum of the Green function on the inner boundary \cite[Prop.~4.1]{Gri}, and on model manifolds the capacity of a ball is given by a radial integral \cite[Example~4.2]{Gri}, the method used in Sections~\ref{sec:corni}--\ref{sec:risonanze}; for metric spaces with a doubling measure and a Poincaré inequality, see \cite[Theorem~1.5]{BBL20} and \cite[Theorem~7.1]{BBL23}. Here the pole is the vertex, outside the regular part, and only the annular condition (H$_b$) is used, without doubling or Poincaré assumptions on balls centered at $0$.
\end{remark}

\begin{remark}\label{oss:HS}
If $X$ satisfies only (HS), with radii $r_k$ and constant $C$, after discarding finitely many radii so that $br_k<R$, the same proof gives the conclusion for $r=r_k$ with $C$ in place of $C_H$, that is, along the Harnack sequence.
\end{remark}

\begin{corollary}\label{cor:ordine-capacita}
Under the assumptions of Theorem~\ref{thm:capacita-profilo}, if there exist $p_{\rm cap}\ge0$ and $j\in\R$ with $\Cap_X(r;R)\asymp r^{p_{\rm cap}}\,|\log r|^{-j}$ as $r\to0$, then $G_R\asymp\rho^{-p_{\rm cap}}|\log\rho|^{j}$ on $\Omega_{r_*}$ for a suitable $r_*>0$, and the same holds for every profile $G_X$ of (AD).
\end{corollary}

\begin{proof}
The first assertion follows from Theorem~\ref{thm:capacita-profilo}. Let $G_X$ be any profile of (AD). Applying (AD) with this profile to $G_R$ gives $G_R=qG_X+b$ near $0$, with $b$ bounded and harmonic. By~\ref{R:picard}, $G_R\to+\infty$ uniformly, so $q>0$. Thus $G_X=q^{-1}G_R-q^{-1}b\asymp G_R$ near $0$.
\end{proof}

\section{Growth order of the profile for horns and bi-Lipschitz equivalent germs}\label{sec:corni}

In this section $X=X_{\gamma,L}$ is a generalized horn (Definition~\ref{def:corno}), $\beta=\gamma(d-1)$, $\kappa=2\gamma-2$, $a_L=\max_L|y|$, $\Lambda_L=(1+a_L^2)^{1/2}$, $V_L=\vol_h(L)$. We identify $(0,1)\times L$ with its image under $\Pi$ and write functions on $X\cap\{0<z<1\}$ as functions of $(z,y)$. Set $\zeta_0:=\min\{1,r_0/\Lambda_L\}$; by Lemma~\ref{lemma:metrica-corno}, $\{0<z<\zeta_0\}\subset X^*$. For $0<\varepsilon<\zeta\le\zeta_0$ set
\[
\Cap^\Pi(\varepsilon;\zeta):=\Cap\bigl(\{0<z\le\varepsilon\};\{0<z<\zeta\}\bigr).
\]

\begin{lemma}\label{lemma:metrica-corno}
On $(0,1)\times L$ the induced metric is
\[
\Pi^*g=\bigl(1+\gamma^2z^{\kappa}|y|^2\bigr)dz^2+2\gamma z^{2\gamma-1}\langle y,dy\rangle\,dz+z^{2\gamma}h ,
\]
its volume form is $z^{\beta}\vartheta^{1/2}\,dz\,d\vol_h$, with $\vartheta=1+\gamma^2z^\kappa\mathcal N$ and $\mathcal N=|y|^2-|\nabla_hS|_h^2\ge0$, and for every locally Lipschitz function $f$ we have, almost everywhere, $|\nabla_Xf|^2\ge(\partial_zf)^2/(1+\gamma^2z^\kappa|y|^2)$. Moreover, $z\le\rho\le\Lambda_Lz$.
\end{lemma}

\begin{proof}
The metric formula follows from $d\Pi=(\gamma z^{\gamma-1}y\,dz+z^\gamma dy,\,dz)$, with $\langle y,dy\rangle=dS$. The angular block is $z^{2\gamma}h$, and the Schur complement of the coefficient $g_{zz}$ with respect to it is
\[
1+\gamma^2z^{\kappa}|y|^2-\gamma^2z^{4\gamma-2}z^{-2\gamma}|\nabla_hS|_h^2=\vartheta ;
\]
hence the determinant relative to $dz\,d\vol_h$ is $z^{2\beta}\vartheta$. The gradient inequality is the Cauchy--Schwarz inequality $(\partial_zf)^2=\langle\nabla_Xf,\partial_z\rangle^2\le|\nabla_Xf|^2|\partial_z|^2$, with $|\partial_z|^2=1+\gamma^2z^\kappa|y|^2$. Finally, $\rho^2=z^2(1+z^{\kappa}|y|^2)$ and $z\le1$. That $\mathcal N\ge0$ follows from $|\nabla_hS|_h\le|y|$, since $\nabla_hS$ is the projection of $y$ onto the tangent space of $L$.
\end{proof}

\begin{proposition}\label{prop:capacita-corno}
For $0<\varepsilon<\zeta\le\zeta_0$,
\[
\frac{V_L}{\displaystyle\int_\varepsilon^\zeta\bigl(1+\gamma^2a_L^2t^{\kappa}\bigr)t^{-\beta}\,dt}
\ \le\ \Cap^\Pi(\varepsilon;\zeta)\ \le\
\frac{V_L}{\displaystyle\int_\varepsilon^\zeta t^{-\beta}\,dt}.
\]
\end{proposition}

\begin{proof}
\emph{Upper bound.} Let $f(z,y):=\varphi(z)$ with $\varphi(z)=1$ for $z\le\varepsilon$, $\varphi(z)=0$ for $z\ge\zeta$ and $\varphi(z)=I^{-1}\int_z^\zeta t^{-\beta}dt$ for $\varepsilon\le z\le\zeta$, where $I:=\int_\varepsilon^\zeta t^{-\beta}dt$. The function is admissible. By Lemma~\ref{lemma:metrica-corno}, $|\nabla f|^2=g^{zz}\varphi'^2$ and $g^{zz}=\vartheta^{-1}$, since $g^{zz}$ is the inverse of the Schur complement. Hence
\[
\int|\nabla f|^2\,d\tau=\int_L\int_\varepsilon^\zeta\vartheta^{-1/2}z^\beta\varphi'(z)^2\,dz\,d\vol_h
\le V_L\int_\varepsilon^\zeta z^\beta I^{-2}z^{-2\beta}\,dz=\frac{V_L}{I},
\]
since $\vartheta\ge1$.

\emph{Lower bound.} Let $f$ be admissible. For almost every $y\in L$ the function $t\mapsto f(t,y)$ is Lipschitz on $[\varepsilon,\zeta]$, with $f(\varepsilon,y)-f(\zeta,y)=1$, and the differential estimate holds almost everywhere along the fiber, by Fubini's theorem in local charts. By Lemma~\ref{lemma:metrica-corno}, $|\nabla f|^2\ge(\partial_zf)^2/(1+\gamma^2a_L^2z^\kappa)$ and the volume density is at least $z^\beta$. The Cauchy--Schwarz inequality gives
\[
1\le\Bigl(\int_\varepsilon^\zeta|\partial_zf(t,y)|\,dt\Bigr)^2\le\int_\varepsilon^\zeta\frac{t^\beta(\partial_zf)^2}{1+\gamma^2a_L^2t^\kappa}\,dt\cdot\int_\varepsilon^\zeta(1+\gamma^2a_L^2t^\kappa)t^{-\beta}\,dt .
\]
Integrating in $y$ with respect to $d\vol_h$ yields the lower bound.
\end{proof}

\begin{theorem}\label{thm:ordine-corni}
Let $X=X_{\gamma,L}$ and let $G_X$ be a singular profile in (AD). Then, near $0$,
\[
G_X\asymp\rho^{1-\beta}\ \ \text{if }\beta>1,\qquad G_X\asymp\log(1/\rho)\ \ \text{if }\beta=1 .
\]
The case $\beta=1$ occurs only for $d=2$ and $\gamma=1$.
\end{theorem}

\begin{proof}
By~\ref{R:corni} and Remark~\ref{oss:m-maggiore}, $X$ satisfies (H$_b$) for some $b$. Fix $R\le\zeta_0$. By Lemma~\ref{lemma:metrica-corno} we have the inclusions
\[
\{0<\rho\le r\}\subset\{0<z\le r\}\subset\{0<\rho\le\Lambda_Lr\},\qquad
\{0<z<R/\Lambda_L\}\subset\Omega_R\subset\{0<z<R\},
\]
and monotonicity of capacity gives
\[
\Cap^\Pi(r/\Lambda_L;R)\le\Cap_X(r;R)\le\Cap^\Pi(r;R/\Lambda_L).
\]
If $\beta>1$, by Proposition~\ref{prop:capacita-corno} both extreme terms are comparable to $r^{\beta-1}$ as $r\to0$: indeed $\int_\varepsilon^\zeta t^{-\beta}dt\sim\varepsilon^{1-\beta}/(\beta-1)$, and $\int_\varepsilon^\zeta t^{\kappa-\beta}dt$ is $O(\varepsilon^{1-\beta})$. If $\beta=1$, then $\kappa=0$ and the integrals are comparable to $\log(1/\varepsilon)$. Corollary~\ref{cor:ordine-capacita} gives the conclusion. Finally, $\beta=\gamma(d-1)=1$ with $\gamma\ge1$ and $d\ge2$ implies $d=2$ and $\gamma=1$.
\end{proof}

\begin{corollary}\label{cor:bilip}
Let $X$ be as in \S\ref{sec:prelim} and let $\Psi$ be a homeomorphism from a neighborhood of $0$ in $X$ onto a neighborhood of $0$ in $X_{\gamma,L}$, with $\Psi(0)=0$, bi-Lipschitz with respect to the intrinsic metrics. Then every singular profile in (AD) for $X$ satisfies the same estimates as in Theorem~\ref{thm:ordine-corni}.
\end{corollary}

\begin{proof}
The link of $X$ is connected by the radial product structure and the homeomorphism of punctured germs. By~\ref{R:corni}, $X$ satisfies (H$_b$) for some $b$. Let $Y:=X_{\gamma,L}$ and let $K_\Psi$ be the constant of $\Psi$. By~\ref{R:distanza}, applied to $X$ and $Y$, we have $(2K_\Psi)^{-1}\rho_X(x)\le\rho_Y(\Psi(x))\le2K_\Psi\rho_X(x)$. Setting $c_\Psi:=(2K_\Psi)^{-1}$, for $r$ small relative to $R$ and $R$ small we have the inclusions
\[
\{0<\rho_Y\le c_\Psi r\}\subset\Psi(\{0<\rho_X\le r\})\subset\{0<\rho_Y\le r/c_\Psi\},\quad
\Omega^Y_{c_\Psi R}\subset\Psi(\Omega^X_R)\subset\Omega^Y_{R/c_\Psi}.
\]
Choose $R$ and then $r<c_\Psi^2R$ sufficiently small that the closures of the outer neighborhoods used are contained in the respective domains of $\Psi$ and $\Psi^{-1}$ and within the working radius of each germ. After transport, the admissible functions are extended by zero outside the neighborhood on which the map is defined; this preserves local Lipschitz continuity, since their supports are contained in that neighborhood.
If $f$ is admissible for $(\{0<\rho_Y\le r/c_\Psi\};\Omega^Y_{c_\Psi R})$, then $f\circ\Psi$ is admissible for $(\{0<\rho_X\le r\};\Omega^X_R)$, and Lemma~\ref{lemma:trasporto} gives $\Cap_X(r;R)\le K_\Psi^{d+2}\Cap_Y(r/c_\Psi;c_\Psi R)$. The same argument with $\Psi^{-1}$ gives $\Cap_Y(c_\Psi r;R/c_\Psi)\le K_\Psi^{d+2}\Cap_X(r;R)$. By the proof of Theorem~\ref{thm:ordine-corni}, the two terms for $Y$ have the same order, $r^{\beta-1}$ or $1/\log(1/r)$; hence so does $\Cap_X(r;R)$, and Corollary~\ref{cor:ordine-capacita} gives the conclusion.
\end{proof}

\begin{corollary}\label{cor:coni}
Let $\mathcal C\subset\R^N$ be the cone over a smooth compact connected submanifold $\Gamma\subset S^{N-1}$ of dimension $d-1$, and let $X$ be a germ bi-Lipschitz equivalent to $\mathcal C$ with respect to the intrinsic metrics as in \cite[Cor.~4.12]{V11}. Then $G_X\asymp\rho^{2-d}$ if $d\ge3$ and $G_X\asymp\log(1/\rho)$ if $d=2$.
\end{corollary}

\begin{proof}
In the coordinates $x=t\upsilon$, $\upsilon\in\Gamma$, the metric of $\mathcal C$ is $dt^2+t^2h_\Gamma$, with $h_\Gamma$ the induced metric on $\Gamma$, and $\rho=t$. The proof of Proposition~\ref{prop:capacita-corno}, with $\gamma=1$ and vanishing mixed term, gives $\Cap_{\mathcal C}(r;R)=\vol(\Gamma)\bigl(\int_r^Rt^{1-d}dt\bigr)^{-1}$. The conclusion follows as in Corollary~\ref{cor:bilip}, using~\ref{R:corni} for (H$_b$).
\end{proof}

\begin{example}\label{es:toro-cuspidi}
(a) For the horn with toroidal cross-section of \cite[\S6.1]{V11}, with $d=3$, Theorem~\ref{thm:ordine-corni} gives $G_X\asymp\rho^{1-2\gamma}$. (b) The elliptic cusps of \cite[Example~4.14]{V11} with exponent $2\gamma$, $\gamma>1$ rational, are horns with exponent $\gamma$ and an ellipsoid as cross-section; we have $G_X\asymp\rho^{1-\gamma(d-1)}$, even when the semiaxes are not all equal. This answers question~(Q1) concerning the growth order of the profile in the case of elliptic cusps.
\end{example}

\section{Exact horns: leading coefficient and flux}\label{sec:corni-esatti}

In this section $X=X_{\gamma,L}$ with $\gamma>1$; thus $\kappa>0$ and $\beta>1$. We use the notation of Section~\ref{sec:corni}.

\subsection{Metric components and Laplacian}

\begin{lemma}\label{lemma:inversa}
In local coordinates $(y^i)$ on $L$, the inverse metric of $\Pi^*g$ has components
\[
g^{zz}=\vartheta^{-1},\qquad
g^{zi}=-\frac{\gamma}{z\vartheta}(\nabla_hS)^i,\qquad
g^{ij}=z^{-2\gamma}h^{ij}+\frac{\gamma^2}{z^2\vartheta}(\nabla_hS)^i(\nabla_hS)^j,
\]
and $\sqrt{\det\Pi^*g}=z^\beta\vartheta^{1/2}\sqrt{\det h}$.
\end{lemma}

\begin{proof}
One checks directly that the displayed matrix is the inverse of that in Lemma~\ref{lemma:metrica-corno}, whose components are $g_{zz}=1+\gamma^2z^\kappa|y|^2$, $g_{zi}=\gamma z^{2\gamma-1}\partial_iS$, $g_{ij}=z^{2\gamma}h_{ij}$, using $|\nabla_hS|_h^2=h^{ij}\partial_iS\,\partial_jS$ and $\vartheta=g_{zz}-\gamma^2z^{\kappa}|\nabla_hS|_h^2$.
\end{proof}

\begin{lemma}\label{lemma:laplaciano}
Let $\xi\in\R$ and $\chi\in C^\infty(L)$. As $z\to0$, uniformly in $y\in L$,
\begin{align}
\Delta_X(z^\xi)&=\xi\bigl(\xi+\beta-1-\gamma\Delta_hS\bigr)z^{\xi-2}+O(z^{\xi-2+\kappa}),\label{eq:lap-radiale}\\
\Delta_X(z^{\xi+\kappa}\chi)&=(\Delta_h\chi)\,z^{\xi-2}+O(z^{\xi-2+\kappa}).\label{eq:lap-angolare}
\end{align}
\end{lemma}

\begin{proof}
Write $\Delta_Xv=(\sqrt{\det g})^{-1}\partial_A\bigl(\sqrt{\det g}\,g^{AB}\partial_B v\bigr)$, with indices $A,B\in\{z,1,\dots,d-1\}$ and the components of Lemma~\ref{lemma:inversa}. Note that $\vartheta=1+O(z^\kappa)$, that the derivatives of $\vartheta$ with respect to $y$ are $O(z^\kappa)$, and that $\partial_z\vartheta=O(z^{\kappa-1})$.

For $v=z^\xi$ only $g^{zz}$ and $g^{iz}$ contribute:
\[
\begin{split}\frac{1}{z^\beta\vartheta^{1/2}}\partial_z\bigl(z^{\beta}\vartheta^{-1/2}\xi z^{\xi-1}\bigr)&=\xi(\xi+\beta-1)z^{\xi-2}\vartheta^{-1}+\xi z^{\xi-1}\vartheta^{-1/2}\partial_z\vartheta^{-1/2}\\
&=\xi(\xi+\beta-1)z^{\xi-2}+O(z^{\xi-2+\kappa}),
\end{split}
\]
\[
\frac{1}{\vartheta^{1/2}\sqrt{\det h}}\partial_i\Bigl(\sqrt{\det h}\,\vartheta^{-1/2}\bigl(-\gamma(\nabla_hS)^i\bigr)\Bigr)\xi z^{\xi-2}=-\gamma\xi z^{\xi-2}\Delta_hS+O(z^{\xi-2+\kappa}),
\]
which gives~\eqref{eq:lap-radiale}.

For $v=z^{\xi+\kappa}\chi$ the term with $z^{-2\gamma}h^{ij}$ gives
\[
z^{\xi+\kappa-2\gamma}\bigl(\Delta_h\chi+\langle\nabla_h\log\vartheta^{1/2},\nabla_h\chi\rangle_h\bigr)=z^{\xi-2}\Delta_h\chi+O(z^{\xi-2+\kappa}),
\]
since $\kappa-2\gamma=-2$. Each remaining term contains the factor $z^{\xi+\kappa-2}$, by the powers of $z$ in $g^{zz}$, $g^{zi}$ and the second part of $g^{ij}$, multiplied by a function bounded together with its derivatives; these terms are therefore $O(z^{\xi-2+\kappa})$.
\end{proof}

\subsection{Cross-sections contained in a sphere}

\begin{proposition}\label{prop:sfera}
Suppose $|y|=a_L$ for every $y\in L$, and let $0<z_1\le\zeta_0$. The function
\[
\varsigma_L(z):=\int_z^{z_1}\bigl(1+\gamma^2a_L^2t^{\kappa}\bigr)^{1/2}t^{-\beta}\,dt
\]
is harmonic on $X\cap\{0<z<z_1\}$, positive, vanishes on $\{z=z_1\}$, and generates the cone $\mathcal K_R$ with $R=z_1(1+a_L^2z_1^\kappa)^{1/2}$. We have $\varsigma_L(z)\sim z^{1-\beta}/(\beta-1)$ as $z\to0$.
\end{proposition}

\begin{proof}
If $|y|$ is constant, $\nabla_hS=0$, hence $g_{zi}=0$ and $g=(1+\gamma^2a_L^2z^\kappa)dz^2+z^{2\gamma}h$. For a function of $z$ alone we have $\Delta_X\varphi=(z^\beta F^{1/2})^{-1}\bigl(z^\beta F^{-1/2}\varphi'\bigr)'$, with $F:=1+\gamma^2a_L^2z^\kappa$, and $\varsigma_L$ satisfies $z^\beta F^{-1/2}\varsigma_L'=-1$. Moreover, $\rho=z(1+a_L^2z^\kappa)^{1/2}$ is an increasing function of $z$ alone, so $\{z=z_1\}=L_R$ and $\varsigma_L\in\mathcal K_R$; the conclusion follows from~\ref{R:picard}, applicable by~\ref{R:corni} and Remark~\ref{oss:m-maggiore}. The asymptotic equivalence follows from $F\to1$.
\end{proof}

\begin{example}
If $m=d$, the condition $|y|=a_L$ forces $L$ to be the sphere of radius $a_L$, and we recover the rotational cusps of \cite[Prop.~5.12]{V11}. For $m>d$ the class contains nonspherical cross-sections: for instance, the torus $L=\{(y',y'')\in\R^2\times\R^2:|y'|=|y''|=1\}$, with $d=3$, $m=4$ and $a_L=\sqrt2$, gives a horn in $\R^5$ with profile $\varsigma_L(z)\sim z^{1-2\gamma}/(2\gamma-1)$.
\end{example}

\subsection{Barriers and leading coefficient}

Let $\sigma_0:=1-\beta<0$ and $\eta\in(0,\min\{\kappa,\beta-1\})$, and let $\sigma_1:=\sigma_0+\eta<0$. For $\sigma\in\R$ set
\[
W_\sigma(z,y):=z^\sigma\bigl(1+\gamma\sigma S(y)z^\kappa\bigr).
\]

\begin{lemma}\label{lemma:barriere}
For every $\sigma\in\R$, uniformly in $y\in L$,
\[
\Delta_XW_\sigma=\sigma(\sigma+\beta-1)z^{\sigma-2}+O(z^{\sigma-2+\kappa}).
\]
Consequently there exists $z_0\in(0,1]$ such that, on $X\cap\{0<z<z_0\}$, the functions $W_{\sigma_0}$ and $W_{\sigma_1}$ are positive and
\[
\Theta_+:=W_{\sigma_0}+W_{\sigma_1},\qquad \Theta_-:=W_{\sigma_0}-W_{\sigma_1}
\]
satisfy $\Delta_X\Theta_+\le0$ and $\Delta_X\Theta_-\ge0$.
\end{lemma}

\begin{proof}
By~\eqref{eq:lap-radiale} with $\xi=\sigma$ and~\eqref{eq:lap-angolare} with $\xi=\sigma$ and $\chi=\gamma\sigma S$,
\[
\Delta_XW_\sigma=\sigma(\sigma+\beta-1)z^{\sigma-2}-\gamma\sigma(\Delta_hS)z^{\sigma-2}+\gamma\sigma(\Delta_hS)z^{\sigma-2}+O(z^{\sigma-2+\kappa}),
\]
which proves the first assertion. For $\sigma=\sigma_0$ the leading term vanishes: $\Delta_XW_{\sigma_0}=O(z^{\sigma_0-2+\kappa})$. For $\sigma=\sigma_1$ we have $\sigma_1(\sigma_1+\beta-1)=\sigma_1\eta<0$, hence $\Delta_XW_{\sigma_1}\le\frac12\sigma_1\eta z^{\sigma_1-2}$ for small $z$. Since $\eta<\kappa$, we have $z^{\sigma_0-2+\kappa}=o(z^{\sigma_1-2})$, and the inequalities for $\Theta_\pm$ hold for small $z$. Positivity of $W_{\sigma_0}$ and $W_{\sigma_1}$ for small $z$ is immediate.
\end{proof}

\begin{theorem}\label{thm:coefficiente}
Let $X=X_{\gamma,L}$ with $\gamma>1$ and let $\eta\in(0,\min\{2\gamma-2,\beta-1\})$. For every singular profile $G_X$ in (AD) there exists $c>0$ such that, uniformly in $y\in L$,
\[
G_X(\Pi(z,y))=c\,z^{1-\beta}\bigl(1+O(z^\eta)\bigr),\qquad
G_X=c\,\rho^{1-\beta}\bigl(1+O(\rho^\eta)\bigr).
\]
\end{theorem}

\begin{proof}
Let $z_0$ be as in Lemma~\ref{lemma:barriere}, with $z_0\le\zeta_0$. For $0<\varepsilon<z_0$ let $U_\varepsilon$ be the solution, given by~\ref{T:dirichlet}, of the Dirichlet problem on $\overline{\{\varepsilon<z<z_0\}}$, a compact manifold with smooth boundary, with data $U_\varepsilon=W_{\sigma_0}$ on $\{z=\varepsilon\}$ and $U_\varepsilon=0$ on $\{z=z_0\}$. Let $M_0:=\max_{z=z_0}(W_{\sigma_0}+W_{\sigma_1})$, independent of $\varepsilon$. The functions $\Theta_--M_0$ and $\Theta_++M_0$ are subharmonic and superharmonic, respectively; on $\{z=\varepsilon\}$ we have $\Theta_--M_0\le W_{\sigma_0}\le\Theta_++M_0$, and on $\{z=z_0\}$ we have $\Theta_--M_0\le0\le\Theta_++M_0$. By the maximum principle,
\begin{equation}\label{eq:sandwich}
\Theta_--M_0\le U_\varepsilon\le\Theta_++M_0\qquad\text{on }\{\varepsilon\le z\le z_0\}.
\end{equation}
The functions $U_\varepsilon$ are thus uniformly bounded on every set $\{\delta\le z\le z_0\}$; by interior estimates and a diagonal argument, a sequence $\varepsilon_n\to0$ gives $U_{\varepsilon_n}\to U$ locally uniformly on $\{0<z<z_0\}$, with $U$ harmonic. Fix $\delta\in(0,z_0)$, let $M_\delta:=\max_{z=\delta}(\Theta_++M_0)$ and let $e_\delta$ be the harmonic function on $\{\delta<z<z_0\}$, continuous up to the boundary, equal to $1$ on $\{z=\delta\}$ and to $0$ on $\{z=z_0\}$. By the maximum principle $0\le U_\varepsilon\le M_\delta e_\delta$ on $\{\delta\le z\le z_0\}$ for $\varepsilon<\delta$; hence $0\le U\le M_\delta e_\delta$, and $U$ extends continuously with value $0$ on $\{z=z_0\}$. Passing to the limit in~\eqref{eq:sandwich} and dividing by $z^{\sigma_0}$, since $W_{\sigma_0}/z^{\sigma_0}=1+O(z^\kappa)$, $W_{\sigma_1}/z^{\sigma_0}=O(z^\eta)$ and $M_0/z^{\sigma_0}=O(z^{\beta-1})$, we obtain
\begin{equation}\label{eq:U}
U(z,y)=z^{1-\beta}\bigl(1+O(z^\eta)\bigr)
\end{equation}
uniformly in $y$, since $\eta<\min\{\kappa,\beta-1\}$. In particular, $U>0$ near $0$ and $U\in\mathcal P$. By (AD), which holds by~\ref{R:AD} and~\ref{R:corni}, we have $U=qG_X+u_{\mathrm b}$ near $0$, with $q\ge0$ and $u_{\mathrm b}$ bounded; since $U\to+\infty$, we have $q>0$. Hence $G_X=q^{-1}(U-u_{\mathrm b})=q^{-1}z^{1-\beta}(1+O(z^\eta))+O(1)$, and $O(1)=z^{1-\beta}O(z^{\beta-1})$ is absorbed into the remainder. Set $c:=q^{-1}$; it does not depend on $\eta$, being the uniform limit of $z^{\beta-1}G_X$ as $z\to0$. Finally, $\rho=z(1+z^\kappa|y|^2)^{1/2}=z(1+O(z^\kappa))$, hence $z^{1-\beta}=\rho^{1-\beta}(1+O(\rho^\kappa))$ and $z^\eta\asymp\rho^\eta$.
\end{proof}

\begin{remark}
The leading exponent and its uniqueness do not depend on the topology of $L$ or on its position in $\R^m$. The correction $\gamma\sigma Sz^\kappa$ in the barrier $W_\sigma$ is necessary in this construction of barriers: by~\eqref{eq:lap-radiale}, if $|y|$ is not constant, $\Delta_X(z^{\sigma_0})=-\gamma\sigma_0(\Delta_hS)z^{\sigma_0-2}+O(z^{\sigma_0-2+\kappa})$, and $z^{1-\beta}$ is not approximately harmonic to leading order.
\end{remark}

\subsection{Capacity and flux}

\begin{theorem}\label{thm:capacita-esatta}
For $X=X_{\gamma,L}$ with $\gamma>1$ and fixed $0<\zeta\le\zeta_0$,
\[
\Cap^\Pi(\varepsilon;\zeta)\sim(\beta-1)V_L\,\varepsilon^{\beta-1}\qquad(\varepsilon\to0).
\]
\end{theorem}

\begin{proof}
By Proposition~\ref{prop:capacita-corno} it suffices to note that $\int_\varepsilon^\zeta t^{-\beta}dt\sim\varepsilon^{1-\beta}/(\beta-1)$ and that $\int_\varepsilon^\zeta t^{\kappa-\beta}dt=o(\varepsilon^{1-\beta})$, since $\kappa>0$.
\end{proof}

\begin{corollary}\label{cor:flusso}
Let $0<z_1\le\zeta_0$ and let $G$ be positive and harmonic on $X\cap\{0<z<z_1\}$, continuous up to $\{z=z_1\}$ with $G=0$ on $\{z=z_1\}$, and such that $G(z,y)=c\,z^{1-\beta}(1+o(1))$ uniformly in $y$, with $c>0$. Then
\[
\mathfrak F(G)=c\,(\beta-1)V_L,
\]
where $\mathfrak F(G)$ is the flux through the levels $\{z=\varepsilon\}$ with the normal directed toward $0$. In particular, the function $U$ in the proof of Theorem~\ref{thm:coefficiente} has flux $(\beta-1)V_L$, and the unit-flux profile $G^{(1)}:=U/((\beta-1)V_L)$ satisfies
\[
G^{(1)}=\frac{z^{1-\beta}}{(\beta-1)V_L}\bigl(1+O(z^\eta)\bigr).
\]
\end{corollary}

\begin{proof}
Let $e$ be the harmonic function on $\{\varepsilon<z<z_1\}$, continuous up to the boundary, equal to $1$ on $\{z=\varepsilon\}$ and to $0$ on $\{z=z_1\}$. The levels of $z$ are smooth submanifolds, since $\Pi$ is an immersion for $z>0$; the argument of Lemma~\ref{lemma:dirichlet} gives $\partial_\nu e>0$ on $\{z=\varepsilon\}$ and $\int_{z=\varepsilon}\partial_\nu e\,dS=\Cap^\Pi(\varepsilon;z_1)$, and that of Theorem~\ref{thm:capacita-profilo}, applied to $G$, which is smooth up to $\{z=z_1\}$, gives $\int_{z=\varepsilon}G\,\partial_\nu e\,dS=\mathfrak F(G)$. Thus $\mathfrak F(G)/\Cap^\Pi(\varepsilon;z_1)$ is an average of $G$ on $\{z=\varepsilon\}$, equal to $c\,\varepsilon^{1-\beta}(1+o(1))$. By Theorem~\ref{thm:capacita-esatta}, $\mathfrak F(G)=c\,(\beta-1)V_L(1+o(1))$; since $\mathfrak F(G)$ is independent of $\varepsilon$, the equality holds. The function $U$ is nonnegative, not identically zero, and harmonic, hence positive by the strong maximum principle; it satisfies the assumptions with $c=1$ by~\eqref{eq:U}.
\end{proof}

\begin{remark}\label{oss:flusso-generale}
The formula $\mathfrak F(G_X)=c(\beta-1)V_L$ holds for every singular profile in Theorem~\ref{thm:coefficiente}. Fix a level $z_1$ in the interior of the domain of harmonicity of $G_X$ and let $e$ be the equilibrium potential between $\{z=\varepsilon\}$ and $\{z=z_1\}$. Green's second identity gives $\mathfrak F(G_X)=\int_{z=\varepsilon}G_X\partial_\nu e\,dS+\int_{z=z_1}G_X\partial_\nu e\,dS$, with $\nu$ the outward normal to the annulus. The first summand tends to $c(\beta-1)V_L$ by Theorem~\ref{thm:capacita-esatta}; the second is bounded in absolute value by $\max_{z=z_1}|G_X|\,\Cap^\Pi(\varepsilon;z_1)$, which tends to zero.
\end{remark}

\begin{remark}
The coefficient $((\beta-1)V_L)^{-1}$ of the unit-flux profile depends on the cross-section only through its volume. Since $\rho=z(1+O(z^\kappa))$, the same coefficient holds with respect to $\rho$. For a germ merely bi-Lipschitz equivalent to a horn, Lemma~\ref{lemma:trasporto} controls capacities up to constant factors, and Corollary~\ref{cor:bilip} gives comparability but not the existence of a coefficient.
\end{remark}

\section{First angular correction}\label{sec:correzione}

In this section $X=X_{\gamma,L}$ with $\gamma>1$ and $d\ge3$, and we use the notation of Section~\ref{sec:corni-esatti}. Set $T_L:=|\nabla_hS|_h^2$, so that $\mathcal N=2S-T_L$, and, for $f\in C(L)$, $\langle f\rangle:=V_L^{-1}\int_Lf\,d\vol_h$. Note that
\[
\sigma_0+\kappa=\gamma(3-d)-1<0 .
\]

\subsection{Expansion of the Laplacian}

For $\xi\in\R$ and $\chi\in C^\infty(L)$ set
\begin{multline*}
\mathfrak L_\xi\chi:=\xi(\xi+\beta-1)\chi-\gamma(\xi+\beta-1)\langle\nabla_hS,\nabla_h\chi\rangle_h-\gamma\xi\,\mathrm{div}_h(\chi\nabla_hS)\\
+\gamma^2\,\mathrm{div}_h\bigl(\langle\nabla_hS,\nabla_h\chi\rangle_h\nabla_hS\bigr)+\tfrac{\gamma^2}{2}\langle\nabla_h\mathcal N,\nabla_h\chi\rangle_h .
\end{multline*}

\begin{lemma}\label{lemma:sviluppo}
Let $\xi\in\R$ and $\chi\in C^\infty(L)$. There exists a function $\Psi_{\xi,\chi}$, smooth on $[0,w_0)\times L$ for some $w_0>0$, with $\Psi_{\xi,\chi}(0,\cdot)=\mathfrak L_\xi\chi$, such that
\[
\Delta_X(z^\xi\chi)=z^{\xi-2-\kappa}\Delta_h\chi+z^{\xi-2}\,\Psi_{\xi,\chi}(z^\kappa,y)\qquad(0<z^\kappa<w_0).
\]
In particular, for every integer $n\ge1$, $\Delta_X(z^\xi\chi)=z^{\xi-2-\kappa}\Delta_h\chi+z^{\xi-2}\sum_{i<n}z^{i\kappa}\partial_w^i\Psi_{\xi,\chi}(0,y)/i!+O(z^{\xi-2+n\kappa})$ uniformly in $y$.
\end{lemma}

\begin{proof}
By Lemma~\ref{lemma:inversa}, $\vartheta$, $g^{zz}=\vartheta^{-1}$, $z\,g^{zi}=-\gamma\vartheta^{-1}(\nabla_hS)^i$ and $z^2g^{ij}-z^{-\kappa}h^{ij}=\gamma^2\vartheta^{-1}(\nabla_hS)^i(\nabla_hS)^j$ are smooth functions of $(z^\kappa,y)$, and $\sqrt{\det g}=z^\beta\vartheta^{1/2}\sqrt{\det h}$. The operator $z\partial_z$ maps a smooth function of $(z^\kappa,y)$ to a function of the same type, vanishing at $z^\kappa=0$. Writing $\Delta_X=(\sqrt{\det g})^{-1}\partial_A(\sqrt{\det g}\,g^{AB}\partial_B)$ and applying it to $z^\xi\chi$, each summand is $z^{\xi-2}$ times a smooth function of $(z^\kappa,y)$, except for the one containing $z^{-2\gamma}h^{ij}$, which equals
\[
z^{\xi-2-\kappa}\bigl(\Delta_h\chi+\langle\nabla_h\log\vartheta^{1/2},\nabla_h\chi\rangle_h\bigr),
\qquad
\nabla_h\log\vartheta^{1/2}=\tfrac12\gamma^2z^\kappa\vartheta^{-1}\nabla_h\mathcal N .
\]
This proves the first assertion. To compute $\Psi_{\xi,\chi}(0,\cdot)$, set $\vartheta=1$ and discard the terms containing $z\partial_z\vartheta$: the block $g^{zz}$ gives $\xi(\xi+\beta-1)\chi$; the term with $g^{zi}$ differentiated in $z$ gives $-\gamma(\xi+\beta-1)\langle\nabla_hS,\nabla_h\chi\rangle_h$; the term with $g^{iz}$ differentiated in $y$ gives $-\gamma\xi\,\mathrm{div}_h(\chi\nabla_hS)$; the second part of $g^{ij}$ gives $\gamma^2\mathrm{div}_h(\langle\nabla_hS,\nabla_h\chi\rangle_h\nabla_hS)$; the derivative of $\log\vartheta^{1/2}$ gives $\frac{\gamma^2}{2}\langle\nabla_h\mathcal N,\nabla_h\chi\rangle_h$.
\end{proof}

For $\chi=1$ and for $\xi$ replaced by $\xi+\kappa$, Lemma~\ref{lemma:sviluppo} recovers formulas~\eqref{eq:lap-radiale} and~\eqref{eq:lap-angolare}.

\begin{lemma}[flux]\label{lemma:flusso-coord}
For $v$ smooth on $\{0<z<1\}$ set
\[
\mathfrak J_v(z):=\int_Lz^\beta\vartheta^{-1/2}\bigl(\partial_zv-\gamma z^{-1}\langle\nabla_hS,\nabla_hv\rangle_h\bigr)\,d\vol_h .
\]
Then $\int_L(\Delta_Xv)\,z^\beta\vartheta^{1/2}\,d\vol_h=\mathfrak J_v'(z)$. In particular, if $v$ is harmonic, $\mathfrak J_v$ is constant and $\mathfrak F(v)=-\mathfrak J_v$.
\end{lemma}

\begin{proof}
Integrate over $L$ the identity $\sqrt{\det g}\,\Delta_Xv=\partial_A(\sqrt{\det g}\,g^{AB}\partial_Bv)$. The terms involving $\partial_i$ are divergences of vector fields on $(L,h)$ and have zero integral; the term involving $\partial_z$ gives $\mathfrak J_v'$, since $\vartheta^{1/2}g^{zz}=\vartheta^{-1/2}$ and $\vartheta^{1/2}g^{zi}\partial_iv=-\gamma z^{-1}\vartheta^{-1/2}\langle\nabla_hS,\nabla_hv\rangle_h$.
\end{proof}

\subsection{The correction coefficient}

\begin{theorem}\label{thm:correzione}
Let $d\ge3$, $\gamma>1$ and $0<\eta<\min\{\kappa,-(\sigma_0+\kappa)\}$. Set
\[
f_1:=\gamma\sigma_0\bigl(S-\langle S\rangle\bigr)+\frac{\gamma^2\sigma_0}{\sigma_0+\kappa}\Bigl(\langle S\rangle+\tfrac12\langle T_L\rangle\Bigr).
\]
For every singular profile $G_X$ in (AD) there exists $c>0$ such that, uniformly in $y\in L$,
\[
G_X(\Pi(z,y))=c\,z^{\sigma_0}\bigl(1+z^\kappa f_1(y)+O(z^{\kappa+\eta})\bigr).
\]
There also exists a positive harmonic function $G^{(1)}$ on $\{0<z<z_0\}$, vanishing on $\{z=z_0\}$, with $\mathfrak F(G^{(1)})=1$ and
\[
G^{(1)}=\frac{z^{\sigma_0}}{(\beta-1)V_L}\bigl(1+z^\kappa f_1+O(z^{\kappa+\eta})\bigr).
\]
\end{theorem}

\begin{proof}
\emph{First corrector.} Let $c_0\in\R$ be determined below and let $P_1:=z^{\sigma_0}\bigl(1+z^\kappa(\gamma\sigma_0S+c_0)\bigr)$. By Lemma~\ref{lemma:sviluppo}, the coefficient of $z^{\sigma_0-2}$ in $\Delta_XP_1$ is $\mathfrak L_{\sigma_0}1+\Delta_h(\gamma\sigma_0S)=-\gamma\sigma_0\Delta_hS+\gamma\sigma_0\Delta_hS=0$, since $\sigma_0(\sigma_0+\beta-1)=0$. Hence
\[
\Delta_XP_1=z^{\sigma_0+\kappa-2}R_1+O(z^{\sigma_0+2\kappa-2}),
\]
with $R_1\in C^\infty(L)$.

\emph{Solvability.} By Lemma~\ref{lemma:flusso-coord} and $\sigma_0+\beta=1$, $\mathfrak J_{P_1}$ is a smooth function of $z^\kappa$: $\mathfrak J_{P_1}=\mathfrak J^{(0)}+\mathfrak J^{(1)}z^\kappa+O(z^{2\kappa})$. Since $\sigma_0+\beta=1$, Lemma~\ref{lemma:flusso-coord} gives
\[
z^{\kappa-1}\int_LR_1\,d\vol_h+O(z^{2\kappa-1})=\kappa\,\mathfrak J^{(1)}z^{\kappa-1}+O(z^{2\kappa-1}),
\]
hence $\int_LR_1\,d\vol_h=\kappa\,\mathfrak J^{(1)}$. We compute $\mathfrak J^{(1)}$ using $\vartheta^{-1/2}=1-\frac{\gamma^2}{2}z^\kappa\mathcal N+O(z^{2\kappa})$, $\partial_zP_1=\sigma_0z^{\sigma_0-1}+(\sigma_0+\kappa)z^{\sigma_0+\kappa-1}(\gamma\sigma_0S+c_0)$ and $\langle\nabla_hS,\nabla_hP_1\rangle_h=\gamma\sigma_0z^{\sigma_0+\kappa}T_L$:
\[
\frac{\mathfrak J^{(1)}}{V_L}=(\sigma_0+\kappa)\bigl(\gamma\sigma_0\langle S\rangle+c_0\bigr)-\gamma^2\sigma_0\langle T_L\rangle-\frac{\gamma^2\sigma_0}{2}\langle\mathcal N\rangle
=(\sigma_0+\kappa)\bigl(\gamma\sigma_0\langle S\rangle+c_0\bigr)-\gamma^2\sigma_0\Bigl(\langle S\rangle+\tfrac12\langle T_L\rangle\Bigr).
\]
Since $\sigma_0+\kappa\ne0$, there exists a unique $c_0$ with $\mathfrak J^{(1)}=0$, and with this choice $\gamma\sigma_0S+c_0=f_1$ and $\int_LR_1\,d\vol_h=0$.

\emph{Second corrector.} Since $L$ is compact and connected and $\int_LR_1\,d\vol_h=0$, by~\ref{T:poisson} there exists $f_2\in C^\infty(L)$ with $\Delta_hf_2=-R_1$. Let $P_2:=P_1+z^{\sigma_0+2\kappa}f_2$. By Lemma~\ref{lemma:sviluppo} with $\xi=\sigma_0+2\kappa$, $\Delta_X(z^{\sigma_0+2\kappa}f_2)=z^{\sigma_0+\kappa-2}\Delta_hf_2+O(z^{\sigma_0+2\kappa-2})$; hence $\Delta_XP_2=O(z^{\sigma_0+2\kappa-2})$.

\emph{Barriers.} Let $\sigma_2:=\sigma_0+\kappa+\eta<0$. By Lemma~\ref{lemma:barriere}, $\Delta_XW_{\sigma_2}=\sigma_2(\kappa+\eta)z^{\sigma_2-2}+O(z^{\sigma_2-2+\kappa})$, with $\sigma_2(\kappa+\eta)<0$. Since $\eta<\kappa$, we have $z^{\sigma_0+2\kappa-2}=o(z^{\sigma_2-2})$. Thus there exists $z_0\in(0,\zeta_0]$ such that, on $\{0<z<z_0\}$, $W_{\sigma_2}>0$, $P_2>0$, $\Delta_X(P_2+W_{\sigma_2})\le0$ and $\Delta_X(P_2-W_{\sigma_2})\ge0$.

\emph{Construction.} Let $M_1:=\max_{z=z_0}(|P_2|+W_{\sigma_2})$. For $0<\varepsilon<z_0$ let $U_\varepsilon$ be harmonic on $\{\varepsilon<z<z_0\}$, with $U_\varepsilon=P_2$ on $\{z=\varepsilon\}$ and $U_\varepsilon=0$ on $\{z=z_0\}$. The maximum principle gives $P_2-W_{\sigma_2}-M_1\le U_\varepsilon\le P_2+W_{\sigma_2}+M_1$. As in the proof of Theorem~\ref{thm:coefficiente}, we extract a limit $U$, harmonic, nonnegative, vanishing continuously on $\{z=z_0\}$, with
\[
U=P_2+O(z^{\sigma_2})+O(1)=z^{\sigma_0}\bigl(1+z^\kappa f_1+O(z^{\kappa+\eta})\bigr),
\]
since $z^{2\kappa}=O(z^{\kappa+\eta})$ and $O(1)=z^{\sigma_0}O(z^{\beta-1})$ with $\beta-1=-\sigma_0>\kappa+\eta$. By the strong maximum principle $U>0$. By Corollary~\ref{cor:flusso}, $\mathfrak F(U)=(\beta-1)V_L$, and set $G^{(1)}:=U/((\beta-1)V_L)$. Finally, by (AD), $U=qG_X+u_{\mathrm b}$ near $0$, with $u_{\mathrm b}$ bounded; since $U\to+\infty$ uniformly, $q>0$, and $G_X=c(U-u_{\mathrm b})$ with $c=q^{-1}$, and $u_{\mathrm b}=z^{\sigma_0}O(z^{\beta-1})$.
\end{proof}

\begin{corollary}\label{cor:differenza-angolare}
Under the assumptions of Theorem~\ref{thm:correzione}, for $y_1,y_2\in L$,
\[
G_X(\Pi(z,y_1))-G_X(\Pi(z,y_2))=c\,\gamma\sigma_0\bigl(S(y_1)-S(y_2)\bigr)z^{\sigma_0+\kappa}+O(z^{\sigma_0+\kappa+\eta}).
\]
If $|y|$ is not constant on $L$, the difference between two points at the same height $z$ can diverge, since $\sigma_0+\kappa<0$. Points at the same height are not generally on the same link $L_r$; since $z=\rho(1-S\rho^\kappa+O(\rho^{2\kappa}))$, in the variable $\rho$ the corrector becomes $f_1-\sigma_0S$.
\end{corollary}

\begin{example}[toroidal horn]\label{es:toro-correzione}
Let $L$ be the torus of \cite[\S6.1]{V11}, $y(\theta_1,\theta_2)=\bigl((2+\cos\theta_2)\cos\theta_1,(2+\cos\theta_2)\sin\theta_1,\sin\theta_2\bigr)$, with $d=3$ and $\gamma=2$. Then $h=(2+\cos\theta_2)^2d\theta_1^2+d\theta_2^2$, $\beta=4$, $\sigma_0=-3$, $\kappa=2$, $V_L=8\pi^2$, $S=\frac52+2\cos\theta_2$, $T_L=4\sin^2\theta_2$, $\langle S\rangle=3$, $\langle T_L\rangle=2$, and
\[
f_1=51-12\cos\theta_2,\qquad
G^{(1)}=\frac{z^{-3}}{24\pi^2}\bigl(1+(51-12\cos\theta_2)z^2+O(z^{2+\eta})\bigr),\quad0<\eta<1 .
\]
The residual in the proof is, with $c_2:=\cos\theta_2$,
\[
R_1=\frac{672c_2^3+1008c_2^2-258c_2-504}{2+c_2},
\]
which has zero average with respect to $(2+c_2)\,d\theta_1d\theta_2$.
\end{example}

\begin{proposition}\label{prop:non-invarianza}
Let $d\ge3$, $\gamma>1$ rational, $\alpha>0$, $b_L\in\R^d$ with $0<|b_L|<\alpha$, $L:=\alpha S^{d-1}$ and $L_{b_L}:=L+b_L$. The horns $X_{\gamma,L}$ and $X_{\gamma,L_{b_L}}$ are bi-Lipschitz equivalent with respect to the intrinsic metrics; the coefficient $f_1$ is constant for $X_{\gamma,L}$ and nonconstant for $X_{\gamma,L_{b_L}}$.
\end{proposition}

\begin{proof}
For any horn, from $\Pi^*g=dz^2+|\gamma z^{\gamma-1}y\,dz+z^\gamma dy|^2$ we obtain $\Pi^*g\le(1+2\gamma^2a_L^2)dz^2+2z^{2\gamma}h$ and $z^{2\gamma}h\le2\,\Pi^*g+2\gamma^2a_L^2dz^2\le(2+2\gamma^2a_L^2)\Pi^*g$, for $z\le1$; hence $\Pi^*g$ is uniformly comparable to $dz^2+z^{2\gamma}h$. The cross-sections $L$ and $L_{b_L}$ are isometric under $y\mapsto y+b_L$, so the map $\Pi(z,y)\mapsto\Pi(z,y+b_L)$ is bi-Lipschitz for the induced metrics, and for the intrinsic distances. On $L$ we have $S\equiv\alpha^2/2$, and $f_1$ is constant. On $L_{b_L}$, writing $y=\alpha\upsilon+b_L$ with $\upsilon\in S^{d-1}$, we have $S=\frac12(\alpha^2+|b_L|^2)+\alpha\langle b_L,\upsilon\rangle$, and the nonconstant part of $f_1$ is $\gamma\sigma_0\alpha\langle b_L,\upsilon\rangle\ne0$.
\end{proof}

\begin{remark}
By Proposition~\ref{prop:non-invarianza}, the term $z^\kappa f_1$ at fixed height $z$ is not an invariant of the bi-Lipschitz equivalence class; it depends on the embedding of the cross-section. The leading exponent of the profile, on the other hand, is preserved in the class of germs bi-Lipschitz equivalent with respect to the intrinsic metrics to the horns under consideration, by Corollary~\ref{cor:bilip}.
\end{remark}

\section{Multi-rate collapsing horns}\label{sec:collasso}

\subsection{The family}

Let $p\ge1$ and, for $i=1,\dots,p$, let $n_i\ge2$, $\alpha_i>0$, let $Y_i\subset\alpha_iS^{n_i-1}\subset\R^{n_i}$ be compact connected semialgebraic analytic submanifolds of dimension $d_i\ge1$, with induced metric $h_i$, and let $\gamma_i\ge1$ be rational. Set
\[
X_{\boldsymbol\gamma}:=\bigl\{(z^{\gamma_1}y_1,\dots,z^{\gamma_p}y_p,z):\ y_i\in Y_i,\ z\ge0\bigr\}\subset\R^{n_1+\dots+n_p+1},
\]
of dimension $d=1+\sum_id_i$, and
\[
\beta:=\sum_id_i\gamma_i,\qquad
\Xi(z):=\Bigl(1+\sum_i\gamma_i^2\alpha_i^2z^{2\gamma_i-2}\Bigr)^{1/2},\qquad
\Xi_0:=\lim_{z\to0}\Xi(z),\qquad
V_Y:=\prod_i\vol_{h_i}(Y_i).
\]
For $p=1$ we recover horns with cross-section contained in a sphere, and $\beta$ coincides with that in Section~\ref{sec:corni}. If the exponents are not all equal, $X_{\boldsymbol\gamma}$ is not a horn in the sense of Definition~\ref{def:corno}: the link factors contract at different rates. The link of $X_{\boldsymbol\gamma}$ is diffeomorphic to $Y_1\times\dots\times Y_p$ and hence connected.

\begin{lemma}\label{lemma:metrica-multipla}
The induced metric is $g=\Xi(z)^2dz^2+\sum_iz^{2\gamma_i}h_i$.\linebreak The measure is $\Xi z^\beta\,dz\,d\vol_{h_1}\cdots d\vol_{h_p}$, and
\[
\Delta_Xu=\frac{1}{\Xi z^\beta}\partial_z\Bigl(\frac{z^\beta}{\Xi}\partial_zu\Bigr)+\sum_iz^{-2\gamma_i}\Delta_{h_i}u .
\]
The distance from the vertex $\rho=z(1+\sum_i\alpha_i^2z^{2\gamma_i-2})^{1/2}$ depends only on $z$, and $\rho/z\to\Xi_0$. The set $X_{\boldsymbol\gamma}$ satisfies the assumptions of \S\ref{sec:prelim}.
\end{lemma}

\begin{proof}
Since $|y_i|=\alpha_i$, we have $\langle y_i,dy_i\rangle=0$, and the mixed term in the metric vanishes; the rest follows by direct computation. We have $\Xi_0^2=1+\sum_{\gamma_i=1}\alpha_i^2=\lim(\rho/z)^2$. The gradient of $\rho$ is parallel to $\partial_z$ and $|\nabla_X\rho|=\rho'(z)/\Xi(z)\to1$, which gives the Whitney condition at $0$; for subanalyticity, let $q$ be a common denominator of the $\gamma_i$ and write $z=s^q$: the graph of the resulting parametrization on $[0,1]\times\prod_iY_i$ is compact and semianalytic, its projection is subanalytic near $0$, the image is closed near $0$ by compactness, and for $z>0$ the parametrization is an analytic embedding.
\end{proof}

For fixed $0<z_1\le1$ sufficiently small that $\rho(z_1)\le r_0$, set
\[
\varsigma_{\boldsymbol\gamma}(z):=\int_z^{z_1}\Xi(t)\,t^{-\beta}\,dt .
\]
The function $\varsigma_{\boldsymbol\gamma}$ is harmonic, positive on $\{0<z<z_1\}$, and vanishes at $z=z_1$; if $\beta>1$, $\varsigma_{\boldsymbol\gamma}(z)\sim\Xi_0z^{1-\beta}/(\beta-1)$, and if $\beta=1$, which forces $p=1$, $d=2$ and $\gamma_1=1$, we have $\varsigma_{\boldsymbol\gamma}(z)=\Xi_0\log(z_1/z)$.

\subsection{Decomposition}

\begin{lemma}\label{lemma:edo}
Let $\mathsf q$ be continuous on $[w_1,\infty)$, $w_1>0$, with $\mathsf q(w)\ge c_{\mathsf q}w^{-2}$ and $c_{\mathsf q}>0$. If $\mathsf u''=\mathsf q\,\mathsf u$ on $[w_1,\infty)$ and $\mathsf u(w)=O(w)$, then $\mathsf u$ is bounded, $\mathsf u(w)\to0$, and for every $\varkappa>0$ with $\varkappa(\varkappa+1)\le c_{\mathsf q}$
\[
|\mathsf u(w)|\le|\mathsf u(w_1)|\,(w_1/w)^{\varkappa}\qquad(w\ge w_1).
\]
\end{lemma}

\begin{proof}
The function $\mathsf u^2$ is convex, since $(\mathsf u^2)''=2\mathsf u'^2+2\mathsf q\mathsf u^2\ge0$. If $\mathsf u$ were unbounded, on a half-line $[w_2,\infty)$ we would have $\mathsf u\ne0$ of constant sign and $|\mathsf u|$ increasing; comparison with the solution of $\mathsf v''=c_{\mathsf q}w^{-2}\mathsf v$ with the same data at $w_2$, where $w_2$ is chosen with $|\mathsf u(w_2)|>0$ and $(|\mathsf u|)'(w_2)>0$, so that the coefficient of the growing Euler solution is strictly positive and the comparison follows from the integral equation for the difference, with zero initial data and nonnegative forcing, would give $|\mathsf u(w)|\ge cw^{a}$, with $a(a-1)=c_{\mathsf q}$, $a>1$, contradicting $\mathsf u=O(w)$. Thus $\mathsf u^2$ is convex and bounded, hence nonincreasing, and $\mathsf u$ has a limit $\ell_\infty$. If $\ell_\infty\ne0$, changing sign if necessary gives $\mathsf u\ge\ell_\infty/2>0$ for large $w$. Then $\mathsf u''=\mathsf q\mathsf u>0$ and $\mathsf u'\le0$ because $\mathsf u^2$ is nonincreasing. The limit of $\mathsf u'$ is zero: a negative limit would contradict boundedness of $\mathsf u$. Integrating $\mathsf u''=\mathsf q\mathsf u$ therefore gives
\[
 -\mathsf u'(w)=\int_w^\infty\mathsf q(t)\mathsf u(t)\dd t
 \ge\frac{c_{\mathsf q}\ell_\infty}{2w},
\]
which contradicts convergence of $\mathsf u$. Finally, $\mathsf b_{\rm ODE}(w):=|\mathsf u(w_1)|(w_1/w)^\varkappa$ satisfies $\mathsf b_{\rm ODE}''\le\mathsf q\mathsf b_{\rm ODE}$; the maximum principle for the operator $\mathsf v\mapsto\mathsf v''-\mathsf q\mathsf v$ on $[w_1,w_3]$, applied to $\pm\mathsf u-\mathsf b_{\rm ODE}-\varepsilon$ with $w_3$ so large that $|\mathsf u(w_3)|\le\varepsilon$, which is possible since $\mathsf u\to0$, and the limits $w_3\to\infty$, $\varepsilon\to0$, give the estimate.
\end{proof}

\begin{theorem}\label{thm:multipla}
Let $u$ be harmonic on $X_{\boldsymbol\gamma}\cap\{0<z<z_1\}$, positive or satisfying $|u|\le C(1+\varsigma_{\boldsymbol\gamma})$. There exist $A\in\R$, $B\in\R$ and $\varkappa>0$, where $\varkappa$ depends only on $X_{\boldsymbol\gamma}$, such that, uniformly on the link,
\[
u=A+B\,\varsigma_{\boldsymbol\gamma}+O\bigl(\varsigma_{\boldsymbol\gamma}^{-\varkappa}\bigr)\qquad(z\to0),
\]
with $B\ge0$ if $u>0$. Consequently, for $X_{\boldsymbol\gamma}$: bounded harmonic functions extend continuously to $0$; the Picard principle holds, with $\mathcal K_R=\{t\varsigma_{\boldsymbol\gamma}:t>0\}$ for $R=\rho(z_1)$; the Bôcher property holds.
\end{theorem}

\begin{proof}
On each closed manifold $Y_i$ the operator $(1-\Delta_{h_i})^{-1}$, defined through the coercive form $\int_{Y_i}(fg+\langle\nabla f,\nabla g\rangle_{h_i})$ and the Lax--Milgram theorem \cite[Theorem~5.8]{GT}, maps $L^2(Y_i)$ boundedly into $W^{1,2}(Y_i)$; it is therefore compact on $L^2(Y_i)$ by the compact embedding used in~\ref{T:dirichlet}, and it is self-adjoint, positive and injective. Maximizing successively its quadratic form on the unit sphere of $L^2(Y_i)$ and on the orthogonal complements of the eigenvectors already found yields an orthonormal basis of eigenvectors: every positive maximum is attained by compactness, the eigenvalues tend to zero if infinitely many, and a nonzero remaining complement would contradict the bound by all preceding maxima. The eigenfunctions are smooth by \cite[Theorem~8.10]{GT} iterated in charts, and by the energy identity and connectedness the kernel of $\Delta_{h_i}$ consists of the constants. Finite sums of simple tensors are dense in $L^2(\prod_iY_i)$, so the products of these bases are complete. Let $\mathsf E_{i,j}$, $j\ge0$, be these orthonormal eigenfunctions of $-\Delta_{h_i}$ with eigenvalues $\lambda_{i,j}$, $\lambda_{i,0}=0<\lambda_{i,1}$ (by connectedness of $Y_i$). For a multi-index $\mathbf j=(j_1,\dots,j_p)$ set $\mathsf E_{\mathbf j}:=\prod_i\mathsf E_{i,j_i}$ and $\lambda_{\mathbf j}:=\sum_i\lambda_{i,j_i}$, and write $u_{\mathbf j}(z):=\int u\,\mathsf E_{\mathbf j}$ on the link, with respect to the product measure. In the variable $w:=\varsigma_{\boldsymbol\gamma}(z)$, by Lemma~\ref{lemma:metrica-multipla}, we have $z^\beta\Xi^{-1}\partial_z=-\partial_w$ and hence
\[
\partial_w^2u_{\mathbf j}=\mathsf q_{\mathbf j}\,u_{\mathbf j},\qquad \mathsf q_{\mathbf j}:=\sum_i\lambda_{i,j_i}z^{2\beta-2\gamma_i}\ge\lambda_{\mathbf j}z^{2\beta-2},
\]
since $\gamma_i\ge1$ and $z<1$. If $\beta>1$, $w\asymp z^{1-\beta}$ and $z^{2\beta-2}\ge cw^{-2}$; if $\beta=1$, $\mathsf q_{\mathbf j}=\lambda_{\mathbf j}$. In both cases, for every $\mathbf j\ne\mathbf0$ we have $\mathsf q_{\mathbf j}\ge c\lambda_{\min}w^{-2}$ for $w\ge w_1$, with $\lambda_{\min}:=\min_i\lambda_{i,1}>0$ and $c,w_1$ independent of $\mathbf j$. With $\mathsf E_{i,0}:=\vol_{h_i}(Y_i)^{-1/2}$, the constant coefficient satisfies $\partial_w^2u_{\mathbf 0}=0$, hence $u_{\mathbf 0}=A_0+B_0w$, and we set $A:=V_Y^{-1/2}A_0$, $B:=V_Y^{-1/2}B_0$. If $u>0$, $|u_{\mathbf j}|\le\|\mathsf E_{\mathbf j}\|_\infty V_Y^{1/2}u_{\mathbf 0}=O(1+w)$ and $B_0\ge0$; in the other case the same estimate follows from the assumption. Lemma~\ref{lemma:edo} gives $|u_{\mathbf j}(w)|\le|u_{\mathbf j}(w_1)|(w_1/w)^\varkappa$ with $\varkappa$ independent of $\mathbf j$. Let $Y:=\prod_iY_i$ with the product metric, so that $\Delta_Y=\sum_i\Delta_{h_i}$ and $(1-\Delta_Y)\mathsf E_{\mathbf j}=(1+\lambda_{\mathbf j})\mathsf E_{\mathbf j}$, and set $v(w,\cdot):=u(w,\cdot)-V_Y^{-1/2}u_{\mathbf 0}(w)$. Fix an integer $M$ with $4M>\dim Y$. By~\ref{T:sobolev}, by Parseval's identity for the complete orthonormal system $(\mathsf E_{\mathbf j})$ and by the estimate above, for $w\ge w_1$
\[
\sup_Y|v(w,\cdot)|\le C\Bigl(\sum_{\mathbf j\ne\mathbf 0}(1+\lambda_{\mathbf j})^{2M}|u_{\mathbf j}(w)|^2\Bigr)^{1/2}\le C\Bigl(\frac{w_1}{w}\Bigr)^{\varkappa}\bigl\|(1-\Delta_Y)^Mv(w_1,\cdot)\bigr\|_{L^2(Y)},
\]
and the last norm is finite since $u$ is smooth on the cross-section $\{w=w_1\}$. Hence
\[
\bigl|u-V_Y^{-1/2}(A_0+B_0w)\bigr|\le Cw^{-\varkappa}
\]
uniformly on the link. If $u$ is bounded, $B_0=0$ and $u$ has a limit. If $u\in\mathcal K_R$, each $u_{\mathbf j}$, $\mathbf j\ne\mathbf 0$, vanishes at $w=0$, tends to $0$ and satisfies $\partial_w^2u_{\mathbf j}=\mathsf q_{\mathbf j}u_{\mathbf j}$ with $\mathsf q_{\mathbf j}>0$; by the maximum principle $u_{\mathbf j}\equiv0$, and $u=Bw$. If $u_1,u_2>0$ diverge, with coefficients $B_1,B_2>0$, the function $B_2u_1-B_1u_2$ has a limit at $0$.
\end{proof}

The proof does not use conditions (H$_b$) or (HS) of \cite{V11}: it relies only on exact separation of variables.

\begin{proposition}\label{prop:capacita-multipla}
For $0<\varepsilon<z_1$, $\Cap\bigl(\{0<z\le\varepsilon\};\{0<z<z_1\}\bigr)=V_Y/\varsigma_{\boldsymbol\gamma}(\varepsilon)$, and $\mathfrak F(\varsigma_{\boldsymbol\gamma})=V_Y$. The unit-flux profile $\varsigma_{\boldsymbol\gamma}/V_Y$ satisfies, if $\beta>1$,
\[
\frac{\varsigma_{\boldsymbol\gamma}}{V_Y}\sim\frac{\Xi_0}{(\beta-1)V_Y}z^{1-\beta}\sim\frac{\Xi_0^{\beta}}{(\beta-1)V_Y}\rho^{1-\beta}.
\]
\end{proposition}

\begin{proof}
The function $\varsigma_{\boldsymbol\gamma}/\varsigma_{\boldsymbol\gamma}(\varepsilon)$ is the equilibrium potential of the condenser, and the argument of Lemma~\ref{lemma:dirichlet} identifies capacity with its flux. On the surface $\{z=\varepsilon\}$ the unit normal directed toward $0$ is $-\Xi^{-1}\partial_z$, the normal derivative of $\varsigma_{\boldsymbol\gamma}$ is $\varepsilon^{-\beta}$ and the area element is $\varepsilon^\beta d\vol_{h_1}\cdots d\vol_{h_p}$; the flux is thus $V_Y$. The last asymptotic equivalence follows from $\rho\sim\Xi_0z$.
\end{proof}

\subsection{Germs with the same leading-order term}

\begin{proposition}\label{prop:x24-x33}
In $\R^5=\R^2\times\R^2\times\R$ let
\[
X_{2,4}:=\{|x'|^2=z^4,\ |x''|^2=z^8,\ z\ge0\},\qquad
X_{3,3}:=\{|x'|^2=z^6,\ |x''|^2=z^6,\ z\ge0\}.
\]
These are semialgebraic germs of dimension $3$ with link $S^1\times S^1$, with $\beta=6$, $V_Y=4\pi^2$, $\Xi_0=1$, and both unit-flux profiles satisfy $G^{(1)}\sim\rho^{-5}/(20\pi^2)$. The two germs are not bi-Lipschitz equivalent with respect to the intrinsic metrics. Their profiles differ at second order:
\[
\varsigma_{2,4}=\tfrac15\rho^{-5}+\tfrac76\rho^{-3}+O(\rho^{-1}),\qquad
\varsigma_{3,3}=\tfrac15\rho^{-5}+O(\rho^{-1}).
\]
\end{proposition}

\begin{proof}
They are cases of $X_{\boldsymbol\gamma}$ with $p=2$, $Y_1=Y_2=S^1$, $\boldsymbol\gamma=(2,4)$ and $(3,3)$; the first assertions follow from Proposition~\ref{prop:capacita-multipla}. The expansions follow from $\Xi_{2,4}=(1+4z^2+16z^6)^{1/2}=1+2z^2-2z^4+O(z^6)$, $\Xi_{3,3}=(1+18z^4)^{1/2}=1+9z^4+O(z^8)$, from $\rho_{2,4}=z(1+z^2+z^6)^{1/2}$, $\rho_{3,3}=z(1+2z^4)^{1/2}$ and inversion of the series.

Let $\Psi$ be a $K_\Psi$-bi-Lipschitz homeomorphism with respect to the intrinsic metrics from a neighborhood of $0$ in $X_{2,4}$ onto a neighborhood of $0$ in $X_{3,3}$, with $\Psi(0)=0$. For small $r$ let $\varpi_r\subset X_{2,4}$ be the curve obtained by rotating the second factor at $z=r$ and keeping the first fixed; it has length $2\pi r^4$ and represents a nonzero element of $\pi_1$ of the punctured germ, isomorphic to $\mathbb Z^2$. Using the cofinal product collars $(0,\delta)\times\mathbb T^2$, $\Psi$ induces an isomorphism of the local fundamental groups; hence the curve $\Psi(\varpi_r)$ is nontrivial in $\pi_1$, and by~\ref{R:distanza} its points have $z\ge c\,r$, with $c>0$ depending on $K_\Psi$. In the metric $\Xi^2dz^2+z^6(d\theta_1^2+d\theta_2^2)$ of $X_{3,3}$, a nontrivial curve contained in $\{z\ge cr\}$ has nonzero winding number in at least one of the angles, and hence length at least $2\pi(cr)^3$. On the other hand its length is at most $2\pi K_\Psi r^4$, which is impossible for small $r$.
\end{proof}

\begin{remark}
The leading exponent of the profile and the volume exponent, equal to $\beta+1=7$ for both, do not distinguish the bi-Lipschitz equivalence classes; subsequent terms of the profile and scales of the cycles of the link provide further information.
\end{remark}

\begin{proposition}\label{prop:log-multipla}
Since the $\gamma_i$ are rational, $\Xi$ is the sum of a convergent Puiseux series in $z$ near $0$; let $[\Xi]_{\beta-1}$ be the coefficient of $z^{\beta-1}$. The expansion of $\varsigma_{\boldsymbol\gamma}$ as $z\to0$ contains a term $c_{\log}\log(1/z)$ with $c_{\log}\ne0$ if and only if $[\Xi]_{\beta-1}\ne0$, and in this case $c_{\log}=[\Xi]_{\beta-1}$. The coefficient is the same with respect to $\log(1/\rho)$.
\end{proposition}

\begin{proof}
Integrate the series of $\Xi(t)t^{-\beta}$ term by term on $[z,\delta]$, with $\delta<z_1$ inside the disc of convergence, the integral over $[\delta,z_1]$ being a constant: each monomial $t^{e-\beta}$ with $e\ne\beta-1$ yields a power of $z$, whereas the one with $e=\beta-1$ yields $\log(\delta/z)$; replacing $\delta$ by $z_1$ changes only the constant term. Since $\rho/z$ is a convergent Puiseux series with constant term $\Xi_0$, the change from $z$ to $\rho$ transforms powers into series of powers and changes $\log(1/z)$ into $\log(1/\rho)$ plus a constant and terms tending to zero.
\end{proof}

\section{Logarithmic resonances of rotational cusps}\label{sec:risonanze}

Consider the rotational cusps of \cite[Prop.~5.12]{V11},
\[
X:=\{(x,z)\in\R^d\times\R:\ |x|^2=z^{2\gamma},\ z\ge0\},\qquad d\ge3,\ \gamma>1\ \text{rational},
\]
which are the case $p=1$, $Y_1=S^{d-1}$, $\alpha_1=1$ of Section~\ref{sec:collasso}; here $\beta=\gamma(d-1)$ and $\Xi(z)=(1+\gamma^2z^{2\gamma-2})^{1/2}$.

\begin{theorem}\label{thm:risonanza}
The expansion of the profile $\varsigma_{\boldsymbol\gamma}$ contains a term $c_{\log}\log(1/\rho)$ with $c_{\log}\ne0$ if and only if there exists an integer $\mathsf n\ge1$ with
\[
\mathsf n(2\gamma-2)=\beta-1,\qquad\text{that is}\qquad \gamma=\frac{2\mathsf n-1}{2\mathsf n-d+1},\quad 2\mathsf n>d-1 .
\]
In this case $c_{\log}=\binom{1/2}{\mathsf n}\gamma^{2\mathsf n}$. For every $d\ge3$ there are infinitely many values of $\gamma>1$ for which this occurs; for comparison, when $d=2$ and $\gamma>1$, there is no such resonance.
\end{theorem}

\begin{proof}
We have $\Xi(t)=\sum_{\mathsf n\ge0}\binom{1/2}{\mathsf n}\gamma^{2\mathsf n}t^{\mathsf n(2\gamma-2)}$ for small $t$; the exponents are distinct and the coefficients nonzero. Proposition~\ref{prop:log-multipla} gives the first assertion. The equation $\mathsf n(2\gamma-2)=\gamma(d-1)-1$ is equivalent to $\gamma(2\mathsf n-d+1)=2\mathsf n-1$; for $2\mathsf n>d-1$ the solution is $\gamma=(2\mathsf n-1)/(2\mathsf n-d+1)$, which is greater than $1$ if and only if $d>2$.
\end{proof}

\begin{example}
For $d=3$ and $\mathsf n=2$ we have $\gamma=3/2$, that is, the cusp $|x|^2=z^3$ in $\R^4$, with $\beta=3$ and
\[
\varsigma_{\boldsymbol\gamma}(z)=\tfrac12z^{-2}+\tfrac98z^{-1}-\tfrac{81}{128}\log(1/z)+O(1).
\]
The tangent cone of this cusp is a half-line, and the germ does not satisfy the regular spherical blow-up condition (R) of \cite[Def.~5.7]{V11}. The possible presence of a logarithmic term, predicted in the expansion for germs with regular blow-up in \cite[Theorem~5.9]{V11}, does not characterize condition (R): it also occurs for the cusp just considered.
\end{example}

\section{A semialgebraic germ with connected link and failure of continuous extension}\label{sec:collo}

In this section we construct, for every integer $k\ge1$, a semialgebraic germ $X_k\subset\R^4$ of dimension $3$, with link diffeomorphic to $S^2$, satisfying the assumptions of \S\ref{sec:prelim} and admitting a bounded harmonic function with no limit at $0$. The link consists of two lobes of diameter comparable to $r$ joined by a neck of radius comparable to $r^{k+1}$. This gives a negative answer to the first part of question~(Q4): connectedness of the link does not suffice for continuous extension.

\subsection{The model}

Use the coordinates $(x_1,x_2,x_3,z)$ of $\R^4$. For $z>0$ and $t\in[-2,2]$ set
\[
\mathsf H(t,z):=\frac{t^2-1+\sqrt{(t^2-1)^2+4z^{2k}}}{2},\qquad
\mathsf B:=(4-t^2)\mathsf H,\qquad
\mathsf r:=z\sqrt{\mathsf B},
\]
\[
\Phi(z,t,\theta):=\bigl(zt,\ \mathsf r(z,t)\cos\theta,\ \mathsf r(z,t)\sin\theta,\ z\bigr),\qquad\theta\in\R/2\pi\mathbb Z,
\]
and $X_k:=\{0\}\cup\Phi\bigl((0,z_*)\times[-2,2]\times\R/2\pi\mathbb Z\bigr)$, with $z_*\in(0,1)$ fixed. The function $\mathsf H$ is the unique positive root of $\mathsf H^2-(t^2-1)\mathsf H-z^{2k}=0$, and
\[
X_k\setminus\{0\}=\bigl\{x_2^2+x_3^2=z^2\,\mathsf B(x_1/z,z),\ |x_1|\le2z,\ 0<z<z_*\bigr\}.
\]
Denote by $g_{zz},g_{zt},g_{tt}$ the coefficients of the block of the induced metric in the variables $(z,t)$, with $\mathsf D:=g_{zz}g_{tt}-g_{zt}^2$, and by $\Sigma_z:=X_k\cap\{x_4=z\}$ the fixed-height cross-sections.

\begin{lemma}\label{lemma:modello}
\begin{enumerate}[label=(\alph*)]
\item $X_k$ is semialgebraic, closed in a neighborhood of $0$, and $X_k\setminus\{0\}$ is an analytic submanifold of dimension $3$; the cross-sections $\Sigma_z$ and the links $L_r$, for small $r$, are diffeomorphic to $S^2$.
\item In the coordinates $(z,t,\theta)$, for $|t|<2$, we have $g_{zz}=1+t^2+\mathsf r_z^2$, $g_{zt}=zt+\mathsf r_z\mathsf r_t$, $g_{tt}=z^2+\mathsf r_t^2$, $g_{\theta\theta}=\mathsf r^2$, the other coefficients vanish, and
\[
\mathsf D=z^2+\mathsf r_t^2+(z\mathsf r_z-t\mathsf r_t)^2,\qquad d\tau=\mathsf r\sqrt{\mathsf D}\,dz\,dt\,d\theta .
\]
\item $0\le z\mathsf H_z\le 2k\min\{\mathsf H,z^k\}$ and $|\mathsf H_t|\le4$; hence $\mathsf B$, $\mathsf B_t$, $\mathsf r_z$ are bounded, $|z(\sqrt{\mathsf B})_z|\le Cz^{k/2}$ and $\mathsf r\sqrt{\mathsf D}\le Cz^2$.
\item $\rho=z(1+t^2+\mathsf B)^{1/2}$, $\partial_z\rho\ge1$, $z\le\rho\le C_\rho z$, and $\tau(\Omega_r)\le Cr^3$.
\item $\bigl|\nabla_X\rho-\Phi/|\Phi|\bigr|\le C\rho^{k/2}$; in particular $|\nabla_X\rho|\to1$ and $X_k$ satisfies the Whitney condition at $0$.
\end{enumerate}
\end{lemma}

\begin{proof}
(a) Semialgebraicity follows from the description by the equation for $\mathsf H$. For $|t|<2$ we have $\mathsf B>0$, and $\Phi$ is an immersion; at $t=\pm2$ the function $\mathsf B$ has a simple zero in $t$, and the implicit function theorem, applied to the equation $x_2^2+x_3^2=z^2\mathsf B(x_1/z,z)$, gives analytic charts in the variables $(z,x_2,x_3)$. Each cross-section $\Sigma_z$ is a surface of revolution whose squared radius has two simple zeros at the endpoints $t=\pm2$; in the regular charts described above it closes smoothly at the two poles and is diffeomorphic to $S^2$; by (d) and (e), $\rho$ has no critical points near $0$ and the links are spheres. The point $0$ is singular: the tangent cone at $0$ is contained in $\{x_4\ge0\}$ and contains vectors with $x_4>0$, so it is not a vector space.

(b) We have $\Phi_z=(t,\mathsf r_z\cos\theta,\mathsf r_z\sin\theta,1)$, $\Phi_t=(z,\mathsf r_t\cos\theta,\mathsf r_t\sin\theta,0)$, $\Phi_\theta=\mathsf r(0,-\sin\theta,\cos\theta,0)$; $\mathsf D$ is the squared norm of the cross product of $(t,\mathsf r_z,1)$ and $(z,\mathsf r_t,0)$.

(c) Differentiating the equation for $\mathsf H$ gives $z\mathsf H_z=2kz^{2k}/\sqrt{(t^2-1)^2+4z^{2k}}$ and $\mathsf H_t=2t\mathsf H/\sqrt{(t^2-1)^2+4z^{2k}}$. Since $\sqrt{(t^2-1)^2+4z^{2k}}=2\mathsf H-(t^2-1)\ge\max\{\mathsf H,2z^k\}$, the first line follows. Then $0\le z\mathsf B_z\le C\min\{\mathsf B,z^k\}$ and $|z(\sqrt{\mathsf B})_z|=z\mathsf B_z/(2\sqrt{\mathsf B})\le C\min\{\mathsf B,z^k\}/\sqrt{\mathsf B}\le Cz^{k/2}$. We have $\mathsf r_z=\sqrt{\mathsf B}+z(\sqrt{\mathsf B})_z$, $\mathsf r\mathsf r_t=z^2\mathsf B_t/2$, and $\mathsf r\sqrt{\mathsf D}\le C(z\mathsf r+\mathsf r|\mathsf r_t|)\le Cz^2$.

(d) $\rho^2=z^2(t^2+\mathsf B+1)$ and $\partial_z\rho=(1+t^2+\mathsf B)^{1/2}+z\mathsf B_z/(2(1+t^2+\mathsf B)^{1/2})\ge1$. The volume bound follows from (c) and $\Omega_r\subset\{z<r\}$.

(e) We have $\Phi_z=\Phi/z+\bigl(0,z(\sqrt{\mathsf B})_z\cos\theta,z(\sqrt{\mathsf B})_z\sin\theta,0\bigr)$, with $|\Phi/z|=\rho/z\ge1$. The vector $\mathbf V:=\Phi_z/\partial_z\rho$ is tangent to $X_k$ and $|\mathbf V-\Phi/|\Phi||\le Cz^{k/2}$; since $\nabla_X\rho$ is the orthogonal projection of $\widehat\Phi=\Phi/|\Phi|$ onto the tangent space, and $\mathbf V$ is tangent, we have
\[
|\mathbf V-\nabla_X\rho|\le|\mathbf V-\widehat\Phi|,
\qquad |\nabla_X\rho-\widehat\Phi|\le|\mathbf V-\widehat\Phi|.
\]
Both left-hand sides are thus $O(\rho^{k/2})$, uniformly; the estimates extend to the poles by continuity in regular charts. The second implies the Whitney condition at zero.
\end{proof}

\subsection{Neck conductance and spectral gap of the cross-sections}

These two results quantify the weakness of the neck in the cross-sections $\Sigma_z$; they are not used in the proof of Theorem~\ref{thm:discontinuita}.

\begin{proposition}\label{prop:conduttanza-collo}
Let $0<\tau_0<1$ and let $\mathsf K_{\tau_0}(z)$ be the infimum of the Dirichlet energy on $\Sigma_z\cap\{|t|<\tau_0\}$ of functions equal to $0$ on $\{t=-\tau_0\}$ and to $1$ on $\{t=\tau_0\}$. Then
\[
\mathsf K_{\tau_0}(z)=2\pi\Bigl(\int_{-\tau_0}^{\tau_0}\frac{(z^2+\mathsf r_t^2)^{1/2}}{\mathsf r}\,dt\Bigr)^{-1}\sim
2\pi\Bigl(\int_{-\tau_0}^{\tau_0}\Bigl(\frac{1-t^2}{4-t^2}\Bigr)^{1/2}dt\Bigr)^{-1}z^k .
\]
\end{proposition}

\begin{proof}
The metric of $\Sigma_z$ is $(z^2+\mathsf r_t^2)dt^2+\mathsf r^2d\theta^2$; averaging in $\theta$ does not increase the energy, and for functions of $t$ the Cauchy--Schwarz inequality gives the formula, with equality for the function whose derivative is proportional to $(z^2+\mathsf r_t^2)^{1/2}/\mathsf r$. For $|t|\le\tau_0$ we have $\mathsf H=z^{2k}(1-t^2)^{-1}(1+O(z^{2k}))$ with derivatives in $t$, hence $\mathsf r=z^{k+1}((4-t^2)/(1-t^2))^{1/2}(1+o(1))$ and $\mathsf r_t=O(z^{k+1})$.
\end{proof}

Set $\varpi:=z^k$, $\mathsf w_\varpi:=(\mathsf B+\mathsf B_t^2/4)^{1/2}$ and $\mathsf c_\varpi:=\mathsf B/\mathsf w_\varpi$, functions of $t$; moreover,
\[
I_0:=\int_{-1}^1\Bigl(\frac{1-t^2}{4-t^2}\Bigr)^{1/2}dt,\qquad
J_0:=\int_1^2\bigl((4-t^2)(t^2-1)+t^2(5-2t^2)^2\bigr)^{1/2}dt,\qquad V_0:=2\pi J_0 .
\]

\begin{theorem}\label{thm:gap}
Let $\lambda_1(\Sigma_z)$ be the first positive eigenvalue of $-\Delta_{\Sigma_z}$. Then
\[
\lim_{z\to0}z^{2-k}\lambda_1(\Sigma_z)=\frac{2}{I_0J_0}=0.911638\ldots
\]
For small $z$ the eigenspace is one-dimensional, generated by a rotationally invariant eigenfunction odd under $t\mapsto-t$; normalized by $z^{-2}\int_{\Sigma_z}f^2\,dS=1$, with sign fixed by $f(3/2)>0$, it converges locally uniformly on $\{1<|t|<2\}$ to the constants $\pm(2V_0)^{-1/2}$.
\end{theorem}

\begin{proof}
The metric of $\Sigma_z$ is $z^2\bigl[(1+\mathsf B_t^2/(4\mathsf B))dt^2+\mathsf B\,d\theta^2\bigr]$; for $f=f(t)$ the energy is $\mathcal E_z(f)=2\pi\int\mathsf c_\varpi|f'|^2dt$ and the mass is $z^2\cdot2\pi\int\mathsf w_\varpi f^2dt$. The formulas in Lemma~\ref{lemma:modello}(c) give $\mathsf B,|\mathsf B_t|,\mathsf w_\varpi\le C$ and $\mathsf c_\varpi^{-1}\le C(\varpi^{-1}+(2-|t|)^{-1})$. Indeed, on $|t|\le3/2$ rationalizing $\mathsf H$ gives $\mathsf B\ge c\varpi^2$ and $|\mathsf B_t|/\mathsf B\le C/\varpi$, while for $3/2<|t|<2$ we have $\mathsf B\asymp2-|t|$ and $|\mathsf B_t|\le C$. We use the identity $\mathsf c_\varpi^{-1}=\sqrt{\mathsf B^{-1}+\mathsf B_t^2/(4\mathsf B^2)}$. In particular, dominated convergence for the resistance applies on intervals separated from the poles. For $|t|<1$, $\varpi\mathsf c_\varpi^{-1}\to((1-t^2)/(4-t^2))^{1/2}$; for $1<|t|<2$, $\varpi\mathsf c_\varpi^{-1}\to0$ and $\mathsf w_\varpi\to\mathsf w_0:=((4-t^2)(t^2-1)+t^2(5-2t^2)^2)^{1/2}$; for $|t|<1$, $\mathsf w_\varpi\to0$. With $t_*:=3/2$ dominated convergence gives $\varpi\int_{-t_*}^{t_*}\mathsf c_\varpi^{-1}\to I_0$ and $\int\mathsf w_\varpi\to2J_0$.

\emph{Compactness.} Let $f_\varpi$ be functions of $t$ with $\mathcal E_z(f_\varpi)\le C\varpi$ and $\int\mathsf w_\varpi f_\varpi^2\le C$. On compact subsets of $\{1<|t|<2\}$ the weights are bounded and bounded away from zero, and $|f_\varpi(t)-f_\varpi(t')|\le C\varpi^{1/2}|t-t'|^{1/2}$; along subsequences $f_\varpi\to c_\pm$ locally uniformly on $\{\pm t>1\}$. Near the pole $t=2$, the inequality
\[
\int_{t_*}^2\mathsf w_\varpi|f-f(t_*)|^2dt\le\Bigl(\int_{t_*}^2\mathsf c_\varpi|f'|^2dt\Bigr)\int_{t_*}^2\mathsf c_\varpi^{-1}(x)\int_x^2\mathsf w_\varpi\,dt\,dx\le C\int_{t_*}^2\mathsf c_\varpi|f'|^2dt,
\]
follows from Cauchy--Schwarz and $\mathsf c_\varpi^{-1}(x)\le C/(2-x)$, $\int_x^2\mathsf w_\varpi\le C(2-x)$. The inequality excludes concentration of mass near the poles, and the same estimate holds at the opposite pole. On the interval $[-t_*,t_*]$ the bound $\int\mathsf c_\varpi^{-1}\le C/\varpi$, combined with $\mathcal E_z(f_\varpi)\le C\varpi$ and boundedness of $f_\varpi(t_*)$, gives a uniform bound for $f_\varpi$ by Cauchy--Schwarz. The mass of the neck thus tends to zero by dominated convergence. Setting $f_{\rm lim}$ equal to $c_\pm$ on the two lobes and zero in the neck, we obtain
\[
\int_{-2}^2\mathsf w_\varpi|f_\varpi-f_{\rm lim}|^2dt\to0 .
\]
In particular, $\int\mathsf w_\varpi f_\varpi^2\to J_0(c_-^2+c_+^2)$, $\int\mathsf w_\varpi f_\varpi\to J_0(c_-+c_+)$, and, by Cauchy--Schwarz, the inner products of two sequences with the same bounds also converge.

\emph{Lower bound and recovery.} By Cauchy--Schwarz on $[-t_*,t_*]$, $\mathcal E_z(f_\varpi)\ge2\pi|f_\varpi(t_*)-f_\varpi(-t_*)|^2/\int_{-t_*}^{t_*}\mathsf c_\varpi^{-1}$, hence $\liminf\varpi^{-1}\mathcal E_z(f_\varpi)\ge2\pi(c_+-c_-)^2/I_0$. The function equal to $c_-$ for $t\le-t_*$, to $c_+$ for $t\ge t_*$ and to $c_-+(c_+-c_-)\int_{-t_*}^t\mathsf c_\varpi^{-1}/\int_{-t_*}^{t_*}\mathsf c_\varpi^{-1}$ in between attains equality. The recovery function converges to the two constant values on compact subsets of the lobes, since the resistance on each interval separated from the neck and the poles is bounded, whereas the denominator is asymptotic to $I_0/\varpi$. To estimate the first eigenvalue choose $c_+=-c_-\ne0$: the recovery function is odd by evenness of $\mathsf c_\varpi$ and has exactly zero mean for every $z$. Its normalized mass tends to $V_0(c_-^2+c_+^2)$.

\emph{Eigenvalue.} The average in $\theta$ and the remainder are orthogonal for energy and mass, and the remainder has Rayleigh quotient at least $cz^{-2}$, since $\mathsf B\le C$; for small $z$ the first eigenspace is therefore rotationally invariant. Normalizing by $2\pi\int\mathsf w_\varpi f^2=1$ and zero mean, every limit satisfies $V_0(c_-^2+c_+^2)=1$ and $c_-+c_+=0$, hence $c_\pm=\pm(2V_0)^{-1/2}$, and the lower bound and recovery give $z^2\lambda_1/\varpi\to4\pi/(I_0V_0)=2/(I_0J_0)$. If for a sequence $z_n\to0$ there were two orthonormal eigenfunctions, their limits would be proportional to $(-1,1)$ and orthonormal in $\R^2$ with the same mass, which is impossible; hence the eigenspace is simple. The reflection $t\mapsto-t$ is an isometry, the eigenfunction is even or odd, and evenness would contradict $c_-=-c_+\ne0$. The numerical value of $2/(I_0J_0)$ is obtained by quadrature.
With the sign chosen at $t_*$, every subsequential limit has the positive value on the lobe $t>1$; the compactness established above therefore implies locally uniform convergence of the entire family on the two open lobes. The result concerns the spectrum of the fixed-height cross-sections, not that of the Euclidean links.
\end{proof}

\subsection{The discontinuity theorem}

Fix an odd $\mathsf f\in C^\infty([-2,2])$, with $0\le\mathsf f\le1$ on $[0,2]$ and $\mathsf f=1$ on $[1/2,2]$. For $0<Z<z_*$ set $\mathcal D_Z:=X_k\cap\{0<z<Z\}$.

\begin{theorem}\label{thm:discontinuita}
For every integer $k\ge1$ the germ $X_k$ satisfies the assumptions of \S\ref{sec:prelim}, has link diffeomorphic to $S^2$, and for sufficiently small $Z$ there exist a function $u_Z$ harmonic on $\mathcal D_Z$, continuous on $\mathcal D_Z\cup\Sigma_Z$ with $u_Z=\mathsf f(t)$ on $\Sigma_Z$, with $|u_Z|<1$ in $\mathcal D_Z$, odd under $t\mapsto-t$, rotationally invariant, and a constant $c_Z>0$ such that, for every $\tau_1\in(1,2)$,
\[
\lim_{\rho\to0}\ \sup_{\tau_1\le t\le2}|u_Z-c_Z|=0,\qquad\lim_{\rho\to0}\ \sup_{-2\le t\le-\tau_1}|u_Z+c_Z|=0 .
\]
In particular, $u_Z$ is bounded and has no limit at $0$; every value in $[-c_Z,c_Z]$ is a cluster value of $u_Z$ at $0$.
\end{theorem}

The proof occupies the remainder of the subsection. The constants denoted by $C$ and $c$ depend only on $k$, $\mathsf f$ and, where indicated, on fixed compact intervals; they are independent of $Z$ and of the heights or radii under consideration.

\emph{Step 1: construction and energy.} For $0<\varepsilon<Z$ let $u_{Z,\varepsilon}$ be the solution, given by~\ref{T:dirichlet}, of the Dirichlet problem on $\overline{\mathcal D_Z}\cap\{z\ge\varepsilon\}$, a compact manifold with smooth boundary, with data $\mathsf f(t)$ on both boundary components. It minimizes energy among functions in $H^1$ with these traces and is therefore unique; the isometries $\theta\mapsto\theta+\theta_0$ and $t\mapsto-t$ leave it invariant, up to a sign change in the second case; thus $u_{Z,\varepsilon}$ is rotationally invariant, odd, and vanishes on $\{t=0\}$. The maximum principle on the half-domain $\{t>0\}$ gives $0\le u_{Z,\varepsilon}\le1$ for $t\ge0$. On the support of $\mathsf f'$, contained in $\{|t|\le1/2\}$, we have $\mathsf r\asymp z^{k+1}$, $\mathsf r_t=O(z^{k+1})$, $\mathsf r_z=O(z^k)$, $g_{zz}\asymp1$, $\mathsf D\asymp z^2$; the energy of the function $\mathsf f(t)$ is
\[
2\pi\int_\varepsilon^Z\!\!\int_{-2}^2\frac{\mathsf r\,g_{zz}}{\sqrt{\mathsf D}}\,\mathsf f'(t)^2\,dt\,dz\le C\int_0^Zz^k\,dz\le CZ^{k+1},
\]
and the same upper bound holds for $u_{Z,\varepsilon}$. By interior and boundary estimates at $\{z=Z\}$, a sequence $\varepsilon_n\to0$ yields a limit $u_Z$, harmonic on $\mathcal D_Z$, with the same symmetries, $0\le u_Z\le1$ for $t\ge0$, data $\mathsf f$ on $\{z=Z\}$ and $\int_{\mathcal D_Z}|\nabla_Xu_Z|^2d\tau\le CZ^{k+1}$. By the strong maximum principle $|u_Z|<1$. The limit $u_Z$ is also an energy minimizer among finite-energy functions with trace $\mathsf f$ on $\Sigma_Z$. Indeed, a competitor can be truncated to $[-1,1]$ without increasing energy. On the band $\varepsilon<z<2\varepsilon$ it is patched to $\mathsf f$ by a cutoff with derivative $O(\varepsilon^{-1})$, so as to have trace $\mathsf f$ on the inner boundary. The energy of the patched function in the band is at most a constant times the sum of the energies of the competitor and $\mathsf f$ in the band, plus $C\varepsilon^{-2}\tau(\{z<2\varepsilon\})$. The first two tend to zero by integrability and the third is $O(\varepsilon)$, by Lemma~\ref{lemma:modello} and $|\nabla_Xz|\le1$. Apply minimality on the annuli and pass to the limit by lower semicontinuity. No trace is prescribed at the vertex.

\emph{Step 2: gradient estimate in the neck.} Let $0<s<Z/2$ and $z=s\hat z$. On the rectangle $\mathcal R:=\{3/4\le\hat z\le5/4,\ |t|\le1/2\}$ the metric block in the variables $(\hat z,t)$, divided by $s^2$, is
\[
\mathsf M_s:=\begin{pmatrix}g_{zz}&g_{zt}/s\\ g_{zt}/s&g_{tt}/s^2\end{pmatrix}\Big|_{z=s\hat z}\longrightarrow\begin{pmatrix}1+t^2&\hat zt\\ \hat zt&\hat z^2\end{pmatrix}
\]
in $C^m(\mathcal R)$ for every $m$, since $\mathsf r(s\hat z,t)/s^{k+1}\to\hat z^{k+1}((4-t^2)/(1-t^2))^{1/2}$ with all derivatives. For rotationally invariant functions the Dirichlet form is $2\pi s^{k+1}\int\langle\mathsf Q_s\nabla u,\nabla u\rangle\,d\hat z\,dt$, with $\mathsf Q_s:=(\mathsf r/s^{k+1})(\det\mathsf M_s)^{1/2}\mathsf M_s^{-1}$, a matrix converging in $C^m(\mathcal R)$ to a uniformly positive definite matrix. The equation $\mathrm{div}(\mathsf Q_s\nabla u_Z)=0$ is written in nondivergence form with coefficients bounded in $C^{0,1}(\mathcal R)$ uniformly in $s$, and the interior Schauder estimate \cite[Theorem~6.2]{GT} on two concentric rectangles gives $|\partial_tu_Z(s,0)|\le C\sup|u_Z|\le C$.

\emph{Step 3: fluxes.} Let $\mathfrak F^z(s)$ be the flux of $u_Z$ through $\Sigma_s\cap\{t>0\}$ in the direction of increasing $z$. On $\{t=0\}$ we have $\mathsf r_t=0$, $g_{zt}=0$ and $u_Z=0$; the area element of the surface is $\mathsf r(1+\mathsf r_z^2)^{1/2}dz\,d\theta$ and the outward unit normal of the half-domain $t>0$ is $-z^{-1}\partial_t$. Set, with the sign opposite to that of the outward flux through this surface,
\[
\ell(z):=2\pi\,\mathsf r(z,0)(1+\mathsf r_z(z,0)^2)^{1/2}z^{-1}\partial_t u_Z(z,0),
\qquad 0\le\ell(z)\le Cz^k .
\]
Nonnegativity follows from $u_Z\ge0$ for $t\ge0$. The divergence theorem gives $\mathfrak F^z(s_2)-\mathfrak F^z(s_1)=\int_{s_1}^{s_2}\ell$, so $\mathfrak F^z$ has a limit as $s\to0$. Setting $\mathcal A(s):=\int_{\Sigma_s\cap\{t>0\}}|\nabla_Xz|\,dS=2\pi\int_0^2\mathsf r\,g_{tt}\mathsf D^{-1/2}dt\le Cs^2$, the Cauchy--Schwarz inequality and the coarea formula give
\[
\int_0^{Z/2}\frac{\mathfrak F^z(s)^2}{\mathcal A(s)}\,ds\le\int_{\mathcal D_Z}|\nabla_Xu_Z|^2d\tau<\infty,
\]
and since $\int_0s^{-2}ds=\infty$ the limit is zero. Hence $0\le\mathfrak F^z(s)=\int_0^s\ell\le Cs^{k+1}$. Let $\mathfrak F^\rho(r)$ be the flux through $L_r\cap\{t>0\}$ in the direction of increasing $\rho$; the divergence theorem on $\{s<z,\ \rho<r,\ t>0\}$ and the limit $s\to0$ give $0\le\mathfrak F^\rho(r)\le Cr^{k+1}$ for $r<Z/2$.

\emph{Step 4: averages on the links.} In the variables $(r,t,\theta)$, with $z=z(r,t)$ the inverse of $\rho(z,t)$, write $\hat\Phi:=\Phi/|\Phi|$, $\mathbf V:=\Phi_z/\partial_z\rho$ and $d\tau=r^2\mathsf j\,dr\,dt\,d\theta$, with $\mathsf j:=\mathsf r\sqrt{\mathsf D}/(r^2\partial_z\rho)$. By Lemma~\ref{lemma:modello}, $0\le\mathsf j\le C$ and $|\mathbf V-\nabla_X\rho|\le Cr^{k/2}$. For fixed $1<t<2$, $\mathsf j(r,t)\to\mathsf j_0(t)>0$, and for $0<t<1$, $\mathsf j(r,t)\to0$. Set
\[
\mathfrak m(r):=\int_0^2\mathsf j(r,t)\,dt,\qquad \bar u(r):=\mathfrak m(r)^{-1}\int_0^2\mathsf j(r,t)\,u_Z(z(r,t),t)\,dt ;
\]
by dominated convergence $\mathfrak m(r)\to\int_1^2\mathsf j_0>0$, hence $\mathfrak m\ge c$ for small $r$. Since $\partial_r[u_Z(z(r,t),t)]=\langle\nabla_Xu_Z,\mathbf V\rangle$ and, by the coarea formula, $\mathfrak F^\rho(r)=2\pi r^2\int_0^2\mathsf j\langle\nabla_Xu_Z,\nabla_X\rho\rangle dt$, we have, for $r$ in a compact subset of $(0,Z/2)$,
\[
\mathfrak m\,\bar u'=\frac{\mathfrak F^\rho(r)}{2\pi r^2}+\int_0^2\mathsf j\langle\nabla_Xu_Z,\mathbf V-\nabla_X\rho\rangle\,dt+\int_0^2\partial_r\mathsf j\,(u_Z-\bar u)\,dt .
\]
Differentiation under the integral sign is justified: near the poles $\mathsf r\sqrt{\mathsf D}=(\mathsf r^2z^2+(\mathsf r\mathsf r_t)^2+(z\mathsf r\mathsf r_z-t\mathsf r\mathsf r_t)^2)^{1/2}$ is the square root of a positive smooth function, since $\mathsf r\mathsf r_t=z^2\mathsf B_t/2\ne0$ at $t=\pm2$. Moreover, in the charts at the poles the smooth function $u_Z$ is rotationally invariant, hence can be expressed as a smooth function of $z$ and $x_2^2+x_3^2$ up to the axis. After substituting $x_2^2+x_3^2=z^2\mathsf B$, the required derivatives are continuous up to the poles for $r$ in a compact subset of $(0,Z/2)$. This justifies differentiating the average. For $s<Z/8$, the Caccioppoli estimate on $\{s/2<\rho<4s\}$, with $|u_Z|\le1$ and $\tau(\Omega_{4s})\le Cs^3$, gives $\int_{s<\rho<2s}|\nabla_Xu_Z|^2d\tau\le Cs$; by Cauchy--Schwarz,
\[
\int_s^{2s}\!\!\int_0^2\mathsf j\,|\nabla_Xu_Z|\,r^{k/2}\,dt\,dr\le C\Bigl(\int_{s<\rho<2s}|\nabla_Xu_Z|^2d\tau\Bigr)^{1/2}\Bigl(\int_s^{2s}r^{k-2}dr\Bigr)^{1/2}\le Cs^{k/2}.
\]
The flux term contributes $\int_s^{2s}Cr^{k-1}dr\le Cs^k$.

\emph{Step 5: variation of the density.} The function $\mathsf j$ is semialgebraic on $(0,r_1)\times(0,2)$, since so are $\mathsf H$, the roots, derivatives of smooth semialgebraic functions, and the inverse $z(r,t)$ \cite[Corollary~2.9]{Coste}. The same is true of $\partial_r\mathsf j$: its graph is defined by a first-order formula in $\mathsf j$, namely the $\varepsilon$--$\delta$ definition of the derivative, and is therefore semialgebraic \cite[Theorem~2.6]{Coste} (see also \cite[Exercise~2.10]{Coste}). The sets $\{\partial_r\mathsf j>0\}$, $\{\partial_r\mathsf j=0\}$, $\{\partial_r\mathsf j<0\}$ are semialgebraic, hence \cite[Proposition~2.1]{Coste} finite unions of sets defined by the signs of finitely many polynomials $\mathsf p_1,\dots,\mathsf p_{N_0}$ in $(r,t)$; for fixed $t$, each $\mathsf p_i(\cdot,t)$ not identically zero has at most $\deg_r\mathsf p_i$ zeros. Thus $(0,r_1)$ is divided into at most $N_1:=1+\sum_i\deg_r\mathsf p_i$ intervals, whose number is independent of $t$, on each of which $\mathsf j(\cdot,t)$ is monotone; since $0\le\mathsf j\le C$, $\int_0^{r_1}|\partial_r\mathsf j(r,t)|dr\le CN_1$. By dominated convergence,
\[
\mathcal V(s):=\int_0^2\!\!\int_0^s|\partial_r\mathsf j(r,t)|\,dr\,dt\longrightarrow0\qquad(s\to0).
\]
Summing the estimates of Step~4 over dyadic annuli, with $\mathfrak m\ge c$ and $|u_Z-\bar u|\le2$, we obtain that $\bar u$ has a limit $\bar u(0)$ and
\begin{equation}\label{eq:coda}
|\bar u(s)-\bar u(0)|\le C\bigl(s^{k/2}+\mathcal V(s)\bigr)\qquad(s<Z/8),
\end{equation}
with $C$ and $\mathcal V$ independent of $Z$.

\emph{Step 6: the limit is nonzero.} Rescale by a factor $Z^{-1}$. For $I\Subset(1,2)$ and $0<\eta_1<1$, on the region $\{\eta_1\le z/Z\le1,\ t\in I\}$ the rescaled metric converges smoothly to a nondegenerate metric, and the data on $\{z=Z\}$ equal $1$. The rescaled energy of $u_Z$, in dimension $3$, is $Z^{-1}\int|\nabla_Xu_Z|^2d\tau\le CZ^k$. For $\mathsf v:=1-u_Z$ and $\hat z=z/Z$, $|\mathsf v(\hat z,t)|^2\le(1-\hat z)\int_{\hat z}^1|\partial_{\hat z}\mathsf v|^2$; integrating, $\|\mathsf v\|_{L^2}\le C_{I,\eta_1}Z^{k/2}$, and the interior $L^2$--$L^\infty$ estimate \cite[Theorem~8.17]{GT} gives $\sup|\mathsf v|\le CZ^{k/2}$ on interior compact subsets. On the link $r=Z/16$ we have $z/Z\in[1/(16C_\rho),1/16]$; hence $u_Z\to1$ for each fixed $t\in(1,2)$, whereas in the neck $\mathsf j\to0$. By dominated convergence $\bar u(Z/16)\to1$ as $Z\to0$, and by~\eqref{eq:coda} $\bar u(0)\to1$. Therefore there exists a positive threshold such that for every $Z$ below it we have $c_Z:=\bar u(0)>0$.

\emph{Step 7: limits on the lobes.} Let $s_n\to0$. On the rescaled manifolds $s_n^{-1}X_k$ consider the functions $x\mapsto u_Z(s_n x)$. We identify them locally, for $1<t<2$, by the parameterizations $(\widehat z,t,\theta)\mapsto s_n^{-1}\Phi(s_n\widehat z,t,\theta)$; near $t=2$ we use the corresponding rescaled regular charts $(\widehat z,x_2/s_n,x_3/s_n)$. On each compact set separated from the vertex and the half-line $t=1$ the metrics converge smoothly to the metric of the cone $\mathcal L_+$ over the surface $\{x_4=1,\ x_2^2+x_3^2=(4-t^2)(t^2-1),\ 1\le t\le2\}$. By interior estimates and a diagonal argument, along subsequences we obtain a bounded harmonic function $U_\infty$ on this cone with the half-line $t=1$ and the vertex removed. The surface is smooth and diffeomorphic to $S^2$ in the variables $(t,x_2,x_3)$, since the derivatives of $(4-t^2)(t^2-1)$ at $t=1$ and $t=2$ are $6$ and $-12$; hence $\mathcal L_+$ is a Riemannian cone $dr^2+r^2h_+$ over a smooth compact connected link. The half-line $\{t=1\}$ has zero local capacity in dimension $3$: in the tube of radius $\delta$ around a compact segment of it, the function $\mathsf t_\delta$ equal to $0$ for distance $\le\delta^2$, to $1$ for distance $\ge\delta$ and to $\log(\mathrm{dist}/\delta^2)/\log(1/\delta)$ in between has energy $O(1/|\log\delta|)$. For a compactly supported cutoff separated from the vertex, the Caccioppoli inequality with the additional factor $\mathsf t_\delta$ gives, by boundedness of $U_\infty$ and the energy estimate for $\mathsf t_\delta$, uniformly bounded energy on compact subsets as $\delta$ tends to zero. Thus $\nabla U_\infty$ is locally square integrable across the half-line. The almost-everywhere extension of $U_\infty$ has this weak gradient: in integration by parts with cutoff $\mathsf t_\delta$, the term containing $\nabla\mathsf t_\delta$ tends to zero by its energy and boundedness of the function. In the weak formulation of the equation, test with $\varphi\mathsf t_\delta$ for every smooth $\varphi$ compactly supported away from the vertex. The term with $\varphi\nabla\mathsf t_\delta$ tends to zero by Cauchy--Schwarz and the one with $\mathsf t_\delta\nabla\varphi$ converges by dominated convergence. We obtain $\int\langle\nabla U_\infty,\nabla\varphi\rangle\,d\tau=0$ across the half-line. By elliptic regularity $U_\infty$ extends harmonically on the punctured cone $\mathcal L_+\setminus\{0\}$. Expanding $U_\infty$ in eigenfunctions of the Laplacian of $h_+$, with eigenvalues $\lambda_i$, each coefficient $b_i$ satisfies $b_i''+2r^{-1}b_i'-\lambda_ir^{-2}b_i=0$; for $\lambda_i>0$ the solutions are linear combinations of $r^{a_+}$ and $r^{a_-}$, with $a_+>0>a_-$, and the nonzero ones are unbounded on $(0,\infty)$; for $\lambda_0=0$ only the constant remains, since $r^{-1}$ is unbounded. Hence $U_\infty$ is constant. To identify it, compute $\bar u(s_n)$: on compact subsets of $\{1<t\le2\}$ convergence is uniform, in the neck $\mathsf j\to0$ and the remaining bands have small measure, with $\mathsf j\le C$ and $|u_Z|\le1$; by dominated convergence $\bar u(s_n)\to U_\infty$, and by~\eqref{eq:coda} $U_\infty=c_Z$. Since every sequence has a subsequence with the same limit, convergence to $c_Z$ is uniform on $\{\tau_1\le t\le2\}$; by oddness $u_Z\to-c_Z$ on $\{-2\le t\le-\tau_1\}$. Finally, $u_Z=0$ on $\{t=0\}$, and by the intermediate value theorem on the connected links every value in $[-c_Z,c_Z]$ is a cluster value. \qed

\begin{corollary}\label{cor:harnack-collo}
For every $b\ge2$ there exists $c_b>0$ such that the optimal constant $C_{H,b}(r)$ of (H$_b$) for $X_k$ satisfies $C_{H,b}(r)\ge c_b\,r^{-k/2}$ for small $r$. In particular, $X_k$ satisfies neither (H$_b$) nor (HS) for any $b$, nor conditions (U) and (T) of Definitions~\ref{def:U} and~\ref{def:T}.
\end{corollary}

\begin{proof}
Choose $Z:=16br$ in Theorem~\ref{thm:discontinuita}; since $z\le\rho$, the annulus $A(r/b,br)$ is contained in $\mathcal D_Z$. Let $p_r\in L_r$ with $t=3/2$ and let $p_r^-$ be the symmetric point, with $t=-3/2$. Under rescaling by the factor $Z^{-1}$ the point $p_r$ remains in a fixed compact subset of the positive lobe, and Step~6 gives $0<1-u_Z(p_r)\le C_bZ^{k/2}$. The function $\mathsf v=1-u_Z$ is positive and harmonic on $A(r/b,br)$, with $\mathsf v(p_r^-)=1+u_Z(p_r)\ge1$. Thus $\sup_{L_r}\mathsf v/\inf_{L_r}\mathsf v\ge c_bZ^{-k/2}$. Conditions (T) and (U) imply (H) for connected links (Proposition~\ref{prop:T-U} and Lemma~\ref{lemma:U-H}).
\end{proof}

\begin{remark}[reduced model]
A purely illustrative scalar two-channel model is associated with the form
\[
 \frac12\int_0^1\bigl(r^2(\mathsf a'^2+\mathsf b'^2)
                    +r^k(\mathsf a-\mathsf b)^2\bigr)\dd r,
\]
in which $r^2$ represents the radial weight of the lobes and $r^k$ the coupling. The half-difference $\mathsf h=(\mathsf a-\mathsf b)/2$ satisfies $(r^2\mathsf h')'=2r^k\mathsf h$. For $k>0$, with fixed $\mathsf h(0)$, the equation
\[
 \mathsf h(r)=\mathsf h(0)+2\int_0^rs^{-2}\int_0^st^k\mathsf h(t)\dd t\dd s
\]
has a continuous solution on an initial interval, by contraction in the uniform norm: the norm of the integral operator on $[0,R_{\rm net}]$ is at most $2R_{\rm net}^k/(k(k+1))$. Taking $\mathsf h(0)\ne0$ and a constant sum gives distinct limits. For $k=0$ the Euler equation has solutions that are combinations of $r$ and $r^{-2}$; boundedness forces the half-difference to tend to zero.
\end{remark}

\section{Minimal harmonic functions over the vertex and the Martin boundary}\label{sec:martin}

In this section $X$ satisfies the assumptions of \S\ref{sec:prelim}, with connected link, $d\ge3$ and $\tau(\Omega_r)\le Cr^d$; fix $R\in(0,r_0)$. The case of interest is $X=X_k$.

\begin{lemma}\label{lemma:unicita-limitata}
Let $v$ be harmonic on $\Omega_R$, bounded near $0$, continuous on $\Omega_R\cup L_R$ with $v=0$ on $L_R$. Then $v\equiv0$.
\end{lemma}

\begin{proof}
The function $v$ is bounded on $\Omega_R$. For $0<\varepsilon<R/2$ the radial cutoff linear in $\rho$ between $\varepsilon$ and $2\varepsilon$ has energy at most $C\varepsilon^{-2}\tau(\Omega_{2\varepsilon})\le C\varepsilon^{d-2}$, hence $\Cap_X(\varepsilon;R)\to0$. On the fixed collar $A(r_1,R)$ the Poincaré inequality with zero trace on $L_R$ gives $\|e_{\varepsilon,R}\|_{L^2(A(r_1,R))}\to0$, and interior estimates give $e_{\varepsilon,R}\to0$ locally uniformly. By the maximum principle on $A(\varepsilon,R)$, $|v|\le\sup|v|\,e_{\varepsilon,R}$; as $\varepsilon\to0$ we obtain $v=0$.
\end{proof}

\begin{proposition}\label{prop:green}
There exists a Green function $\mathscr G$ for $-\Delta_X$ on $\Omega_R$, vanishing on $L_R$: $\mathscr G(x,x')$ is positive, symmetric, harmonic in $x$ away from $x'$, with $-\Delta_X\mathscr G(\cdot,x')=\delta_{x'}$, and is the increasing limit of the Green functions $\mathscr G_\varepsilon$ of the annuli $A(\varepsilon,R)$. Fix $r_c\in(0,R)$ and $o\in L_{r_c}$, and set $\mathsf p=e_{r_c,R}$. There exists $C>0$ such that for every $x'\in\Omega_{r_c/2}$ we have
\[
 \mathscr G(o,x')\le C,\qquad
 \mathscr G(x,x')\le C\mathscr G(o,x')\mathsf p(x)
 \quad(x\in A(r_c,R)).
\]
Moreover, $\mathscr G$ is continuous off the diagonal and for every compact set $K\Subset\Omega_R$ there exists $C_K>0$ such that
\[
0<\mathscr G(x,x')\le C_K d_X(x,x')^{2-d}+C_K\qquad (x,x'\in K,\ x\ne x').
\]
\end{proposition}

\begin{proof}
\emph{Green functions on annuli.} On the compact manifold with smooth boundary $\overline{A(\varepsilon,R)}$ there exists the Dirichlet Green function $\mathscr G_\varepsilon$ for $-\Delta_X$. To construct it, in a chart around the pole use the local Green function of the divergence-form operator with matrix $\sqrt{\det g}\,g^{-1}$ (cf. \cite[Theorem~1.1]{GW}). In this chart the equation with Dirac mass is precisely the weak equation associated with the Riemannian energy. Multiplied by a cutoff equal to one near the pole and zero outside the chart, the local Green function has, for $-\Delta_X$, source equal to the Dirac mass plus a smooth function supported away from the pole. This second source is canceled by solving, by~\ref{T:dirichlet}, a Dirichlet problem with zero trace on the annulus. Positivity and symmetry are those of the Dirichlet Green function (maximum principle and self-adjointness of the operator). If $\varepsilon'<\varepsilon<\rho(x')$, comparison on the common annulus, applied to the difference in which the Dirac masses cancel, gives $\mathscr G_{\varepsilon'}(\cdot,x')\ge\mathscr G_\varepsilon(\cdot,x')$.

\emph{Estimate on the collar.} Fix $r_c$ and $o\in L_{r_c}$, and suppose $\rho(x')<r_c/2$ and $\varepsilon<\rho(x')$. The sum of the outward fluxes on the two boundary components is $-1$: this follows from the equation with unit mass, or from Green's formula on an annulus with a small ball around the pole removed. The normal derivatives on both boundary components are nonpositive; in particular $\int_{L_R}|\partial_\nu\mathscr G_\varepsilon(\cdot,x')|dS\le1$. On a compact connected neighborhood of $L_{r_c}$ separated from the poles and the boundary components, Harnack's inequality gives, with constants independent of $\varepsilon$ and $x'$,
\[
c\mathscr G_\varepsilon(o,x')\le \mathscr G_\varepsilon(x,x')\le C\mathscr G_\varepsilon(o,x')\quad(x\in L_{r_c}).
\]
Comparison with $\mathsf p=e_{r_c,R}$ extends the two inequalities to $A(r_c,R)$ by multiplying the bounds by $\mathsf p(x)$. On the boundary $L_R$ the lower inequality implies
\[
-\partial_\nu\mathscr G_\varepsilon(\cdot,x')\ge c\mathscr G_\varepsilon(o,x')(-\partial_\nu\mathsf p).
\]
Integrating and using Hopf's lemma gives
\[
c\mathscr G_\varepsilon(o,x')\int_{L_R}|\partial_\nu\mathsf p|\,dS\le1,
\]
where the integral is strictly positive. This yields the two estimates of the statement for the Green functions on annuli. The same argument applies to any family of poles with $\rho(x')\le r_c-\delta$ for fixed $\delta>0$, by choosing the neighborhood of the link to have thickness less than $\delta$.

\emph{Uniform local estimates and limit.} Fix $K\Subset\Omega_R$ and choose $r_-,r_+$ with
\[
0<2r_-<\min_K\rho\le\max_K\rho<r_+<R.
\]
Let $A_K:=A(r_-,r_+)$ and let $\mathscr G_{A_K}$ be the Dirichlet Green function of $A_K$. For each pole $x'\in K$, the collar estimate just proved, with $r_c=r_+$ and separation $r_+-\max_K\rho>0$, gives $\mathscr G_\varepsilon(x,x')\le C_K$ on $L_{r_+}$, for $\varepsilon<r_-$. On $L_{r_-}$ use symmetry: for $x\in L_{r_-}$ consider the pole $x$, choose a fixed radius $r_c$ with $2r_-<r_c<\min_K\rho$ and evaluate the collar estimate at $x'\in K$. We again have $\mathscr G_\varepsilon(x,x')\le C_K$ uniformly.

For completeness, the Green function of the fixed domain $A_K$ satisfies
\[
0<\mathscr G_{A_K}(x,x')\le C_K d_X(x,x')^{2-d}+C_K\quad(x,x'\in K,\ x\ne x').
\]
Indeed, cover $K$ by finitely many regular charts relatively compact in $A_K$, and choose for each a larger chart, with a cutoff identically one on a neighborhood of the smaller chart. The coefficients are smooth and uniformly elliptic on the fixed charts. Denoting by $\xi$ the coordinate map of the chart, Theorem 1.1 of \cite{GW} gives the upper bound $C|\xi(x)-\xi(x')|^{2-d}$ for the local Green function with variable pole in the smaller chart. On the support of the derivatives of the cutoff, uniformly separated from the pole, interior estimates also uniformly control the required derivatives; the residual smooth source is therefore uniformly bounded. The Dirichlet correction on $A_K$ is uniformly bounded: compare its absolute value with a multiple of the nonnegative solution of $-\Delta_X w=1$ on $A_K$, $w=0$ on $\partial A_K$. Choose the larger charts with image a Euclidean ball and the cutoffs supported inside them. For $x'$ in the smaller chart and $x$ in the support of the cutoff, the length of the coordinate segment gives $d_X(x,x')\le C|\xi(x)-\xi(x')|$, where $\xi$ denotes the chart coordinates. It follows that $|\xi(x)-\xi(x')|^{2-d}\le C^{d-2}d_X(x,x')^{2-d}$. Outside the support of the cutoff only the uniformly bounded correction remains. Finiteness of the family of charts gives the upper bound with a common constant. No control of charts near the vertex is required.

The difference $\mathscr G_\varepsilon(\cdot,x')-\mathscr G_{A_K}(\cdot,x')$ is harmonic on $A_K$ also at $x'$, by distributional cancellation of the sources and elliptic regularity \cite[Theorem~8.3.1]{Hor}, applied in charts. It has nonnegative, uniformly bounded values on $\partial A_K$, hence lies between zero and $C_K$ in $A_K$. It follows, with constants independent of $\varepsilon<r_-$ and of the pole $x'\in K$,
\[
0<\mathscr G_\varepsilon(x,x')\le C_K d_X(x,x')^{2-d}+C_K\qquad (x,x'\in K,\ x\ne x').
\]
The increasing limit is therefore finite off the diagonal, since any pair of distinct points is contained in a compact subset of $\Omega_R$. It is positive and symmetric. Interior estimates on compact sets separated from the pole give harmonicity and locally uniform convergence; the same estimate, applied in both variables by symmetry, gives equicontinuity on compact sets separated from the diagonal and joint continuity of the limit. The singularity $d_X(x,x')^{2-d}$ is locally integrable with respect to $d\tau(x)$; domination on each compact set allows passage to the limit in the distributional equation for a fixed pole. We obtain $-\Delta_X\mathscr G(\cdot,x')=\delta_{x'}$. For each fixed pole, choosing an outer collar separated from the pole, comparison with the corresponding equilibrium potential gives continuous vanishing on $L_R$. The estimates for poles with $\rho(x')<r_c/2$ pass to the limit, completing the proof.
\end{proof}

\begin{theorem}\label{thm:trasferimento}
Suppose there exists a bounded harmonic function on $\Omega_R$ with two distinct cluster values at $0$. Then the cone $\mathcal K_R$ contains two nonproportional elements; the Picard principle, condition (AD) and the Bôcher property do not hold for $X$.
\end{theorem}

\begin{proof}
Let $u_\star$ be the harmonic function in the assumption, and let $\mathsf y_n^{\pm}\to0$ with $u_\star(\mathsf y_n^{\pm})\to\ell_\pm$, $\ell_+\ne\ell_-$. Let $\Upsilon$ be a smooth function of $\rho$, equal to $1$ near $0$ and to $0$ near $L_R$; the functions $\mathsf s_1:=-\Delta_X\Upsilon$ and $\mathsf s_2:=-\Delta_X(\Upsilon u_\star)$ have compact support in $\Omega_R\setminus\{0\}$. By Green's formula on $A(\varepsilon,R)$,
\[
\Upsilon u_\star(x')=\int\mathscr G_\varepsilon(x,x')\,\mathsf s_2(x)\,d\tau(x)+\mathsf b_\varepsilon(x'),\qquad|\mathsf b_\varepsilon|\le\sup|u_\star|\,e_{\varepsilon,R},
\]
The analogous formula holds for $\Upsilon$. Fix $x'$ sufficiently close to zero that $\Upsilon(x')=1$ and $x'$ is separated from the supports of the sources. The Green functions on the annuli converge on their support with an integrable upper bound, by Proposition~\ref{prop:green}. Lemma~\ref{lemma:unicita-limitata} and its proof give $e_{\varepsilon,R}(x')\to0$. Passing to the limit gives, for $x'$ near zero,
\[
\int\mathscr G(x,x')\,\mathsf s_1(x)\,d\tau(x)=1,\qquad\int\mathscr G(x,x')\,\mathsf s_2(x)\,d\tau(x)=u_\star(x').
\]
For $x'\ne o$ and $x\ne x'$ set $\mathscr K(x,x'):=\mathscr G(x,x')/\mathscr G(o,x')$. Discard any terms of the sequences for which $\mathsf y_n^\pm=o$. By Harnack's inequality and interior estimates, along subsequences $\mathscr K(\cdot,\mathsf y_n^{\pm})\to\mathscr K^\pm$ locally uniformly on $\Omega_R$; the $\mathscr K^\pm$ are positive harmonic functions, equal to $1$ at $o$, and by Proposition~\ref{prop:green} satisfy $\mathscr K^\pm\le C\mathsf p$ on $A(r_c,R)$, hence belong to $\mathcal K_R$. Since $\int\mathscr K(\cdot,\mathsf y_n^{\pm})\mathsf s_1=1/\mathscr G(o,\mathsf y_n^{\pm})\ge1/C$ and the supports of $\mathsf s_1,\mathsf s_2$ are compact,
\[
\int\mathscr K^\pm\mathsf s_1\ge1/C,\qquad\frac{\int\mathscr K^\pm\mathsf s_2}{\int\mathscr K^\pm\mathsf s_1}=\ell_\pm .
\]
Thus $\mathscr K^+\ne\mathscr K^-$, and since both equal $1$ at $o$, they are not proportional: the Picard principle of \cite[Prop.~5.16]{V11} fails. By Lemma~\ref{lemma:unicita-limitata} both are unbounded near $0$.

If (AD) held, we would have $\mathscr K^\pm=q_\pm G_X+\mathsf b_\pm$ near $0$ with $q_\pm>0$ and $\mathsf b_\pm$ bounded; the function $q_-\mathscr K^+-q_+\mathscr K^-$ would be bounded near $0$ and vanish on $L_R$, hence would vanish by Lemma~\ref{lemma:unicita-limitata}, and $\mathscr K^+$, $\mathscr K^-$ would be proportional. For the Bôcher property: if one of the two functions does not tend uniformly to $+\infty$, the first clause fails; otherwise the second clause would give a nontrivial combination continuous at $0$, hence bounded and vanishing on $L_R$, and again identically zero.
\end{proof}

\begin{theorem}\label{thm:minimali}
Under the assumptions of Theorem~\ref{thm:trasferimento}, the set $\mathcal B:=\{\omega\in\mathcal K_R:\ \omega(o)=1\}$ is convex and compact for locally uniform convergence, has at least two extreme points, and each extreme point $\omega$ is a minimal positive harmonic function: if $\upsilon$ is harmonic and $0\le\upsilon\le\omega$, then $\upsilon=\upsilon(o)\,\omega$. Each extreme point of $\mathcal B$ is a locally uniform limit of $\mathscr K(\cdot,x'_n)$ for a sequence $x'_n\to0$.
\end{theorem}

\begin{proof}
\emph{Compactness and extreme points.} By Harnack's inequality, interior estimates and a diagonal argument, every sequence in $\mathcal B$ has a locally uniformly convergent subsequence; on the collar $A(r_c,R)$ we have uniformly $\omega\le C\mathsf p$, by Harnack's inequality on $L_{r_c}$ and comparison, and the limit belongs to $\mathcal B$. The topology is metrizable. By Theorem~\ref{thm:trasferimento}, $\mathcal B$ contains two distinct elements, and by the Krein--Milman theorem \cite{KM}, in the form \cite[Theorem~V.7.4]{Conway}, it has at least two extreme points.

\emph{Minimality.} Let $\omega$ be extreme and let $0\le\upsilon\le\omega$ be harmonic. If $\upsilon(o)=0$ or $\upsilon(o)=1$, the strong maximum principle gives $\upsilon=0$ or $\upsilon=\omega$. Otherwise, with $a_\upsilon:=\upsilon(o)\in(0,1)$, the functions $\upsilon/a_\upsilon$ and $(\omega-\upsilon)/(1-a_\upsilon)$ belong to $\mathcal B$ since they are dominated by multiples of $\omega$, and $\omega=a_\upsilon(\upsilon/a_\upsilon)+(1-a_\upsilon)(\omega-\upsilon)/(1-a_\upsilon)$; extremality gives $\upsilon=a_\upsilon\omega$.

\emph{Representation.} Let $\omega\in\mathcal B$ and $0<\varepsilon<\rho(o)$. For $0<\delta<\varepsilon$ solve, by~\ref{T:dirichlet}, the Dirichlet problem on $A(\delta,\varepsilon)$ with data $\omega$ on $L_\varepsilon$ and zero on $L_\delta$. The solutions increase as $\delta\downarrow0$ and are dominated both by $\omega$ and by $\max_{L_\varepsilon}\omega$. The limit $\upsilon_\varepsilon$ is harmonic on $\Omega_\varepsilon$ and has trace $\omega$ on $L_\varepsilon$: interior and boundary estimates apply on a fixed collar of $L_\varepsilon$, separated from the moving boundary. In particular the limit is smooth up to $L_\varepsilon$. Set
\[
\widetilde\omega_\varepsilon(x)=
\begin{cases}
\upsilon_\varepsilon(x),&0<\rho(x)<\varepsilon,\\
\omega(x),&\varepsilon\le\rho(x)<R.
\end{cases}
\]
The function is continuous, bounded, and vanishes continuously on $L_R$. Integrating by parts on both sides of $L_\varepsilon$, with $\nu$ normal toward increasing $\rho$, we have
\[
-\Delta_X\widetilde\omega_\varepsilon=\mu_\varepsilon,
\qquad d\mu_\varepsilon=(\partial_\nu\upsilon_\varepsilon-\partial_\nu\omega)\,dS\big|_{L_\varepsilon}\ge0.
\]
Indeed, $\omega-\upsilon_\varepsilon\ge0$ in the interior, vanishes on $L_\varepsilon$ and has nonpositive outward normal derivative. The density is smooth on the fixed link.

For $0<\delta<\varepsilon$, let
\[
\mathscr P_{\delta,\varepsilon}(x)=\int_{L_\varepsilon}\mathscr G_\delta(x,x')\,d\mu_\varepsilon(x').
\]
The potential is continuous across $L_\varepsilon$. Indeed, in regular charts of a compact neighborhood of the link, the area of $L_\varepsilon$ contained in a ball of radius $s$ is at most $C_\varepsilon s^{d-1}$, uniformly in the center. The density of $\mu_\varepsilon$ is bounded. By a dyadic decomposition and the local estimate in Proposition~\ref{prop:green}, also valid for the Green functions on annuli as $\delta$ decreases,
\[
 \int_{L_\varepsilon\cap\{d_X(x,x')<s\}}
 \mathscr G_\delta(x,x')\,d\mu_\varepsilon(x')
 \le C_\varepsilon(s+s^{d-1}).
\]
For $\delta$ not sufficiently small, one may bound $\mathscr G_\delta$ above by the Green function of a larger annulus. By introducing a continuous cutoff zero for $d_X(x,x')\le s/2$ and equal to one for $d_X(x,x')\ge s$, we approximate the potential uniformly in a neighborhood of the link by integrals of continuous kernels. The same argument applies to the global Green function. The constants may depend on $\varepsilon$. For fixed $\delta$, the support of the measure is separated from the two boundary components and continuity of the Dirichlet Green function up to the boundary off the diagonal gives zero trace of the potential.

By Fubini's theorem and the Green function equation, $-\Delta_X\mathscr P_{\delta,\varepsilon}=\mu_\varepsilon$ on $A(\delta,R)$. Subtracting the potential from $\widetilde\omega_\varepsilon$ gives a distributionally harmonic function, hence smooth in the interior by \cite[Theorem~8.3.1]{Hor} in charts, with nonnegative traces. By the maximum principle,
\[
 0\le\mathscr P_{\delta,\varepsilon}\le\widetilde\omega_\varepsilon.
\]
As $\delta$ decreases the potentials increase and converge by monotone convergence to
\[
 \mathscr P_\varepsilon(x):=\int_{L_\varepsilon}
 \mathscr G(x,x')\,d\mu_\varepsilon(x')\le\widetilde\omega_\varepsilon(x).
\]
The estimate for the part near the diagonal just proved gives continuity also on the link. Integrating against $-\Delta_X\varphi$, with $\varphi\in C_c^\infty(\Omega_R)$, domination permits passage to the limit in the distributional identity and we obtain $-\Delta_X\mathscr P_\varepsilon=\mu_\varepsilon$. The difference $\widetilde\omega_\varepsilon-\mathscr P_\varepsilon$ is harmonic on $\Omega_R$, bounded, and vanishes continuously on $L_R$ by domination. Lemma~\ref{lemma:unicita-limitata} implies that it vanishes. Therefore
\[
\widetilde\omega_\varepsilon(x)=\int_{L_\varepsilon}\mathscr G(x,x')\,d\mu_\varepsilon(x').
\]
Evaluating at $o$, we see that $d\pi_\varepsilon(x')=\mathscr G(o,x')\,d\mu_\varepsilon(x')$ is a probability measure. In particular, for $\rho(x)>\varepsilon$,
\[
\omega(x)=\int_{L_\varepsilon}\mathscr K(x,x')\,d\pi_\varepsilon(x').
\]

\emph{Limits of kernels and extreme points.} By Harnack's inequality on connected compact subsets of $\Omega_R$ and interior estimates, every sequence of poles $x'_n\to0$ admits a subsequence whose normalized kernels converge uniformly on each compact set, eventually free of poles. The limits belong to $\mathcal B$ by the collar comparison of Proposition~\ref{prop:green}. Denote by $\mathcal T\subset\mathcal B$ the set of all such limits. It is nonempty and closed: if $\omega_n\in\mathcal T$ converges locally uniformly to $\omega\in\mathcal B$, fix an increasing exhaustion $K_n$ by compact subsets of $\Omega_R$, with every compact subset of $\Omega_R$ contained in some $K_n$; we can take $K_n=\{a_n\le\rho\le b_n\}$ with $a_n\downarrow0$ and $b_n\uparrow R$, $0<a_n<b_n<R$. By the definition of $\mathcal T$ choose for each $n$ a pole $x'_n$, with $\rho(x'_n)<\min(1/n,\min_{K_n}\rho/2,\rho(o)/2)$, such that
\[
\sup_{K_n}|\mathscr K(\cdot,x'_n)-\omega_n|<1/n.
\]
Then $x'_n\to0$ and the kernels converge locally uniformly to $\omega$, so $\omega\in\mathcal T$. By compactness of $\mathcal B$, also $\mathcal T$ is compact. Moreover, for every compact set $K\Subset\Omega_R$,
\[
\lim_{\epsilon\downarrow0}\ \sup_{0<\rho(x')<\epsilon}\ 
\inf_{v\in\mathcal T}\sup_{x\in K}|\mathscr K(x,x')-v(x)|=0.
\]
If this were not true we would obtain a sequence of poles tending to zero with uniform distance from $\mathcal T|_K$ bounded away from zero, contradicting the preceding subsequence argument.

For $\varepsilon$ sufficiently small relative to $K_n$, the representation of $\omega$ by $\pi_\varepsilon$ is approximated uniformly on $K_n$ by finite convex combinations of restrictions of kernels: the kernel is continuous on $K_n\times L_\varepsilon$. By the property just proved, choosing $\varepsilon$ still smaller, one may replace each kernel by an element of $\mathcal T$ with uniform error less than $1/n$ on $K_n$. This gives a sequence of convex combinations of elements of $\mathcal T$ converging locally uniformly to $\omega$. It follows that $\mathcal B=\overline{\operatorname{conv}}\mathcal T$, since $\mathcal T\subset\mathcal B$ and $\mathcal B$ is closed and convex.

The discrete measures associated with the convex combinations have a subsequence converging weakly to a probability measure $\pi$ on the compact metrizable set $\mathcal T$. For each $x\in\Omega_R$ evaluation at $x$ is continuous on $\mathcal T$, hence $\omega(x)=\int_{\mathcal T}v(x)\,d\pi(v)$. The barycenters of every probability measure on $\mathcal T$ belong to $\mathcal B$: approximate the integral on each compact set by convex sums, using compactness of $\mathcal T$ in the topology of locally uniform convergence. If $\omega$ is extreme and $\pi$ is not concentrated on $\omega$, there exists a point of a countable dense subset of $\Omega_R$ at which evaluation is not almost surely constant; otherwise, by continuity, almost all functions would coincide with $\omega$. Separating the values of this evaluation by a threshold assigning strictly positive mass to both sides expresses $\omega$ as a nontrivial convex combination of two distinct barycenters in $\mathcal B$, contradicting extremality. Therefore $\pi$ is concentrated on $\omega$, in particular $\omega\in\mathcal T$. This is the argument of Milman's converse \cite[Theorem~V.7.8]{Conway} in the present function space.
\end{proof}

\begin{corollary}\label{cor:martin-collo}
For the germ $X_k$ and every sufficiently small $R>0$, the cone $\mathcal K_R$ has at least two extreme rays, generated by minimal positive harmonic functions that are limits of normalized Green kernels with poles approaching $0$; the Picard principle, (AD), and the Bôcher property fail.
\end{corollary}

\begin{proof}
By Lemma~\ref{lemma:modello}, $X_k$ satisfies the assumptions of this section with $d=3$: connected link and $\tau(\Omega_r)\le Cr^3$. By Theorem~\ref{thm:discontinuita}, for $R<Z$ the restriction of $u_Z$ to $\Omega_R$ satisfies the assumption of Theorem~\ref{thm:trasferimento}, which gives the failure of the Picard principle, (AD), and the Bôcher property. Theorem~\ref{thm:minimali} gives at least two extreme points of the normalized base, hence at least two extreme rays of $\mathcal K_R\cup\{0\}$, whose generators are minimal and arise as locally uniform limits of normalized Green kernels with poles approaching the vertex.
\end{proof}

\begin{remark}\label{oss:letteratura-martin}
It has not been proved that there are exactly two extreme rays, or that the kernels $\mathscr K^\pm$ in the proof of Theorem~\ref{thm:trasferimento} are extreme. For Euclidean domains with Dirichlet boundary condition, Benedicks \cite[Theorem~2 and~4]{Ben} establishes a dichotomy between one and two Martin points at infinity for the complement of a proper closed set contained in a hyperplane, whose points are Dirichlet regular, with an integral criterion; Ancona \cite[Theorem~3.1]{AncD} characterizes double multiplicity by a series in the symmetric, uniformly regular Denjoy domains considered in that paper. Ioffe and Pinsky \cite[Theorem~1 and Remark~3]{IP} obtain, under their assumptions on horn-shaped domains, one point or a sphere of minimal points at infinity; Ancona \cite[Theorem~1.1--1.3]{AncC} distinguishes uniqueness of positive harmonic functions vanishing on the boundary of a cone from the possible multiplicity of Martin points at infinity in dimension at least four. The cited results do not apply directly to the induced metric on a germ without boundary and do not provide an example equivalent to Theorem 9.4.
\end{remark}

\section{Singularities of type \texorpdfstring{$A_{\mathsf k}$}{Ak}: partial results}\label{sec:Ak}
\subsection{Conductance of the balanced band of \texorpdfstring{$A_{\mathsf k}$}{Ak}}
Consider the complex surface of real dimension four,
\[
 A_{\mathsf k}=\{(\mathsf X,\mathsf Y,\mathsf Z)\in\C^3:\mathsf X^2+\mathsf Y^2+\mathsf Z^{\mathsf k+1}=0\},\qquad \mathsf k\in\Z,\ \mathsf k\ge2.
\]
The index $\mathsf k$ is distinct from the parameter $k$ of the model in Section~\ref{sec:collo}. With $\mathsf U=\mathsf X+i\mathsf Y$, $\mathsf V=\mathsf X-i\mathsf Y$ we have $\mathsf U\mathsf V=-\mathsf Z^{\mathsf k+1}$ and the ambient metric is
\[
 |d\mathsf X|^2+|d\mathsf Y|^2+|d\mathsf Z|^2
 =\tfrac12(|d\mathsf U|^2+|d\mathsf V|^2)+|d\mathsf Z|^2.
\]
Let $\mathsf m=(\mathsf k+1)/2>1$. On the locus $\mathsf Z\ne0$ introduce $\varrho=|\mathsf Z|>0$, $\psi\in\R$ and the angles $\phi_1,\phi_2\in\R/(2\pi\Z)$. The global parameterization is
\begin{align*}
 \mathsf Z&=\varrho e^{i\phi_2},\\
 \mathsf U&=\varrho^{\mathsf m}e^{\psi+i\phi_1},\\
 \mathsf V&=-\varrho^{\mathsf m}e^{-\psi+i((\mathsf k+1)\phi_2-\phi_1)}.
\end{align*}
Indeed, $\mathsf U$ does not vanish and uniquely determines $\psi=\log(|\mathsf U|/\varrho^{\mathsf m})$ and $\phi_1$ modulo $2\pi$; $\mathsf V$ is determined by the equation. The expression is well defined even when $\mathsf m$ is a half-integer, since $\mathsf k+1$ is an integer and the radial power is positive real.

\begin{proposition}[Conductance of the balanced band]
Let $0<\varrho<\varrho_0<1$, with $\varrho_0$ fixed, and let
\[
 \psi_*(\varrho)=(\mathsf m-1)\log(1/\varrho).
\]
In the cross-section $\mathsf S_{\varrho}=A_{\mathsf k}\cap\{|\mathsf Z|=\varrho\}$, the Dirichlet capacity $K_{A_{\mathsf k}}(\varrho)$ of the band $|\psi|<\psi_*$ with values $0$ and $1$ on the respective boundary tori is
\[
 K_{A_{\mathsf k}}(\varrho)=\frac{4\pi^2\varrho}{\displaystyle\int_{-\psi_*}^{\psi_*}
 \left(1+\frac{\mathsf m^2\varrho^{2\mathsf m-2}}{\cosh(2\psi)}\right)^{-1/2}\dd \psi}
 \sim\frac{4\pi^2}{\mathsf k-1}\frac{\varrho}{\log(1/\varrho)}.
\]
\end{proposition}
\begin{proof}
Set $\mathsf P=\varrho^{2\mathsf m}\cosh(2\psi)$. Before restriction to the cross-section, the radial and meridional part of the metric is
\[
\bigl(1+\mathsf m^2\varrho^{-2}\mathsf P\bigr)d\varrho^2
+2\mathsf m\varrho^{2\mathsf m-1}\sinh(2\psi)d\varrho\,d\psi
+\mathsf P\,d\psi^2 .
\]
The angular block is orthogonal to this part. At fixed height $\varrho$, the squared differentials of $\mathsf U,\mathsf V$ are respectively
\begin{align*}
 \varrho^{2\mathsf m} e^{2\psi}(d\psi^2+d\phi_1^2),\\
 \varrho^{2\mathsf m} e^{-2\psi}(d\psi^2+(2\mathsf md\phi_2-d\phi_1)^2).
\end{align*}
The metric of the cross-section is thus
\begin{align*}
 g_{\psi\psi}=g_{\phi_1\phi_1}&=\mathsf P,\\
 g_{\psi\phi_1}=g_{\psi\phi_2}&=0,\\
 g_{\phi_1\phi_2}&=-\mathsf m\varrho^{2\mathsf m}e^{-2\psi},\\
 g_{\phi_2\phi_2}&=\varrho^2+2\mathsf m^2\varrho^{2\mathsf m}e^{-2\psi}.
\end{align*}
The determinant $\mathsf J$ of the angular block is
\[
 \mathsf J=\varrho^{2\mathsf m+2}\cosh(2\psi)+\mathsf m^2\varrho^{4\mathsf m},
\]
since $2e^{-2\psi}\cosh(2\psi)-e^{-4\psi}=1$. The energy density for a function depending only on $\psi$ is therefore
\[
 \mathsf W=\sqrt{\mathsf J/\mathsf P}
 =\varrho\left(1+\frac{\mathsf m^2\varrho^{2\mathsf m-2}}{\cosh(2\psi)}\right)^{1/2}.
\]
The variational problem is on a compact manifold with two boundary tori, with constant traces. Angular translations are isometries and preserve these traces. Averaging with respect to Haar measure on the torus does not increase the energy, by convexity of the quadratic form. It thus suffices to minimize $4\pi^2\int \mathsf W|\partial_{\psi} f|^2\dd \psi$ among functions with values $0$ and $1$ at the endpoints. Cauchy--Schwarz implies
\[
 1\le\left(\int \mathsf W|\partial_{\psi} f|^2\dd \psi\right)
       \left(\int \mathsf W^{-1}\dd \psi\right),
\]
where both integrals have endpoints $-\psi_*,\psi_*$. Equality is attained by the function
\[
f(\psi)=\frac{\displaystyle\int_{-\psi_*}^{\psi}\mathsf W(a)^{-1}\dd a}
{\displaystyle\int_{-\psi_*}^{\psi_*}\mathsf W(a)^{-1}\dd a},
\]
with $\varrho$ fixed and $a$ a dummy variable in the logarithmic coordinate. This gives the exact formula. Since $\cosh(2\psi)\ge1$, the factor in the denominator integral is $1+O(\varrho^{2\mathsf m-2})$ uniformly over the entire band, even though its width diverges. The denominator is therefore $2\psi_*(1+O(\varrho^{2\mathsf m-2}))$. The relation $2(\mathsf m-1)=\mathsf k-1$ gives the conclusion.
\end{proof}

The Euclidean distance $\rho$ satisfies $\rho^2=\varrho^2+\varrho^{2\mathsf m}\cosh(2\psi)$. In the band $|\mathsf U|,|\mathsf V|\le \varrho$, hence $\varrho\le\rho\le\sqrt2\,\varrho$. On each boundary torus one of the two moduli equals $\varrho$ and the other $\varrho^{\mathsf k}$. The cross-section $\mathsf S_{\varrho}$ is not, however, the Euclidean link $L_r$: the capacity just computed is not directly a capacity between regions of $L_r$.

\subsection{The two-channel model}
\begin{remark}
The preceding calculation suggests, in the network model, a coupling of order $\mathfrak r/\log(1/\mathfrak r)$, with radial weight $\mathfrak r^3$ (real dimension four). The model admits no bounded solutions with distinct limits in the two channels.
\end{remark}
More precisely, let $0<\mathfrak r<\mathfrak r_0<1$ and $\mathfrak c>0$, and let $\mathfrak v_+,\mathfrak v_-$ be the two scalar functions of the model, solutions of
\[
 (\mathfrak r^3\mathfrak v_\pm')'=\frac{\mathfrak c\mathfrak r}{\log(1/\mathfrak r)}(\mathfrak v_\pm-\mathfrak v_\mp).
\]
Here the variable $\mathfrak r$ is the abstract radial scale of the model. Setting $\mathfrak h=(\mathfrak v_+-\mathfrak v_-)/2$, we have
\[
 (\mathfrak r^3\mathfrak h')'=\frac{2\mathfrak c\mathfrak r}{\log(1/\mathfrak r)}\mathfrak h.
\]
If $\mathfrak h$ is bounded, the right-hand side is integrable near zero, so $\mathfrak r^3\mathfrak h'$ has a finite limit. This limit must be zero: a nonzero limit would force a fixed sign and growth of $|\mathfrak h'|$ at least as fast as $\mathfrak r^{-3}$, incompatible with boundedness. If moreover $\mathfrak h(\mathfrak r)\to \mathfrak h_0\ne0$, then, replacing $\mathfrak h$ by its negative if necessary, $\mathfrak h(\mathfrak r)\ge \mathfrak h_0/2>0$ for small $\mathfrak r$, and
\begin{align*}
 \mathfrak h'(s)&=2\mathfrak c s^{-3}\int_0^s\frac{t \mathfrak h(t)}{\log(1/t)}\dd t\\
 &\ge \mathfrak c\mathfrak h_0s^{-3}\int_{s/2}^{s}\frac{t}{\log(2/s)}\dd t
 =\frac{3\mathfrak c\mathfrak h_0}{8s\log(2/s)}.
\end{align*}
Consequently the integral equation with datum $\mathfrak h_0$,
\[
 \mathfrak h(\mathfrak r)=\mathfrak h_0+2\mathfrak c\int_0^{\mathfrak r} s^{-3}\int_0^s
 \frac{t \mathfrak h(t)}{\log(1/t)}\dd t\dd s,
\]
cannot be satisfied by a bounded solution with that limit: the integral $\int_0^{\mathfrak r} ds/(s\log(2/s))$ diverges, as follows by the substitution $\log(2/s)$. If the two limits are assumed to exist, they must therefore coincide. The equation for the sum is
\[
 (\mathfrak r^3(\mathfrak v_++\mathfrak v_-)')'=0,
\]
so the sum of a bounded pair is constant.

A possible objective for the PDE on the germ $A_{\mathsf k}$ is an upper bound $C_H(r)\le C\log(1/r)$ on the Harnack constant on the Euclidean links for positive harmonic functions on $A(r/4,4r)$. If this upper bound were proved, the series of reciprocals along $r_j=r_1 8^{1-j}$ would diverge and \cite[Theorem~3.3]{V11} would give continuity of bounded harmonic functions. The conductance of the balanced band and the network equation prove neither this upper bound nor continuity for the PDE. The radial part of the operator, the exchanges between regions, and the relation with the Euclidean links remain to be controlled.

\appendix
\section{Known results and elementary consequences}\label{app:noti}
We collect here direct consequences of previous results and elementary calculations for specific models. The source attribution distinguishes, where necessary, the classical result used from its application to the model under consideration.

\subsection{Uniform charts and thickness condition}
The following two conditions reduce (H) to Moser's Harnack inequality for divergence-form operators with measurable coefficients, in the form of \cite[Theorem~8.20]{GT}: if $w\in W^{1,2}_{\rm loc}(B(0,5))$, $w\ge0$, is a weak solution of $\operatorname{div}(\mathsf a_{\rm M}\nabla w)=0$ in $B(0,5)\subset\R^d$, with $\mathsf a_{\rm M}$ measurable, symmetric and $\lambda_{\rm M}|\xi|^2\le\mathsf a_{\rm M}\xi\cdot\xi\le\Lambda_{\rm M}|\xi|^2$, then $\sup_{B(0,1)}w\le C\inf_{B(0,1)}w$ with $C$ depending only on $d$ and $\Lambda_{\rm M}/\lambda_{\rm M}$. The theorem is applied in $B(0,9/2)$ with radius one; there are no lower-order terms.

\begin{definition}\label{def:U}
$X$ satisfies the \emph{uniform chart condition}, abbreviated (U), if there exist $K_{\rm cart}\ge1$, $n_U\in\mathbb N$ and $r_U\in(0,r_0/4)$ such that for every $r\in(0,r_U)$ there are open sets $U^{\rm c}_1,\dots,U^{\rm c}_n\subset A(r/4,4r)$, with $n\le n_U$, numbers $s^{\rm c}_i>0$ and bi-Lipschitz homeomorphisms, with respect to the distance $d_X$ restricted to $U_i^{\rm c}$ and the Euclidean distance on $B(0,5)$, $\varphi^{\rm c}_i:U^{\rm c}_i\to B(0,5)\subset\R^d$ satisfying
\begin{equation}\label{eq:carte-uniformi}
K_{\rm cart}^{-2}|\xi|^2\le(s^{\rm c}_i)^{-2}\,g\bigl(d(\varphi^{\rm c}_i)^{-1}\xi,d(\varphi^{\rm c}_i)^{-1}\xi\bigr)\le K_{\rm cart}^2|\xi|^2\qquad\text{a.e. in }B(0,5),\ \xi\in\R^d,
\end{equation}
and $L_r\subset\bigcup_{i=1}^n(\varphi^{\rm c}_i)^{-1}(B(0,1))$.
\end{definition}

\begin{lemma}\label{lemma:U-H}
If $X$ has connected link and satisfies (U), then $C_H(r)\le C_M^{\,n_U}$ for $0<r<r_U$, where $C_M$ is the Moser constant corresponding to the ellipticity ratio $K_{\rm cart}^{2d+4}$; in particular $X$ satisfies (H).
\end{lemma}
\begin{proof}
Let $v>0$ be harmonic on $A(r/4,4r)$ and let $w_i=v\circ(\varphi^{\rm c}_i)^{-1}$. Since $\varphi^{\rm c}_i$ is bi-Lipschitz, $w_i\in W^{1,2}_{\rm loc}(B(0,5))$ and, by the area formula, $w_i$ is a weak solution of $\operatorname{div}(\mathsf a_{{\rm cart},i}\nabla w_i)=0$, with $\mathsf a_{{\rm cart},i}=\sqrt{\det g_{\rm cart}}\,g_{\rm cart}^{-1}$, where $g_{\rm cart}$ is the matrix of $((\varphi^{\rm c}_i)^{-1})^*g$. The matrix $\mathsf a_{{\rm cart},i}$ changes only by a constant factor if $g_{\rm cart}$ is replaced by $(s^{\rm c}_i)^{-2}g_{\rm cart}$, whose eigenvalues lie in $[K_{\rm cart}^{-2},K_{\rm cart}^2]$ by \eqref{eq:carte-uniformi}; the ellipticity ratio is therefore at most $K_{\rm cart}^{2d+4}$. Setting $\mathcal V^{\rm c}_i=(\varphi^{\rm c}_i)^{-1}(B(0,1))$, Moser's inequality gives $\sup_{\mathcal V^{\rm c}_i}v\le C_M\inf_{\mathcal V^{\rm c}_i}v$.

Let $x,y\in L_r$. The set of points of $L_r$ reachable from $x$ through a sequence $\mathcal V^{\rm c}_{i_1},\dots,\mathcal V^{\rm c}_{i_p}$ of distinct indices, with $x\in\mathcal V^{\rm c}_{i_1}$ and $\mathcal V^{\rm c}_{i_j}\cap\mathcal V^{\rm c}_{i_{j+1}}\cap L_r\ne\emptyset$, is open and closed in $L_r$; since $L_r$ is connected, it coincides with $L_r$, and we may assume $y\in\mathcal V^{\rm c}_{i_p}$ with $p\le n\le n_U$. Choosing points $x_j\in\mathcal V^{\rm c}_{i_j}\cap\mathcal V^{\rm c}_{i_{j+1}}\cap L_r$, we have $v(x)\le C_Mv(x_1)\le\dots\le C_M^{\,p}v(y)$.
\end{proof}

Consider finite families of analytic graphs $G_{\rm gr,\iota}\subset X\cap\overline B(0,r_0)$ of maps $\varphi_{\rm gr,\iota}:D_{\rm gr,\iota}\to U_{\rm gr,\iota}^{\perp}$, where $D_{\rm gr,\iota}$ is open in the coordinate subspace $U_{\rm gr,\iota}$, and suppose that the graphs cover $X\cap\overline B(0,r_0)$. The differentials of the graph parameterizations are required to be bounded. For decompositions of subanalytic sets into graphs see Stasica \cite{Stasica}. Denote by $\pi_{\rm gr,\iota}$ the orthogonal projection onto $U_{\rm gr,\iota}$ and by $K_S\ge1$ a common upper bound for the largest singular value of $d(u\mapsto u+\varphi_{\rm gr,\iota}(u))$; the smallest singular value is at least $1$. Set $V^*(r)=\sup_{r/2\le t\le2r}\mathcal H^{d-1}(L_t)$.

\begin{definition}\label{def:T}
$X$ satisfies the \emph{thickness condition}, abbreviated (T), if there exist a decomposition as above, $c_T\in(0,1]$, $r_T\in(0,r_0/4)$ and a function $s_T:(0,r_T)\to(0,\infty)$ with
\[
 s_T(r)\le\frac r{16K_S},\qquad V^*(r)\le c_T^{-1}s_T(r)^{d-1},
\]
such that for every $r\in(0,r_T)$ and every $x\in L_r$ there is a graph $G_{\rm gr,\iota}$ of dimension $d$ with $x\in G_{\rm gr,\iota}$ and $B_{U_{\rm gr,\iota}}(\pi_{\rm gr,\iota}(x),s_T(r))\subset D_{\rm gr,\iota}$.
\end{definition}
Condition (T) requires more than the mere existence of a finite family of graphs: every point of $L_r$ must be interior, at scale $s_T(r)$, to a graph of dimension $d$, and the same scale must control $V^*(r)$. For the model of Section~\ref{sec:collo}, Corollary~\ref{cor:harnack-collo}, together with Proposition~\ref{prop:T-U} and Lemma~\ref{lemma:U-H}, shows that (T) fails.

\begin{proposition}\label{prop:T-U}
If $X$ satisfies (T), then $X$ satisfies (U) with $K_{\rm cart}=K_S$ and $n_U\le4(10K_S)^d/(\omega_dc_T)$, where $\omega_d$ is the volume of the unit ball of $\R^d$.
\end{proposition}
\begin{proof}
Let $r\in(0,r_T)$, $s=s_T(r)$ and $x\in L_r$, with $G_{\rm gr,\iota}$ as in (T); set $u_x=\pi_{\rm gr,\iota}(x)$ and $G_x=\{u+\varphi_{\rm gr,\iota}(u):|u-u_x|<s\}$. Segments in the ball of the projected domain lift to curves of length at most $K_S$ times the segment length. In particular, the points of $G_x$ have distance from $x$ at most $K_Ss\le r/16$; thus $G_x\subset A(r/2,2r)\subset X^*$. Being a submanifold of dimension $d$ contained in $X^*$, $G_x$ is open in $X^*$. Orthogonal projection and the bound on lifted segments show that the chart $p\mapsto5(\pi_{\rm gr,\iota}(p)-u_x)/s$ is bi-Lipschitz for the distance $d_X$ restricted to $G_x$. The graph metric has eigenvalues in $[1,K_S^2]$, so the inequality of Definition~\ref{def:U} holds with scale $s/5$ and constant $K_S$.

The intrinsic ball $B_X(x,s/5)$ is contained in the preimage $\mathcal V^{\rm c}_x$ of $B(0,1)$. A limit point in $X$ of points of $G_x$ that does not belong to $G_x$ projects onto $\partial B_{U_{\rm gr,\iota}}(u_x,s)$, since $\varphi_{\rm gr,\iota}$ is continuous on $D_{\rm gr,\iota}$. A curve in $X$ starting at $x$, of length less than $s/5$, has projection displaced by less than $s/5$; it therefore does not reach such points, remains in $G_x$ and has projection in $B_{U_{\rm gr,\iota}}(u_x,s/5)$.

Let $\{x_i\}\subset L_r$ be a maximal family with $d_X(x_i,x_j)\ge s/5$ for $i\ne j$. Then $L_r\subset\bigcup_iB_X(x_i,s/5)\subset\bigcup_i\mathcal V^{\rm c}_{x_i}$. The balls $B_X(x_i,s/10)$ are disjoint; each contains the graph over the ball of radius $s/(10K_S)$ centered at $u_{x_i}$, of volume at least $\omega_d(s/(10K_S))^d$, and is contained in $\{r-s\le\rho\le r+s\}$. By the coarea formula and $|\nabla_X\rho|\ge\frac12$ \cite[Lemma~2.2]{V11}, the volume of the latter set is at most $4sV^*(r)$. The number of points is therefore at most
\[
\frac{4sV^*(r)}{\omega_d(s/(10K_S))^d}=\frac{4(10K_S)^d}{\omega_d}\,\frac{V^*(r)}{s^{d-1}}\le\frac{4(10K_S)^d}{\omega_dc_T}.
\]
The bound holds for every finite family of separated centers and ensures that successive selection of the centers terminates after a uniformly bounded number of steps. We can therefore take $r_U=r_T$ and $n_U=\lfloor4(10K_S)^d/(\omega_dc_T)\rfloor$.
\end{proof}

\begin{theorem}\label{th:thick}
If $X$ satisfies the assumptions of \S\ref{sec:prelim}, with connected link, and (T), then every bounded harmonic function near zero extends continuously. For sufficiently small $R$ the Picard principle holds; moreover, (AD) and the Bôcher property hold.
\end{theorem}
\begin{proof}
By Proposition~\ref{prop:T-U} and Lemma~\ref{lemma:U-H}, (H) holds, hence (HS). Theorem 5.17 of \cite{V11} gives continuity, (AD) and the Bôcher property; Proposition 5.16 of \cite{V11} gives the Picard principle. The growth order of $G_X$ is not deduced from (T) alone.
\end{proof}

\subsection{Variable number of charts}
\begin{corollary}
Suppose $X$ has connected link and that, for each small $r$, there exist $n_{\rm cart}(r)$ charts as in Definition~\ref{def:U}, with the same metric constant $K_{\rm cart}$, but without a uniform bound on the number of charts. Let $C_M=C_M(d,K_{\rm cart}^{2d+4})\ge1$ be Moser's Harnack constant. If there exists a sequence $r_j\downarrow0$ such that $r_1<r_0/4$, $r_{j+1}\le r_j/8$, all the charts exist at the radii $r_j$, and
\[
 \sum_{j\ge1} C_M^{-n_{\rm cart}(r_j)}=+\infty,
\]
then every bounded harmonic function near zero has a limit at zero.
\end{corollary}
\begin{proof}
In the bi-Lipschitz charts the divergence-form matrix is $\sqrt{\det g_{\rm cart}}g_{\rm cart}^{-1}$, where $g_{\rm cart}$ denotes the metric matrix in the chart, and the constant scale cancels out of the equation. The ratio of eigenvalues is at most $K_{\rm cart}^{2d+4}$, as in Lemma~\ref{lemma:U-H}. Every subchart of radius one thus has a Harnack inequality with constant $C_M$ for positive harmonic functions on the annulus $A(r/4,4r)$. The subcharts cover the connected link; the graph of their intersections on the link is connected (its disconnection would separate the link into two disjoint open sets). Between two points one may choose a chain of at most $n_{\rm cart}(r)$ charts without repetitions. Applying Harnack's inequality successively at points in the intersections gives $C_H(r)\le C_M^{n_{\rm cart}(r)}$. Thus $\sum_j C_H(r_j)^{-1}=\infty$. For a function defined only in a smaller neighborhood, discard finitely many terms; \cite[Theorem~3.3]{V11} gives the conclusion.
\end{proof}
The criterion ensures continuity; divergence of the series alone does not give a uniform Harnack constant along a subsequence. Conversely, convergence of the indicated series does not imply discontinuity.

\subsection{Pluriharmonic functions on normal complex singularities}
\begin{theorem}\label{th:pluri}
Let $(V_{\rm an},0)\subset(\C^{n_{\rm amb}},0)$ be a normal complex analytic germ of complex dimension $d_{\C}\ge2$ with an isolated singularity. Let $Y_{\rm reg}=(V_{\rm an}\cap B(0,r_{\rm an}))\setminus\{0\}$ be connected, with $r_{\rm an}$ sufficiently small, and suppose $H^1(Y_{\rm reg};\R)=0$. Every real-valued pluriharmonic function $u$ on $Y_{\rm reg}$, even without an a priori boundedness assumption, is, on a possibly smaller punctured neighborhood, the real part of a holomorphic function $f_{\rm hol}$ on the germ $V_{\rm an}$. In particular it extends and $|u(x)-u(0)|\le C\rho(x)$.
\end{theorem}
\begin{proof}
The form $\alpha_{\rm hol}=2\partial u$ is holomorphic and closed: $\bar\partial\partial u=0$ and $\partial^2u=0$. De Rham's theorem and the vanishing of $H^1(Y_{\rm reg};\R)$ imply the vanishing of $H^1_{\rm dR}(Y_{\rm reg};\C)$, so all periods of $\alpha_{\rm hol}$ vanish. Fixing $x_{\rm hol}\in Y_{\rm reg}$, the integral $f_{\rm hol}(x)=\int_{x_{\rm hol}}^x\alpha_{\rm hol}$ is path independent, has differential $2\partial u$ and is holomorphic. Since $u$ is real, $d\Rea f_{\rm hol}=du$; adding a real constant gives $u=\Rea f_{\rm hol}$. The subset $\{0\}$ has complex codimension at least two in $V_{\rm an}$. Since $V_{\rm an}$ is normal, it is locally irreducible, hence locally pure dimensional, and $\{0\}$ is a thin subset of order two. By the Riemann extension theorem on locally pure dimensional complex spaces \cite[7.4.1]{GR}, a holomorphic function on the complement of a thin subset of order two has, without any boundedness assumption, a unique extension belonging to the sheaf $\widehat{\mathcal O}$ of weakly holomorphic functions; on a normal space this sheaf coincides with the structure sheaf. Apply this result to $f_{\rm hol}$. The structure sheaf of the analytic subspace is the quotient of the ambient sheaf: the extended germ has an ambient holomorphic representative. Its real part is smooth and Lipschitz on a smaller ball, yielding the estimate.
\end{proof}
\begin{corollary}\label{cor:pluri-AE}
Theorem~\ref{th:pluri} holds for the surfaces $A_{\mathsf k}$ with equation $\mathsf X^2+\mathsf Y^2+\mathsf Z^{\mathsf k+1}=0$ for every integer $\mathsf k\ge1$.
\end{corollary}
\begin{proof}
They are hypersurfaces with an isolated singularity, hence normal (a hypersurface is Cohen--Macaulay and therefore satisfies $S_2$, and the isolated singularity in complex dimension two gives $R_1$; apply Serre's criterion \cite[Theorem~23.8]{Mats}). The punctured germ retracts onto the link by the product structure. The link of $A_{\mathsf k}$ is diffeomorphic to the lens space $L(\mathsf k+1,\mathsf k)$ (\cite[Prop.~4.1]{Abbrescia}), with finite first homology group. Hence $H^1(Y_{\rm reg};\R)=0$.
\end{proof}
On a complex surface with the induced Kähler metric, pluriharmonic functions are harmonic, but the converse does not hold in general: the trace of the complex tensor of second derivatives can vanish without the tensor vanishing. The corollary does not establish (H), (AD) or the Bôcher property for general harmonic functions.

\subsection{Limit values for germs with several branches}
\begin{theorem}\label{th:rami}
Let $\Omega^{(j)}$ ($1\le j\le n_{\rm br}$) be the connected components of the punctured germ of $X$, finite in number. Suppose that every bounded harmonic function has a finite limit on each branch, for example if each satisfies (H). Let $\mathcal H_{\rm bd}$ be the real space of bounded harmonic germs, $\mathcal H_{\rm zero}$ the subspace of those with all limits zero, and $\mathcal H_{\rm cont}$ that of germs with a common limit. Then
\[
 0\longrightarrow\mathcal H_{\rm zero}\longrightarrow\mathcal H_{\rm bd}
 \xrightarrow{\operatorname{ev}_{\rm br}}\R^{n_{\rm br}}\longrightarrow0,
 \qquad \operatorname{ev}_{\rm br}(u)=(\ell_{{\rm br},j}(u))_{j=1}^{n_{\rm br}},
\]
is exact, where $\ell_{{\rm br},j}$ is the limit on branch $j$. Moreover,
\[
 \mathcal H_{\rm bd}/\mathcal H_{\rm cont}\simeq
 \R^{n_{\rm br}}/\R(1,\ldots,1)\simeq\R^{n_{\rm br}-1}.
\]
\end{theorem}
\begin{proof}
The map is linear and well defined by the assumption on limits. Its kernel is exactly $\mathcal H_{\rm zero}$. Every real vector is realized by a function constant on each component, with the prescribed constant: such a function is harmonic and bounded. Finally, $\mathcal H_{\rm cont}$ is the inverse image of the diagonal, and the isomorphism theorem gives the second assertion.
\end{proof}
If, in addition, (HS) holds on each branch, by \cite[Prop.~5.16 and Theorem~5.17]{V11} the decomposition applies separately on each branch. For a sufficiently small common outer radius, the cone of nonnegative harmonic functions with zero trace on the outer links (including the zero function) has $n_{\rm br}$ extreme rays, represented by the branch profiles with zero outer trace, extended by zero on the other components. Indeed, the restriction of such a harmonic function to each branch is zero or positive by the strong maximum principle and, in the latter case, a multiple of the profile. Functions strictly positive on all branches instead have all coefficients strictly positive in the cone with zero trace.

As an example, the complex quadric $V_{\rm quad}=\{x_1^2+x_2^2-x_3^2=0\}\subset\C^3$ is normal, being a hypersurface of complex dimension two with an isolated singularity. The real locus has two sheets, determined by the sign of $x_3$. The function equal to $1$ on the positive sheet and to $0$ on the negative sheet is harmonic and bounded on the punctured real locus but has no common limit. Complex normality does not ensure harmonic removability for the real locus without further assumptions on the branches.

\subsection{Monomial curves}
\begin{theorem}
Let $2\le m_{\rm P}<n_{\rm P}$ be coprime integers and let $V_{\rm P}=\{y^{m_{\rm P}}=x^{n_{\rm P}}\}\subset\C^2$. The normalization is $\nu_{\rm P}(t_{\rm P})=(t_{\rm P}^{m_{\rm P}},t_{\rm P}^{n_{\rm P}})$; setting $s_{\rm P}=|t_{\rm P}|$, the pullback metric is
\[
 \nu_{\rm P}^*g=(m_{\rm P}^2s_{\rm P}^{2m_{\rm P}-2}+n_{\rm P}^2s_{\rm P}^{2n_{\rm P}-2})|dt_{\rm P}|^2.
\]
Every bounded real-valued harmonic function near the vertex extends continuously and satisfies $|u(x)-u(0)|=O(\rho(x)^{1/m_{\rm P}})$; the exponent $1/m_{\rm P}$ is optimal for this estimate at the vertex. The positive profile, chosen with additive constant fixed by
\[
 G_{\rm P}^{(1)}(\nu_{\rm P}(t_{\rm P}))=\frac1{2\pi}\log\frac1{s_{\rm P}},
\]
has unit flux and
\[
 G_{\rm P}^{(1)}=\frac1{2\pi m_{\rm P}}\log\frac1\rho+
 \frac1{4\pi m_{\rm P}}\rho^{2(n_{\rm P}-m_{\rm P})/m_{\rm P}}
 +o\bigl(\rho^{2(n_{\rm P}-m_{\rm P})/m_{\rm P}}\bigr).
\]
If $u\circ\nu_{\rm P}=\Rea f_{\rm P}$, with $f_{\rm P}$ holomorphic on the disk, then $u$ is the restriction of the real part of an ambient holomorphic function if and only if $f_{\rm P}$ belongs to $\C\{t_{\rm P}^{m_{\rm P}},t_{\rm P}^{n_{\rm P}}\}$ (up to an imaginary constant).
\end{theorem}
\begin{proof}
The parameterization is finite and biholomorphic away from zero: by Bézout there exist integers whose linear combination of $m_{\rm P},n_{\rm P}$ equals one, giving the inverse by a Laurent monomial in nonzero $x,y$. It is surjective by coprimality and the disk is normal; hence it is the normalization. The metric formula follows by differentiating the two components. In dimension two the Laplacian of a metric conformal to the Euclidean metric is a positive multiple of the Euclidean Laplacian. The bounded harmonic function therefore extends across the point of the disk, and the interior estimate gives $|u\circ\nu_{\rm P}-u(0)|\le Cs_{\rm P}$. Since
\[
 \rho=s_{\rm P}^{m_{\rm P}}(1+s_{\rm P}^{2(n_{\rm P}-m_{\rm P})})^{1/2},
\]
the stated estimate follows. The function $\Rea t_{\rm P}$ attains the optimal exponent along $t_{\rm P}>0$ and cannot be the restriction of an ambient $C^1$ function.

For the profile, conformality implies that the product of the inward normal derivative and the line measure on $s_{\rm P}$ constant is $-(\partial_{s_{\rm P}}G_{\rm P}^{(1)})s_{\rm P}\dd\theta_{\rm P}$, with $\theta_{\rm P}=\arg t_{\rm P}$. Integrating gives one. The classical Bôcher theorem on the disk identifies the singular part of positive harmonic functions with a nonnegative multiple of the logarithm. The coefficient and the correction follow from the exact relation
\[
 \log(1/s_{\rm P})=m_{\rm P}^{-1}\log(1/\rho)+
 (2m_{\rm P})^{-1}\log(1+s_{\rm P}^{2(n_{\rm P}-m_{\rm P})}).
\]
Use $\log(1+v)=v+o(v)$ for $v\to0$ and $s_{\rm P}\sim\rho^{1/m_{\rm P}}$. If zero trace on an outer link is prescribed, the corresponding constant must be added to the profile; the same constant term is then not assigned to its expansion. The stated correction concerns the explicitly chosen profile, not an arbitrary profile modulo bounded harmonic functions.

On the full disk the extended harmonic function has a global harmonic conjugate. The ring of ambient holomorphic germs restricted and then pulled back is precisely $\C\{t_{\rm P}^{m_{\rm P}},t_{\rm P}^{n_{\rm P}}\}$. Two holomorphic functions on the disk with the same real part differ by an imaginary constant, which proves the equivalence. In particular, setting $\mathcal S_{\rm P}=\{i m_{\rm P}+j n_{\rm P}:i,j\in\Z_{\ge0}\}$, the monomial $\Rea t_{\rm P}^{j_{\rm P}}$ ($j_{\rm P}\ge1$ integer) does not extend as the real part of an ambient holomorphic function when $j_{\rm P}\notin\mathcal S_{\rm P}$: the series of an element of the ring has only exponents in $\mathcal S_{\rm P}$.
\end{proof}
The logarithmic coefficient identifies $m_{\rm P}$ and, in the monomial model with unit coefficients and the indicated additive normalization, the exponent of the first correction determines $n_{\rm P}$. It is not asserted that this correction is a bi-Lipschitz invariant.

\subsection{A real cone over a sextic curve}
Let $\Gamma_{\rm L}=\{(x,y)\in\R^2:y^2=(1-x^2)(1+x^2)^2\}$. Set
\[
 q_{\rm L}(\theta_{\rm L})=(\cos\theta_{\rm L},(1+\cos^2\theta_{\rm L})\sin\theta_{\rm L},1),
 \qquad \theta_{\rm L}\in\R/(2\pi\Z).
\]
The geometric cone $K_{\rm L}$ is the set $\{s_{\rm L}q_{\rm L}(\theta_{\rm L}):s_{\rm L}\in\R\}$; the half-cones $K_{\rm L,+}$ and $K_{\rm L,-}$ correspond to $s_{\rm L}\ge0$ and $s_{\rm L}\le0$. The curve $q_{\rm L}$ describes exactly $\Gamma_{\rm L}\times\{1\}$: for $x=\cos\theta_{\rm L}$ we have $y^2=(1+x^2)^2(1-x^2)$, and both signs of $y$ are obtained as $\sin\theta_{\rm L}$ varies. The cone does not coincide with the entire real locus of the homogenized polynomial $Q_{\rm L}(x,y,t)=y^2t^4+(x^2-t^2)(x^2+t^2)^2$. Indeed, for $t\ne0$, homogeneity of degree six gives $Q_{\rm L}(x,y,t)=t^6Q_{\rm L}(x/t,y/t,1)$, and $Q_{\rm L}(\cdot,\cdot,1)=0$ is the equation of $\Gamma_{\rm L}$; for $t=0$ we have $Q_{\rm L}=x^6$. The real locus is thus $K_{\rm L}\cup\{x=t=0\}$, with an additional line.
\begin{theorem}
Let $\ell_{\rm L}$ be the length of the unit link of $K_{\rm L,+}$, and $\sigma_{\rm L}\in\R/(\ell_{\rm L}\Z)$ the arc length parameter. The metric on the punctured half-cone is $d\rho^2+\rho^2d\sigma_{\rm L}^2$. Setting $\alpha_{\rm L}=2\pi/\ell_{\rm L}$, every bounded harmonic function on $0<\rho<R$ has the uniform estimate
\[
 u(\rho,\sigma_{\rm L})=u(0)+O((\rho/R)^{\alpha_{\rm L}}),
\]
with optimal exponent. Every positive harmonic function has, near the vertex, decomposition $u=qG_{\rm L,+}^{(1)}+u_b$, where $q\ge0$, $G_{\rm L,+}^{(1)}=\ell_{\rm L}^{-1}\log(R/\rho)$ has unit flux and $u_b$ is harmonic and bounded. Moreover,
\[
 \ell_{\rm L}\le\frac{5\pi}{2\sqrt2}<2\pi,
 \qquad \alpha_{\rm L}\ge\frac{4\sqrt2}{5}>1.
\]
In particular $u-u(0)=o(\rho)$ for every bounded harmonic function, that is, the derivative relative to the vertex vanishes in the sense of this estimate.
\end{theorem}
\begin{proof}
The parameterization $q_{\rm L}$ is regular and injective modulo its period: the first coordinate determines the cosine, the second the sine; when the derivative of the first vanishes, the derivative of the second does not. The normalized curve $\omega_{\rm L}=q_{\rm L}/|q_{\rm L}|$ is still regular and injective, since the third coordinate of $q_{\rm L}$ is one. In the arc length parameter $\omega_{\rm L}\cdot\omega_{\rm L}'=0$ and $|\omega_{\rm L}'|=1$, so that the pullback of $g$ by $\rho\omega_{\rm L}(\sigma_{\rm L})$ is the indicated metric. The coordinate
\[
 \zeta_{\rm L}=(\rho/R)^{\alpha_{\rm L}}e^{i\alpha_{\rm L}\sigma_{\rm L}}
\]
identifies the punctured half-cone conformally with the punctured disk: in coordinates $\log(\rho/R),\sigma_{\rm L}$ the metric is conformal to the flat metric and the exponential has exact period $\ell_{\rm L}$. Removability of the point in the disk and the interior gradient estimate give the remainder $O(|\zeta_{\rm L}|)$; $\Rea\zeta_{\rm L}$ proves optimality. Bôcher's theorem on the disk gives a logarithmic singular part and a harmonic remainder on the full disk; the latter has the same estimate. The inward flux of $\ell_{\rm L}^{-1}\log(R/\rho)$ is $-(\ell_{\rm L}\rho)\partial_\rho G_{\rm L,+}^{(1)}=1$.

For the length, with $c_{\rm L}=\cos\theta_{\rm L}$, we compute
\[
 |q_{\rm L}|^2=2+2c_{\rm L}^2-c_{\rm L}^4-c_{\rm L}^6\ge2,
 \qquad |q_{\rm L}'|^2=1-6c_{\rm L}^4+9c_{\rm L}^6.
\]
The first inequality follows from $2c_{\rm L}^2-c_{\rm L}^4-c_{\rm L}^6=c_{\rm L}^2(1-c_{\rm L}^2)(2+c_{\rm L}^2)$. The derivative of normalization is the projection of $q_{\rm L}'$ onto the plane perpendicular to $q_{\rm L}$, divided by $|q_{\rm L}|$, hence
\begin{align*}
 \ell_{\rm L}&\le\int_0^{2\pi}\frac{|q_{\rm L}'|}{|q_{\rm L}|}\dd\theta_{\rm L}\\
 &\le\left(2\pi-6\frac{3\pi}4+9\frac{5\pi}8\right)^{1/2}\pi^{1/2}
 =\frac{5\pi}{2\sqrt2}.
\end{align*}
We used Cauchy--Schwarz, $\int_0^{2\pi}|q_{\rm L}|^{-2}\dd\theta_{\rm L}\le\pi$ and the usual integrals of $\cos^4,\cos^6$ over the period, respectively $3\pi/4$ and $5\pi/8$ (obtained by expanding the powers into trigonometric harmonics). This gives the bound on $\alpha_{\rm L}$ and hence $u-u(0)=o(\rho)$.
\end{proof}
\begin{corollary}
On the full cone $K_{\rm L}$ every bounded harmonic function has two branchwise limits; it extends continuously if and only if they coincide. The cone of nonnegative harmonic functions with zero trace on the two outer links has two extreme rays.
\end{corollary}
\begin{proof}
The antipodal map identifies the two half-cones isometrically; each has the continuity property and the Picard principle, the latter by the decomposition just proved and the maximum principle for harmonic functions on the full disk. Apply Theorem~\ref{th:rami} and its consequence for the branch profiles. The function constant one on one half-cone and zero on the other shows failure of global continuity. This is not an example with connected link.
\end{proof}

\subsection{Metric induced by a homogeneous map}
\begin{lemma}
Let $Y_{\rm hom}$ be a Riemannian cone smooth away from the vertex, with metric $g_{Y_{\rm hom}}$ and compact unit cross-section, and let $\nu_{\rm hom}:Y_{\rm hom}\to X$ be a map of class $C^1$ away from the vertex, homogeneous of degree $a_{\rm hom}\ge1$ with respect to the cone dilations. Suppose it identifies the regular parts and that, at every point of the unit cross-section, its differential is injective on the entire tangent space of the source cone. If $r_{\rm hom}$ is the radial coordinate of the source cone, there exist $c,C>0$ such that
\[
 c^2r_{\rm hom}^{2(a_{\rm hom}-1)}g_{Y_{\rm hom}}
 \le\nu_{\rm hom}^*g\le
 C^2r_{\rm hom}^{2(a_{\rm hom}-1)}g_{Y_{\rm hom}}.
\]
Here the target metric is that induced by the Euclidean ambient space in which the dilations are taken.
\end{lemma}
\begin{proof}
For $y_{\rm hom}$ on the unit cross-section and for a tangent vector $v_{\rm hom}$ at that point, differentiating the relation $\nu_{\rm hom}\circ\operatorname{dil}_{r_{\rm hom}}=\operatorname{dil}_{r_{\rm hom}^{a_{\rm hom}}}\circ\nu_{\rm hom}$ gives
\[
 d\nu_{\rm hom}|_{r_{\rm hom}y_{\rm hom}}\bigl(d\operatorname{dil}_{r_{\rm hom}}v_{\rm hom}\bigr)
 =r_{\rm hom}^{a_{\rm hom}}d\nu_{\rm hom}|_{y_{\rm hom}}v_{\rm hom}.
\]
The source dilation multiplies norms by $r_{\rm hom}$. By compactness and injectivity, the singular values of $d\nu_{\rm hom}$ along the cross-section lie between two strictly positive constants. The ratio of norms at the dilated point therefore lies between $cr_{\rm hom}^{a_{\rm hom}-1}$ and $Cr_{\rm hom}^{a_{\rm hom}-1}$. Squaring gives the metric comparison.
\end{proof}
An example is given by the quadric cone $\widehat E_{\R}=\{z_1^2+z_2^2=z_3^2\}\subset\R^3$ and the cubic map
\[
 \nu_{\rm hom}(z_1,z_2,z_3)=(z_1z_3^2,(z_3^2+z_1^2)z_2,z_3^3).
\]
This map sends $\widehat E_{\R}$ to $K_{\rm L}$: with $c=z_1/z_3$ we have $x/t=c$ and $y/t=(1+c^2)z_2/z_3$, where $(z_2/z_3)^2=1-c^2$. Away from the vertex $z_3\ne0$ and the real inverse on the image cone is smooth: given $(x,y,t)$, we recover $z_3=\sqrt[3]{t}$, $z_1=x/z_3^2$ and $z_2=y/(z_3^2+z_1^2)$. The assumptions of the lemma hold and $\nu_{\rm hom}^*g_{K_{\rm L}}\asymp r_{\rm hom}^4g_{\widehat E_{\R}}$. The comparison is not an estimate for a Jacobian ideal of the complex normalization; it concerns the differential of the fixed real map. In real dimension two the conformal factor $r_{\rm hom}^4$ cancels from the energy density, leaving coefficients uniformly comparable to those of the conical metric.

\subsection{Counterexample to the angular estimate for the Gauss map}
\begin{proposition}
In general there is no constant $C$ such that, on a subanalytic germ smooth away from the vertex, the estimate
\[
 \delta_{\rm Gr}(T_xX,T_{x'}X)\le C\left|\frac{x}{|x|}-\frac{x'}{|x'|}\right|
\]
holds uniformly for $x,x'$ near zero. Here $\delta_{\rm Gr}$ is a distance inducing the topology of the Grassmannian, for instance the norm of the difference of the orthogonal projections.
\end{proposition}
\begin{proof}
Choose the rotational horn $X_{\gamma,S^{d-1}}$, with $d\ge2$ and $\gamma>1$ rational. At $\Pi(z,\omega_{\rm G})=(z^\gamma\omega_{\rm G},z)$, with $\omega_{\rm G}\in S^{d-1}$, the normalized direction and a unit normal are respectively
\[
 \frac{(z^{\gamma-1}\omega_{\rm G},1)}{\sqrt{1+z^{2\gamma-2}}},\qquad
 n_{\rm G}(z,\omega_{\rm G})=\frac{(\omega_{\rm G},-\gamma z^{\gamma-1})}{\sqrt{1+\gamma^2z^{2\gamma-2}}}.
\]
The second formula follows from the gradient of the defining function $|x|^2-z^{2\gamma}$. Choosing two orthogonal unit vectors $\omega_{{\rm G},1}$ and $\omega_{{\rm G},2}$, the normalized directions both converge to the vertical unit vector, whereas the normals converge to two orthogonal horizontal unit vectors. The limiting tangent planes are thus distinct; by continuity of the Grassmannian distance the left-hand side does not tend to zero whereas the right-hand side does. This is a contradiction.
\end{proof}
The disproved estimate concerns the actual directions $x/|x|$, which do not necessarily belong to the link of the tangent cone. The example rules out a uniform Lipschitz estimate for the Gauss map in the angular variable alone.

\section*{Acknowledgements}
\addcontentsline{toc}{section}{Acknowledgements}
The author thanks Prof.~Fulvio Lazzeri for the interesting discussion on the comparison between real and complex singularities, from which the examples of \S\S A.5--A.7 originate, together with the themes of the notion of a real isolated singularity and of harmonic normality, in particular for germs with several branches.

This work used AI assistance from Claude Opus 5.5 and GPT-6.2, invoked dynamically during the investigation. The author takes full responsibility for the paper.

\end{document}